\documentclass[10pt,a4paper,leqno]{amsart}

\usepackage{amsmath,amssymb,amsthm,euscript,epsfig}
\usepackage{a4wide}
\def\esp{\vspace{5pt}}

\def\C{{\mathbb C}}
\def\R{{\mathbb R}}
\def\N{{\mathbb N}}

\def\F{{\mathcal F}}
\def\S{{\EuScript S}}

\def\a{\alpha}
\def\H{{\tt H}}
\def\v{{\tt v}}
\def\V{{\tt V}}
\def\ta{{\tt a}}
\def\tA{{\tt A}}
\def\D{{\tt D}}

\def\virgp{\raise 2pt\hbox{,}}

\def\<{\langle}
\def\>{\rangle}

\def\Op{\operatorname{Op}}

\def\({\left(}
\def\){\right)}
\def\supp{\operatorname{supp}}
\def\O{\mathcal{O}}

\def\Tend#1#2{\mathop{\longrightarrow}\limits_{#1\rightarrow#2}}

\def\d{{\partial}}
\def\e{\varepsilon}

\def\om{{\omega}}

\theoremstyle{plain}
\newtheorem{theo}{Theorem}[section]
\newtheorem{lem}[theo]{Lemma}

\newtheorem{prop}[theo]{Proposition}

\theoremstyle{definition}
\newtheorem{defin}[theo]{Definition}

\theoremstyle{remark}
\newtheorem{rema}[theo]{Remark}

\numberwithin{equation}{section}

\def\CI{{C_\infty}}

\begin{document}

\title[Schr\"odinger equation with
repulsive potential]{Scattering theory for the Schr\"odinger
equation with repulsive potential} 
\author[J.-F. Bony]{Jean-Fran\c cois Bony}
\author[R. Carles]{R{\'e}mi Carles}
\author[D. H\"afner]{Dietrich H\"afner}
\author[L. Michel]{Laurent Michel}
\address[J.-F. Bony and D. H\"afner]{MAB, UMR CNRS 5466\\
Universit{\'e} Bordeaux 1\\ 351 cours de la Lib{\'e}ration\\ 33~405 Talence
cedex\\ France}
\email[J.-F. Bony]{bony@math.u-bordeaux.fr}
\email[D. H\"afner]{hafner@math.u-bordeaux.fr}
\address[R. Carles\footnote{On leave from MAB, Universit\'e Bordeaux
1.}]{IRMAR\\ 
Universit\'e de Rennes 1\\ Campus de 
Beaulieu\\ 35042 Rennes cedex\\
France} 
\email[R. Carles]{remi.carles@math.univ-rennes1.fr}
\address[L. Michel]{Institut Gallil\'ee, D\'epartement de
Math\'ematiques, Universit\'e Paris XIII, 99 avenue J.-B. Cl\'ement, 
93430 Villetaneuse, France}
\email[L. Michel]{michel@math.univ-paris13.fr}
\begin{abstract}
We consider the scattering theory for the Schr\"odinger equation with
$-\Delta -|x|^{\alpha}$ as a reference Hamiltonian, for $0< \alpha
\leq 2$, in any space dimension. We prove that when this Hamiltonian is
perturbed by a potential, the usual short  range/long range condition
is weakened: the limiting decay for the potential depends on the value
of $\a$, and is  
related to the growth of classical trajectories in the unperturbed
case. The existence of wave operators and their asymptotic
completeness are established thanks to Mourre estimates relying on new
conjugate operators. We construct the asymptotic velocity and
describe its spectrum. Some 
results are generalized to the case where 
$-|x|^{\alpha}$ is replaced by 
a general second order polynomial. 
\end{abstract}
\maketitle


\section{Introduction}

The aim of this paper is to study the scattering theory for a
large class of Hamiltonians with repulsive potential. We
find optimal short range conditions for the perturbation, and prove
asymptotic completeness under these conditions. The family of
Hamiltonians is given by:
\begin{equation}\label{eq:def}
\H_{\alpha ,0} = -\Delta - |x|^{\alpha} \quad 0<\alpha \le 2\quad ; \quad
\H_{\alpha} = \H_{\alpha ,0} +V_{\alpha} (x)\quad ; \quad x\in \R^n,\
n\geq 1. 
\end{equation}
The main new feature with respect to the usual free Schr\"odinger
operator $\H_{0,0}=-\Delta$ is the acceleration due to the potential
$-|x|^{\alpha}$. The case $\alpha=2$ is a borderline case: if
$\alpha>2$ classical trajectories reach infinite speed and $(
\H_{\alpha ,0} , C_0^{\infty}(\R) )$ is not essentially self-adjoint
(see \cite{Dunford}).

The consequence of the acceleration is that the usual position
variable increases faster than $t$ along the
evolution. Roughly speaking, the usual short range condition is:  
\begin{equation}  \label{SRF} 
|V_0(x)|\lesssim \<x\>^{-1-\varepsilon},
\end{equation}
for some $\varepsilon>0$, where $\<x\>=(1+|x|^2)^{1/2}$. One expects
 it to be weakened in the case of  $\H_{\alpha}$. 

For the Stark Hamiltonian, 
associated to a constant electric field $E\in \R^n$ (see \cite{Cycon}),
$-\Delta + E\cdot x$,
it is well known that the short range condition \eqref{SRF} becomes
$|V_s(x)|\lesssim \< E\cdot x\>^{-1/2-\varepsilon}$. We refer
to  the papers by J.~E.~Avron and 
I.~W.~Herbst \cite{AvronHerbst,Herbst} for weaker conditions. The idea
is that the drift caused by $E$ (which may also model gravity, see
e.g. \cite{TW}) accelerates the particles in the direction of the
electric field. 
This phenomenon has been observed for a larger class of Hamiltonians 
 by M.~Ben-Artzi \cite{Be83,Be84}: generalizing the Stark Hamiltonian
($\a =1$), let
\begin{equation*}
\widehat{\H}_{\alpha, 0} = -\Delta - {\rm sgn}(x_{1}) |x_{1}|^{\alpha}\,,
 \quad 0<\alpha \le 2; \quad 
\widehat{\H}_{\alpha} = \widehat{\H}_{\alpha ,0} +
\widehat{V}_{\alpha} (x)\,, 
\end{equation*}
with $x=(x_{1} , x')$. In \cite{Be84}, asymptotic completeness is
proved under the condition 
\begin{equation}  \label{SRS}
| \widehat{V}_{\alpha} | \lesssim M\(x'\) \cdot \left\{
\begin{aligned}
\<& x_{1} \>^{\alpha - \e} \quad &\text{for } x_{1} \leq 0,   \\
\<& x_{1} \>^{-1 + \alpha/2 - \e} \quad &\text{for } x_{1} \geq 0,
\end{aligned}
\right.
\end{equation}
with $\e >0$ and $M(x') \to 0$ as $|x'| \to \infty$. In the one
dimensional case, we obtain similar results for $0< \alpha <2$, and a
weaker condition for $ \alpha =2$. The proofs in \cite{Be83,Be84} rely on 
some specific properties of one-dimensional Hamiltonians, and it seems
they cannot be adapted to \eqref{eq:def} when $n\geq 2$. Our
approach is 
completely different, since it is based on Mourre estimates.

Notice that for $\alpha \in ]0,2]$, the Hamiltonian
$\H_{\a ,0}$ shares an interesting difficulty with the Stark Hamiltonian: its
symbol, $|\xi|^2 -|x|^\a$, is not signed, and can take arbitrarily
large negative values.

The case $\alpha=2$ is in some sense very instructive. A nonlinear
scattering theory is already available in this case. In \cite{CaSIMA},
the second author studied nonlinear perturbations of $\H_{2,0}$ and
showed that all the usual nonlinearities are short range. This is
closely related to the fact that the classical trajectories can be
computed explicitly:  
\begin{equation*}
x(t)=\frac{1}{2}(x_0 + \xi_0)e^{2t}+ \frac{1}{2}(x_0-\xi_0)e^{-2t}.
\end{equation*}
Thus $x(t)$ grows exponentially fast (in general). For $0< \alpha <2$, a formal
computation shows that the classical trajectory $x(t)$ can go to
infinity like $t^{1/(1-\alpha /2)}$;  denoting ${\tt
U}_\a(x)=-|x|^\a$, the equations of motion imply:
\begin{equation*}
0=\ddot x(t) + \nabla {\tt U}_\a\( x(t)\) = \ddot x(t) - \a
|x(t)|^{\a -2}x(t)\, .
\end{equation*}
Seeking a particular solution of the form $x(t) = t^\kappa y$, for a
constant $y\in \R^n$, yields $\kappa -2 = (\a -1) \kappa$, hence $
\kappa = 2/(2-\alpha)$. We will prove that in general, $x(t)$ does
go to infinity like $t^{1/(1-\alpha /2)}$. This shows that the
acceleration causes by $-|x|^\a$ increases progressively as $\a$
ranges $]0,2]$. For a small $\alpha>0$, the particle moves hardly faster than in the free case
$|x(t)|=\O(t)$. As $\a$ increases, the particle goes to infinity
faster and faster, and reaches the maximal exponential growth for
$\a=2$. For $\a>2$, it is known that particles can reach an infinite
speed, which is the reason why $(
\H_{\alpha ,0} , C_0^{\infty}(\R) )$ is not essentially self-adjoint. 
This suggests to define as a new position
variable,
\begin{equation}  \label{tx34}
p_{\alpha} (x) = \left\{
\begin{aligned}
& \ln \< x \> &&\text{for } \alpha =2,  \\
& \< x \>^{1-\alpha/2} &&\text{for } 0< \alpha <2.
\end{aligned}
\right.
\end{equation}
We assume that the multiplication potential $V_{\alpha} (x)$ is real-valued, and writes
$V_{\alpha} (x) = V_{\alpha}^{1} (x) + V_{\alpha}^{2} (x)$,
with:
\begin{equation}  \label{va1}
V_{\alpha}^{1} \text{ is a measurable real-valued function, compactly
supported, and }\Delta \text{--compact}, 
\end{equation}
and $V_{\alpha}^2 \in L^{\infty} (\R^{n};\R)$ satisfies the \emph{short range
condition}:
\begin{equation}  \label{va2}
|V_{\alpha}^{2} (x) | \lesssim p_{\alpha} (x)^{-1-\varepsilon},
 \quad\text{a.e. }x\in\R^n\, ,
\end{equation}
for some  $\varepsilon >0$.

The operator $\H_{\alpha}$ is essentially self-adjoint, with domain the
domain of the harmonic oscillator, and we denote again $\H_{\alpha}$ its
self-adjoint extension. In Section~\ref{ele}, we prove that
$\H_{\alpha}$ has no singular continuous 
spectrum and $\sigma ( \H_{\alpha} ) = \R$. Under general assumptions
(see Theorem~\ref{glm}), we also show that its point spectrum is
empty. We can now state 
the main result of this paper.

\begin{theo}[Asymptotic completeness]    \label{comple}
Let $0<\a\leq 2$, and $\H_{\a ,0}$, $\H_\a$ defined by \eqref{eq:def}. 
Assume that $V_\a=V^1_\a+V^2_\a$ satisfies  \eqref{va1} and
\eqref{va2}. Then the following limits exist:  
\begin{gather}
s-\lim_{t\rightarrow\infty}e^{it\H_{\alpha}}e^{-it\H_{\alpha ,0}},
\label{W1} \\ 
s-\lim_{t\rightarrow\infty} e^{it\H_{\alpha ,0}} e^{-it\H_{\alpha}} {\bf
1}^c (\H_{\alpha}),   \label{W2} 
\end{gather}
where ${\bf 1}^c (\H_{\alpha})$ is the projection on the continuous
spectrum of $\H_{\alpha}$. If we denote \eqref{W1} by $\Omega^+$, then
\eqref{W2} is equal to $(\Omega^+)^*$, and we have: 
\begin{equation}
(\Omega^+)^*\Omega^+={\bf 1} \quad\text{and}\quad \Omega^+(\Omega^+)^*={\bf
1}^c(\H_{\alpha}). 
\end{equation}
\end{theo}
To prove this result, we establish Mourre estimates, relying on new
conjugate operators, adapted to the repulsive potential
$-|x|^\a$  (as
a matter of fact, we work with the smoother repulsive potential
$-\<x\>^\a$; see below). To give an idea of the difficulty at 
this stage, consider the one-dimensional case. For $\a=2$,  the
natural idea for a conjugate operator is to consider the  generator of
dilations $(xD+Dx) /2$: 
\begin{equation*}
i \big[ \H_{2 ,0} , (x D + D x)/2 \big] = -2 \Delta + 2 x^{2} \, . 
\end{equation*}
This is the harmonic oscillator, which is of course positive. This
seems an encouraging point. Nevertheless, it is not
$\H_{2,0}$--bounded; we must find another conjugate
operator. Therefore, we look for a 
pseudo-differential operator $A_{2}$ with symbol $a_{2} (x,\xi)$,  and
try to solve:  
\begin{equation*}
\left\{ \xi^{2} - x^{2} , a_{2} (x,\xi) \right\} = 4,\quad \text{on the
energy level } \left\{(x,\xi)\ 
;\ \xi^{2} -  x^{2} = E\right\}.  
\end{equation*}
A solution to this equation is given by:
$\displaystyle a_{2} (x,\xi) =  \ln (\xi +x) - \ln ( \xi -x)$. 
Now consider the case $0 < \alpha < 2$. For $x>0$, we
try to solve  
\begin{equation*}
\left\{ \xi^{2} - x^{\alpha} ,a_{\a} (x,\xi) \right\}  = 2-\a,\quad
 \text{on } \left\{(x,\xi)\in \R^*_+\times\R\ ;\ \xi^{2} - x^{\alpha}
 = E\right\}.  
\end{equation*}
Plugging $a_{\a} (x,\xi) = \xi x^{1-\alpha} $
into this equation, we get
\begin{equation*}
\{ \xi^{2} - x^{\alpha} , a_{\a} (x,\xi) \} = 2-\a + 2E
(1-\alpha) x^{-\alpha}\, , \quad\text{for } \xi^{2} - x^{\alpha} =
E.
\end{equation*}
The term in $x^{-\alpha}$ should not matter for the Mourre estimate,
since it is compact on the energy level. This formal
discussion is the foundation for the constructions of Section~\ref{iii}.

To apply Mourre's method, truncations in energy are needed, of the form
$\chi(-\Delta -|x|^\a)$. However, $\H_{\a,0}$ is not elliptic,
so it is not clear that this defines a good
pseudo-differential operator. These difficulties are solved in
Section~\ref{sec:mourre}, where we consider the general case.
\esp

Theorem~\ref{comple} shows how the short range condition \eqref{va2} takes 
the acceleration caused by the potential into account. 
Let us note that the condition \eqref{va2} is not necessary the weakest one, but the decay $|V(x)|\lesssim
p_{\alpha}(x)^{-1}$ at infinity is expected to be the borderline case
between long range and short range scattering, because the position
variable increases exactly like $t$ along the evolution (compare with
Theorem~\ref{Th.P.1}). For the case of the Stark Hamiltonian, it is well-known
that the case $\varepsilon =0$ in \eqref{SRS} is the limiting case,
which involves long range effects (see \cite{OzawaStark2}).

We obtain more precise informations by constructing the asymptotic
velocity (see e.g. \cite{DG}). We note $\CI (\R^{n})$ the set of continuous functions which go to $0$ at $\infty$. Let $B^m:=(B_1^m,...,B_n^m)$ be a
sequence of commuting self-adjoint 
operators on a Hilbert space ${\mathcal H}$. Suppose that for every
$g\in \CI(\R^n)$, there exists 
\begin{equation}
\label{P.0.1}
s-\lim_{m\rightarrow\infty}g(B^m).
\end{equation}
Then by \cite[Proposition B.2.1]{DG}, there exists a unique vector
$B=(B_1,...,B_n)$ of 
commuting self-adjoint operators such that
\eqref{P.0.1} equals $g(B)$. $B$ is densely defined if, for some $g\in
\CI(\R^n)$ such that $g(0)=1$, 
\begin{equation*}
s-\lim_{R\rightarrow\infty}(s-\lim_{m\rightarrow\infty}(g(R^{-1}B_m))={\bf
1}. 
\end{equation*}
We denote $\displaystyle B:=s-\CI-\lim_{m\rightarrow\infty}B^m$.

\begin{theo}[Asymptotic velocity]
\label{Th.P.1}
Let $\sigma_\a$ be given by:
\begin{equation*}
\sigma_{\alpha} = \left\{
\begin{aligned}
2 & - \alpha  \quad && \text{if } 0< \alpha <2 \,, \\
2 & \quad && \text{if } \alpha =2\,.
\end{aligned}
\right.
\end{equation*}
There exists a bounded self-adjoint operator $P_{\alpha}^+$, which
commutes with $\H_{\alpha}$, such that: 
\begin{itemize}
\item[(i)] $\displaystyle
P_{\alpha}^+=s-\CI-\lim_{t\rightarrow\infty}
e^{it \H_{\alpha}}\frac{p_{\alpha}(x)}{t}e^{-it \H_{\alpha}}$.
\item[(ii)] The spectrum of the asymptotic velocity is $\displaystyle
\sigma(P_{\alpha}^+)= 
\left\{
\begin{aligned}
\{ & 0,\sigma_{\alpha}\} \quad && \text{if } \sigma_{pp} ( \H_{\alpha})\neq\emptyset\, ,\\ 
\{ & \sigma_{\alpha}\} \quad && \text{if } \sigma_{pp} ( \H_{\alpha})=\emptyset\, . 
\end{aligned}
\right. $
\item[(iii)] $\displaystyle {\bf 1}_{\{0\}}(P_{\alpha}^+)={\bf
1}^{pp}( \H_{\alpha})$.
\esp

\item[(iv)] For any $J\in \CI(\R)$, we have
\begin{align*}
J(P_{\alpha}^+){\bf
1}_{\R\setminus\{0\}}(P_{\alpha}^+)= s-\lim_{t\rightarrow\infty}
e^{it \H_{\alpha}}J( \V_{\alpha} )e^{-it \H_{\alpha}}{\bf 
1}_{\R\setminus\{0\}}(P_{\alpha}^+),
\end{align*}
where $\V_{\alpha}:=\left[iH_{\alpha},p_{\alpha}\(x\)\right]$ is the local velocity.
\end{itemize}
\end{theo}

\begin{rema}
Using the decomposition ${\mathcal H} = {\mathcal H}_{\rm pp} (\H_{\alpha})\oplus{\mathcal H}_c (\H_{\alpha})$ we can thus write
$P_{\alpha}^+=0\oplus\sigma_{\alpha}$. 
\end{rema}

Let us note that the limits we stated in Theorem~\ref{comple} and
Theorem~\ref{Th.P.1} are for $t\to + \infty$;  analogous results
obviously hold for $t\to -\infty$.  
\esp

Notice that computing the asymptotic velocity is all the more
interesting that the free dynamics, $e^{-it\H_{\a ,0}}$, is not
known in the case $0<\a<2$. On the other hand, it is very well
understood in the case $\a=2$, since a generalized Mehler's formula is
available (see \cite{HormanderQuad} and Section~\ref{sec:mehler} below). For
$\a=2$, we also consider more general Hamiltonians:
\begin{equation}\label{eq:H0}
H_0 = -\Delta - \sum_{k=1}^{n_{-}} \omega_{k}^{2} x_{k}^{2} +
\sum_{k=n_{-} + 1}^{n_{-} + n_{+}} \omega_{k}^{2} x_{k}^{2} +
\sum_{k=n_{-} + n_{+} +1}^{n_{-} + n_{+} + n_{E}} E_{k} x_{k}\quad ; \quad H
= H_{0} + V(x) \, ,
\end{equation}
with $\omega_{j} >0$ and $E_{j} \neq 0$. We prove the existence of wave operators in this more general case,
under weaker conditions than \eqref{va2} (see
Section~\ref{sec:mehlercook}). The proof is based on an explicit
formula for the dynamics $e^{-itH_{0}}$ (Mehler's formula, see
Section~\ref{sec:mehler}).  

Asymptotic completeness is shown if $n_{-} + n_{+} =n$, under a condition
similar to \eqref{va2} in Section~\ref{gena2}. In that case, we also
construct asymptotic velocities in each space direction. The
asymptotic velocity, given by Theorem~\ref{Th.P.1},
also exists and is equal to  
$P_\ell^+$, where $\omega_\ell = \max_{1\leq j \leq n_-}
\omega_j$,  and $P_\ell^+$ is the asymptotic velocity in the direction
$x_\ell$.

To our knowledge, there is very little motivation from a physical
point of view to study the above Hamiltonians:  in general,
electromagnetic fields have saddle points, like the potential in 
$H_0$, but the above model should then be valid only
\emph{locally}, in a neighborhood of the saddle point. In
\eqref{eq:def}, the potentials
$-|x|^\a$ are unbounded from below; this does not seem physically
relevant (notice however that the Stark potential $E\cdot x$ is also
unbounded). 

On the other
hand, we believe that these models are mathematically interesting. The
dependence on $\a\in ]0,2]$ is somehow well understood, in particular
thanks to the definition of the position variable $p_\a$ \eqref{tx34}
and to the study of the asymptotic velocity. We also introduce new
conjugate operators in order to obtain Mourre estimates (see
\eqref{nnn} and
\eqref{nji}). Here again, the dependence of the analysis upon $\a$
seems to be interesting (in particular the limiting case $\a=2$ is
better understood than in \cite{Be84}).  
\esp

As mentioned above, in our analysis, we replace $\H_{\alpha ,0}$ and
$\H_{\alpha}$ by 
\begin{equation}\label{eq:vraiH}
H_{\alpha ,0} = -\Delta - \< x \>^{\alpha} \quad\text{and}\quad  
H_{\alpha} = H_{\alpha ,0} +V_{\alpha} (x).
\end{equation}
This does not affect the results, since for large $|x|$, $\left| \< x
\>^{\alpha} - |x|^{\alpha} \right|$ is estimated by $\< x \>^{\alpha
-2}$, which is a short range perturbation for $0< \alpha \leq 2$ (no
 smoothness is required for the perturbative potentials). We
 therefore prove Theorems~\ref{comple} and \ref{Th.P.1}  with
 $\H_{\a,0}$ (resp. $\H_\a$)  replaced by $H_{\a,0}$ (resp. $H_\a$).

The paper is organized as follows.  
\begin{itemize}
\item In Section~\ref{ele} we show some elementary properties of the
Hamiltonians. In particular, we recall  Mehler's formula for $\alpha
=2$, and prove the absence of eigenvalues for $H_\a$ in many cases.  

\item Section~\ref{sec:mourre} is devoted to the Mourre estimate. In
Section~\ref{ddd}, we treat some rather technical features. For
example, $\chi(\xi^2-\<x\>^{\alpha})$ is not a good symbol, and we need
some preparations before being able to use the pseudo-differential
calculus (see Proposition~\ref{pse}). We give the conjugate operator $A_\a$,
which is a pseudo-differential operator, 
in Section~\ref{iii}. Section~\ref{qqq} is devoted to the regularity
results and the Mourre estimate 
is established in Section~\ref{moumou}.  

\item In Section~\ref{sec:complet} we prove asymptotic
completeness. The Mourre estimate yields a minimal velocity estimate
for $A_{\alpha}$. We obtain a minimal velocity
estimate for the observable $p_{\alpha} (x)$ using a lemma due to
C.~G\'erard and F.~Nier (see \cite{GeNi98}). 

\item In Section~\ref{sP}, we construct the asymptotic velocity and
describe its spectrum. 

\item In Section~\ref{www}, we generalize our results in the case of
the Hamiltonians defined in \eqref{eq:H0}.
\end{itemize}

\noindent {\bf Acknowledgments}. The authors are grateful to Vincent
Bruneau and Thierry Jecko for stimulating discussions on this work.

\section{Elementary properties}
\label{ele}

\subsection{Domain and spectrum}  \label{doma}

We begin with some properties on the spectrum of the operator
$H_{\alpha}$. For $0< \alpha \leq 2$, introduce $N_{\alpha} = - \Delta + \<
x \>^{\alpha}$. It is self-adjoint, with domain  
\begin{equation*} 
D(N_{\alpha} ) = \left\{ u \in H^{2} (\R^{n}); \ \< x \>^{\alpha} u \in
L^{2} (\R^{n}) \right\}.  
\end{equation*}
It can be viewed as the ``confining'' counterpart of $H_{\alpha,0}$
(the repulsive potential $-\<x\>^\a$ is replaced by the confining one
$+\<x\>^\a$). 
Since it is not easy to know the domains of $H_{\alpha,0}$, we work on
a core for these operators, $D(N_2)$, the domain of the
harmonic oscillator. We recall an extension of 
Nelson's Theorem due to C.~G\'erard and I.~\L aba
\cite[Lemma~1.2.5]{GeLa02}:

\begin{theo}[Nelson's Theorem] \label{nel}
Let ${\mathcal H}$  be a Hilbert space, $N\geq 1$ a self-adjoint
operator on ${\mathcal H}$, $H$ a symmetric operator such that $D(N)
\subset D(H)$, and 
\begin{gather}
\| H u \| \lesssim \| Nu \| \, , \quad u \in D(N),   \label{nel1}\\
|(Hu,Nu) - (Nu,Hu) | \lesssim \| N^{1/2} u \|^{2} \, , \quad u \in D(N).
 \label{nel2} 
\end{gather}
Then $H$ is essentially self-adjoint on $D(N)$, and we denote
$\overline{H}$ its extension. If $u \in D(\overline{H})$, then
$(1+i\varepsilon N)^{-1} u$ converges to $u$ in the graph topology of
$D(\overline{H})$ as $\varepsilon \to 0$. 
\end{theo}

From this theorem, we deduce the following:
\begin{lem}   \label{lem:essadj0}
For any $\a \in ]0,2]$, the operator $H_{\alpha,0}$ is essentially
self-adjoint on $D(N_2)$. 
\end{lem}

\begin{proof}
For $u\in D(N_2)$, we have
\begin{equation*}
\| H_{\alpha,0} u \| \leq \| \<x\>^{\alpha} u \| + \| -\Delta u \|
\lesssim\|N_2u\| \, ,
\end{equation*}
which proves \eqref{nel1}. Now, let us prove (\ref{nel2}). A straightforward
computation shows that  
\begin{equation*}
[H_{\alpha,0},N_{2}]=\left[\<x\>^2+\<x\>^{\alpha},\Delta\right].
\end{equation*}
Hence, it suffices to show that
$\|[\<x\>^{\alpha},\Delta]\|\lesssim\|N_{2}^{1/2}u\|$ for
$0<\alpha\leq 2$. But we have
\begin{equation*}
[\<x\>^{\alpha},\Delta]=-2i\alpha\<x\>^{\alpha-2}xD -
n\alpha\<x\>^{\alpha-2} - \alpha(\alpha-2)\<x\>^{\alpha-4} x^2\, ,
\end{equation*}
which is clearly bounded by $N_2^{1/2}$ for $0<\alpha\leq 2$.
\end{proof}

Before going further, let us notice the following characterization of
$H_{\alpha}$--compactness:

\begin{lem}  \label{lem:compH}
Let $V_{\alpha} (x) = V_{\alpha}^{1} (x) + V_{\alpha}^{2} (x)$, where
$V_{\alpha}^{1}$ is a compactly supported 
measurable function, and $V_{\alpha}^{2} \in
L^{\infty} (\R^{n})$ with $V_{\alpha}^{2} (x) \to 0$ as $|x| \to
\infty$. Then, $V_{\alpha}$ is $H_{\alpha,0}$--compact if and only if
$V_{\alpha}^{1}$ is $\Delta$--compact. 
\end{lem}

\begin{proof}
Let $\chi_{1}$, $\chi_{2} \in C^{\infty}_{0} (\R^{n})$ be such that $\chi_{2} = 1$ near the support of $\chi_{1}$. We have
\begin{align}
\chi_{1} (x) (H_{\alpha,0} +i)^{-1} =& (\Delta -i)^{-1} (\Delta -i)
\chi_{1} (H_{\alpha,0} +i)^{-1}  \nonumber    \\ 
=& (\Delta -i)^{-1} \chi_{1} (\Delta -i) (H_{\alpha,0} +i)^{-1} +
(\Delta -i)^{-1} [ \Delta , \chi_{1} ] (H_{\alpha,0} +i)^{-1}
\nonumber    \\ 
=& - (\Delta -i)^{-1} \chi_{1} - (\Delta -i)^{-1} \chi_{1} \< x
\>^{\alpha} (H_{\alpha},0 +i)^{-1}    \nonumber   \\ 
&+ (\Delta -i)^{-1} [ \Delta , \chi_{1} ] (\Delta -i)^{-1} (\Delta -i)
\chi_{2} (H_{\alpha,0} +i)^{-1}  \nonumber   \\ 
=& (\Delta -i)^{-1}  {\mathcal O} (1) + (\Delta -i)^{-1} [ \Delta ,
\chi_{1} ] (\Delta -i)^{-1} \chi_{2} \left( -1 - \< x \>^{\alpha}
(H_{\alpha,0} +i)^{-1} \right)  \nonumber   \\ 
&+ (\Delta -i)^{-1} [ \Delta , \chi_{1} ] (\Delta -i)^{-1} [\Delta ,
\chi_{2}] (H_{\alpha,0} +i)^{-1}   \nonumber   \\ 
=& (\Delta -i)^{-1} {\mathcal O} (1),   \label{tx12}
\end{align}
since $[ \Delta , \chi_{1} ] (\Delta -i)^{-1}$ and $[ \Delta ,
\chi_{1} ] (\Delta -i)^{-1} [\Delta , \chi_{2}]$ are bounded. On the
other hand, we have 
\begin{align}
\chi_{1} (x) (\Delta -i)^{-1} =& (H_{\alpha,0} +i)^{-1} (H_{\alpha,0} +i)
\chi_{1} (\Delta -i)^{-1}   \nonumber   \\ 
=& (H_{\alpha,0} +i)^{-1} \chi_{1} (x) (H_{\alpha,0} +i) (\Delta -i)^{-1}
- (H_{\alpha,0} +i)^{-1} [\Delta , \chi_{1}] (\Delta -i)^{-1}
\nonumber  \\ 
=& (H_{\alpha,0} +i)^{-1} {\mathcal O} (1).  \label{tx13}
\end{align}
Since $V_{\alpha}^{2} (x) \to 0$ as $x \to \infty$, we get
\begin{align*}
{\bf 1}_{|x| > R} V_{\alpha} (x) (H_{\alpha,0} +i)^{-1} = {\bf 1}_{|x| >
R} V_{\alpha}^{2} (x) (H_{\alpha,0} +i)^{-1} \to 0\quad\text{as }R\to
\infty, 
\end{align*}
and ${\bf 1}_{|x| > R} V_{\alpha} (x) (\Delta -i)^{-1} \to 0$ for the
norm topology as $R \to + \infty$. So, from \eqref{tx12} and
\eqref{tx13}, 
\begin{equation*}
V_{\alpha} \text{ is }H_{\alpha,0}\text{--compact } \Leftrightarrow
{\bf 1}_{|x| <R} V_{\alpha} \text{ is }H_{\alpha,0}\text{--compact }
  \Leftrightarrow V_{\alpha}^{1}\text{ is
}H_{\alpha,0}\text{--compact } \Leftrightarrow V_{\alpha}^{1}\text{ is
}\Delta\text{--compact.} 
\end{equation*}
\end{proof}

Then Lemma \ref{lem:essadj0} and the Kato--Rellich theorem
\cite[Theorem X.12]{ReedSimon2} imply: 

\begin{lem}   \label{aut}
Let $0< \alpha \leq 2$. For $V_\a=V^1_\a+V^2_\a$ satisfying
\eqref{va1} and \eqref{va2},  the operator $H_{\alpha}$ is
self-adjoint on $D(H_{\alpha ,0})$, and essentially self-adjoint on
$D(N_2)$.  
\end{lem}

Notice that $H_{2,0}+1=-\Delta-x^2$ is conjugated to the generator of
dilations. This is obvious from a glance at the symbols: $\xi^2-x^2$
can be written as $(\xi+x)\cdot(\xi-x)= y\cdot \eta$, with suitable
new variables corresponding to a rotation of angle $\pi/4$ in the
phase space. To make this argument precise, we set, for
$u \in {\mathcal S}' (\R^{n})$:
\begin{equation*}
U u (x) = \frac{1}{(\sqrt{2} \pi )^{n/2}}   e^{-ix^{2}/2} \int e^{i
\sqrt{2} x \cdot y} e^{-iy^{2}/2} u(y) \, dy. 
\end{equation*}
The operator $U$ is an isometry on $L^{2} (\R^{n})$. We have
\begin{equation}\label{pmo}
x U u  = U \(\frac{y - D_{y}}{\sqrt{2}} u\) ,  \quad \text{and}\quad 
D_{x} U u (x) = U \(\frac{D_{y} + y}{\sqrt{2}} u\) .  
\end{equation}
Using these relations, it is easy to see that:
\begin{equation} \label{mmm}
N_2U = UN_2 \text{ and } 
H_{2,0} U = U \widetilde{H}_{2,0}\, ,
\end{equation}
where
\begin{equation}   \label{eq:defN_Htilde}
\widetilde{H}_{2,0} = D_x x + x D_x-1.
\end{equation}
Then, from \cite[Proposition~6.2]{Pe83} on the spectrum of $Dx+xD$
and the Weyl's essential spectrum theorem
\cite[Theorem~XIII.14]{ReedSimon4}, we obtain:

\begin{prop}  \label{sp2}
The spectrum of $H_{2,0}$ is purely absolutely continuous, and $\sigma
(H_{2}) = \R$ if $V_{2}$ is an $H_{2,0}$--compact real-valued
potential. 
\end{prop}

For $0< \alpha <2$, we have the following proposition.

\begin{prop}  \label{spa}
Let $V_{\alpha}$ be a $H_{\alpha ,0}$--compact potential
with $0< \alpha < 2$. Then 
\begin{equation}
\sigma (H_{\alpha}) = \R .
\end{equation}
\end{prop}

\begin{proof}
It is enough to show that $\sigma (H_{\alpha ,0}) = \R$. In that case,
we have $\sigma_{\rm ess} (H_{\alpha ,0}) = \R$ and then, by the Weyl's
theorem \cite[Theorem XIII.14]{ReedSimon4}, $\sigma_{\rm ess} (H_{\alpha
}) = \R$. Since $H_{\alpha}$ is self-adjoint, we get the proposition.

Let $\varphi \in C_{0}^{\infty} (]0, + \infty[; [0,1])$ so that
$\varphi =1$ near $1$. For $E \in \R$, we set 
\begin{gather*}
u (x_{1}) = e^{i x_{1}^{1+\alpha /2} /(1+\alpha /2)} e^{i E
x_{1}^{1-\alpha /2} /(2-\alpha )} ,  \\ 
u_{\varepsilon, \delta} (x) = u (x_{1}) \sqrt{\varepsilon} \varphi
(\varepsilon x_{1}) \delta^{(n-1)/2} \varphi (\delta | x'|), 
\end{gather*}
where $x= (x_{1}, x')$. We first note that $\| u_{\varepsilon, \delta}
\|_{L^{2} (\R^{n})}$ does not depend on $\varepsilon$ and $\delta$. We
have 
\begin{align*}
\Delta u_{\varepsilon, \delta} =& \partial^{2}_{x_{1}} (u)
\sqrt{\varepsilon} \varphi (\varepsilon x_{1}) \delta^{(n-1)/2}
\varphi (\delta | x'|) + 2 \partial_{x_{1}} (u) \partial_{x_{1}} \(
\sqrt{\varepsilon} \varphi (\varepsilon x_{1}) \) \delta^{(n-1)/2}
\varphi (\delta | x'|) \\ 
&+ u \partial^{2}_{x_{1}} \( \sqrt{\varepsilon} \varphi (\varepsilon
x_{1}) \) \delta^{(n-1)/2} \varphi (\delta | x'|) + u
\sqrt{\varepsilon} \varphi (\varepsilon x_{1}) \partial^{2}_{x'} \(
\delta^{(n-1)/2} \varphi (\delta | x'|) \). 
\end{align*}
The second term is equal to ${\mathcal O} (x_{1}^{\alpha /2})
\varepsilon^{3/2} |\varphi '| (\varepsilon x_{1}) \delta^{(n-1)/2}
\varphi (\delta | x'|)$ which is ${\mathcal O} (\varepsilon^{1-\alpha
/2})$ in $L^{2}$-norm.  
The third and fourth terms are ${\mathcal O} (\varepsilon^{2})$ and
${\mathcal O} (\delta^{2})$ respectively in $L^{2}$ norm. Then 
\begin{equation*}
H_{\alpha ,0} u_{\varepsilon, \delta} = \sqrt{\varepsilon} \varphi
(\varepsilon x_{1}) \delta^{(n-1)/2} \varphi (\delta | x'|) H_{\alpha
,0} \( u \) + o(1), 
\end{equation*}
as $\varepsilon, \delta \to 0$.
Since
\begin{align*}
\< x \>^{\alpha} - \< x_{1} \>^{\alpha} = {\mathcal O} (1) \< x_{1}
\>^{\alpha -2} \< x' \>^{2}\, ,
\end{align*}
we get
\begin{equation*}
\sqrt{\varepsilon} \varphi (\varepsilon x_{1}) \delta^{(n-1)/2}
\varphi (\delta | x'|) \big( \< x \>^{\alpha} - \< x_{1} \>^{\alpha}
\big)u = {\mathcal O} \(\varepsilon^{2-\alpha} \delta^{-2}\)\, . 
\end{equation*}
If $\varepsilon^{2-\alpha} \delta^{-2} \to 0$, we have
\begin{align*}
H_{\alpha ,0} u_{\varepsilon, \delta} = \sqrt{\varepsilon} \varphi
(\varepsilon x_{1}) \delta^{(n-1)/2} \varphi (\delta | x'|)  \(
-\partial_{x_{1}}^{2} - \< x_{1} \>^{\alpha} \) u( x_{1}) + o(1). 
\end{align*}
But we have
\begin{align*}
\partial_{x_{1}}^{2} u(x_{1}) = \( -x_{1}^{\alpha} -E -E^{2}
x_{1}^{-\alpha}/4 +i \alpha x_{1}^{\alpha /2-1} /2 - E \alpha
x_{1}^{-\alpha /2 -1} /4 \)u (x_{1}), 
\end{align*}
and then, there is a $\mu >0$ so that
\begin{align*}
(H_{\alpha ,0} -E) u_{\varepsilon, \delta} = \sqrt{\varepsilon}
\varphi (\varepsilon x_{1}) \delta^{(n-1)/2} \varphi (\delta | x'|)
{\mathcal O} (x_{1}^{- \mu}) + o(1) = {\mathcal O} (\varepsilon^{\mu})
+ o(1) = o(1). 
\end{align*}
By the Weyl's criterion \cite[Theorem VII.12]{ReedSimon1}, $E$ is in
$\sigma (H_{\alpha ,0})$. 
\end{proof}

\subsection{Generalized Mehler's formula}
\label{sec:mehler}

In this paragraph, we restrict our attention to the case $\alpha =2$
and drop the index $2$. We consider a more general Hamiltonian on
$L^{2} (\R^{n})$, 
\begin{equation}  \label{tx21}
H_0 = -\Delta - \sum_{k=1}^{n_{-}} \omega_{k}^{2} x_{k}^{2} +
\sum_{k=n_{-} + 1}^{n_{-} + n_{+}} \omega_{k}^{2} x_{k}^{2} +
\sum_{k=n_{-} + n_{+} +1}^{n_{-} + n_{+} + n_{E}} E_{k} x_{k}, 
\end{equation}
with $n_{-} + n_{+} + n_{E} \leq n$, $\omega_{k} >0$ and $E_{k} \neq
0$ if $n_E\not =0$. By convention, $\sum_{j=a}^{b} =0$ if $b<a$.

In this case, ${H}_0$ is essentially self-adjoint on
$C_0^\infty(\R^n)$ from Faris--Lavine Theorem
(\cite[p.~199]{ReedSimon2}). The kernel of $e^{-it{H}_0}$ is known
explicitly, 
through a generalized Mehler's formula (see
e.g. \cite{HormanderQuad}):  
\begin{equation}\label{eq:mehlergen}
e^{-it{H}_0}f=
\prod_{k=1}^n
\left(\frac{1}{2i\pi g_k(2t)}\right)^{1/2}
\int_{\R^n}e^{iS(t,x,y)}f(y)dy\, ,
\end{equation}
where
\begin{equation*}
S(t,x,y)= \sum_{k=1}^n \frac{1}{g_k(2t)}\(
\frac{x_k^2 +y_k^2}{2}h_k(2t) -x_k
y_k\)- \sum_{k=n_-+n_++1}^{n_-+n_++n_{E}} \( \frac{E_k}{2} (x_{k} + y_{k}) t
+\frac{E_k^2}{12} t^3 \),
\end{equation*}
and the functions $g_k$ and $h_k$, related to the classical
trajectories, are given by:
\begin{equation}\label{eq:rays} 
\begin{aligned}
g_k(t)=&\left\{
\begin{aligned}
\frac{\sinh (\om_k t)}{\om_k}\, ,\ &\textrm{ for }1\leq k\leq n_-\, ,\\
\frac{\sin (\om_k t)}{\om_k}\, ,\ &\textrm{ for } n_-+1\leq k\leq
n_-+n_+\, ,\\ 
t\ ,\ &\textrm{ for }k>n_-+n_+ \,, 
\end{aligned}
\right.\\
\, \\
h_k(t)=&\left\{
\begin{aligned}
\cosh (\om_j t)\, ,\ &\textrm{ for }1\leq k\leq n_-\, ,\\
\cos(\om_k t)\, ,\ &\textrm{ for } n_-+1\leq k\leq
n_-+n_+\, ,\\
1\ ,\ &\textrm{ for }k>n_-+n_+ \, .
\end{aligned}
\right.
\end{aligned}
\end{equation}
Recall that if $n_+\geq 1$, then $e^{-it{H}_0}$ has some singularities
(see e.g. \cite{KR}). This affects the 
above formula with phase factors we did not write (which can be
incorporated in the definition of $(ig_k(2t))^{1/2}$), but not the
computations we shall make in Section~\ref{sec:mehlercook}.

The group generated by $H_0$ is given by Mehler's formula
\eqref{eq:mehlergen}, and can
be factored in an agreeable way, in the same fashion as
$e^{it\Delta}$ (see for instance
\cite{Ozawa,JensenOzawa,Ginibre}). Recalling \eqref{eq:mehlergen} and
\eqref{eq:rays}, we have
\begin{equation}  \label{eq:factor}
e^{ -it H_{0} } = {\mathcal M}_t {\mathcal D}_t {\mathcal F} {\mathcal M}_t 
e^{-i \frac{t^3}{12} |E|^2} \, ,
\end{equation}
where $E = (E_{n_{-} + n_{+} +1} ,  \ldots , E_{n_{-} + n_{+} +n_{E}} )$,
\begin{gather*}
{\mathcal M}_t={\mathcal M}_t(x)=
\exp \bigg( i \sum_{k=1}^n x_k^2\frac{h_k(2t)}{2 g_k(2t)}
-i\frac{t}{2} \sum_{k= n_{-} + n_{+} +1}^{n_{-} + n_{+} +n_{E}} E_{k} x_k \bigg) \, ,\\ 
({\mathcal D}_t\varphi)(x)= \prod_{k=1}^n
\left(\frac{1}{i g_k(2t)}\right)^{1/2}\varphi
\left(\frac{x_1}{g_1(2t)}\virgp \ldots \virgp \frac{x_n}{g_n(2t)}
\right)\, ,
\end{gather*}
and 
\begin{equation} \label{eq:Fourier}
{\mathcal F} \varphi (\xi) =\hat \varphi (\xi)=
\frac{1}{(2\pi)^{n/2}}\int_{\R^n}e^{-ix\cdot \xi}\varphi(x)dx,
\end{equation}
denotes the Fourier transform.

\subsection{Absence of eigenvalues}

We prove the absence of embedded eigenvalues under the unique
continuation property. This result is very similar to
\cite[Theorem~XIII.58]{ReedSimon4}.   We recall the notion of unique
continuation property.

\begin{defin}
A Hamiltonian $H$ has the \emph{unique continuation property} if the following holds. Suppose that $Hu=0$ for some $u\in L^2$, and that $u$ vanishes outside a compact subset of $\R^{n}$; then $u$ is identically zero.
\end{defin}

\begin{theo}   \label{glm}
Let $V_{\alpha} = V_{\alpha}^{1} + V_{\alpha}^{2} $ be a real-valued
potential satisfying \eqref{va1} and \eqref{va2}. Assume
that $-\Delta + V_{\alpha} (x) {\bf 1}_{|x | < R}$ has the unique
continuation property for all $R$ large enough and 
\begin{equation}   \label{dec}
\< x \>^{1-\alpha /2} \ln \< x \> | V_{\alpha}^{2} (x) | \to 0\,
,\quad\text{as }|x|\to \infty\, .
\end{equation}
Then $H_{\alpha}$ has no eigenvalue.
\end{theo}

\begin{rema}\label{rema:unique}
The unique continuation property for Schr\"{o}dinger operators is
known in many situations but we recall only two cases. The works of
M.~Schechter and B.~Simon \cite{ScSi80} for $n=1$, $2$ and D.~Jerison
and C.~E.~Kenig \cite{JeKe85} imply:
\begin{itemize}
\item[] If $V_{\alpha}^{2} \in L^{p} (\R^{n})$ with $p>1$ for
$n=1$, $2$, and $p \geq n/2$ for $n \geq 3$, then the unique
continuation property holds for $-\Delta + V_{\alpha} (x) {\bf 1}_{|x |
< R}$. 
\end{itemize}
We recall  \cite[Theorem XIII.57]{ReedSimon4}:
\begin{itemize}
\item[] Assume there is a closed set $S$ of measure zero so that
$\R^{n} \setminus S$ is connected, and so that $V_{\alpha}^{2}$ is
bounded on any compact subset of $\R^{n} \setminus S$, then the unique
continuation property holds for $-\Delta + V_{\alpha} (x) {\bf 1}_{|x |
< R}$. 
\end{itemize}
\end{rema}

\begin{proof}[Proof of Theorem~\ref{glm}]
We follow the proof of \cite[Theorem XIII.58]{ReedSimon4}.
Suppose that $u \in D(H_{\alpha})$ is an eigenfunction for
$H_{\alpha}$ with eigenvalue $E$. As in \cite{ReedSimon4}, we define a
function $w$ from $[0, \infty[$ to $L^{2} ({\mathbb S}^{n-1})$ by: 
\begin{align*}
w(r, \omega ) = r^{(n-1)/2} u(r \omega ).
\end{align*}
We have:
\begin{equation}    \label{ii1}
\int_{0}^{+ \infty} \| w(r)\|^{2}_{L^{2} ({\mathbb S}^{n-1} , d\omega
)} dr = \| u \|^{2}_{L^{2} (\R^{n})} < + \infty\, .
\end{equation}
Since $u \in D(H_{\alpha ,0})$, we have $\< x \>^{-\alpha} \Delta u
\in L^{2} (\R^{n})$ and we get, using the pseudo-differential
calculus, that $\(\partial_{x_{j}}\< x \>^{-\alpha} \partial_{x_{k}}
u\)_{1\leq j,k\leq n}$ and $\nabla \< x \>^{-\alpha /2} u$ are in
$L^{2} (\R^{n})$. We have:
\begin{equation*}
\int_{0}^{+\infty} \( \( \{ {\tt d}_{r} \langle r \rangle^{-\alpha /2}w , {\tt d}_{r} \< r 
\>^{-\alpha /2} w \) - r^{-2} \< r \>^{-\alpha} ( w , Bw) \) dr =
\big\| \nabla \< x \>^{-\alpha /2} u \big\|^2, 
\end{equation*}
where $B$ is the Laplace--Beltrami operator on $L^{2} ({\mathbb
S}^{n-1})$ such that $-B \geq 0$, and 
\begin{equation*}
{\tt d}_{r} f(r) = r^{(n-1)/2}
\partial_{r} \( r^{-(n-1)/2} f(r) \).
\end{equation*}
Here $(\cdot,\cdot)$ is the
scalar product on $L^{2} ({\mathbb S}^{n-1})$. Using this formula, we
get 
\begin{equation}   \label{ii2}
 \int_{0}^{+\infty} \< r \>^{-\alpha} \| w' \|^{2} dr < + \infty\quad
; \quad   
 \int_{0}^{+\infty} r^{-2} \< r \>^{-\alpha} (w, -Bw) dr < + \infty\,
,
\end{equation}
and the quantities $ w'$ and $(w, Bw)$ are defined almost everywhere
on $]0, + \infty[$.

Now we define, for $r$ large enough,
\begin{align*}
F(r) = r^{-\alpha} \| w' \|^{2} + r^{-2 - \alpha} (w, Bw) +
\(r^{-\alpha} \< r \>^{\alpha} +E r^{-\alpha}\) \| w \|^{2}. 
\end{align*}
From \eqref{ii1} and \eqref{ii2}, $F(r)$ is integrable. On the other
hand, we have 
\begin{align*}
\( r \ln (r) F(r) \)' =& 2r^{1-\alpha} \ln (r) (w' , w'') +
r^{-\alpha} ((1-\alpha ) \ln (r) +1) \| w' \|^{2} \\ 
&+ 2 r^{-1-\alpha } \ln (r) (w', Bw) - r^{-2- \alpha } ( (1+\alpha )
\ln (r) -1) (w, Bw)  \\ 
&+ r^{-\alpha} \< r \>^{\alpha} \big( (1-\alpha ) \ln (r) + \alpha \ln
(r) r^{2} \< r \>^{-2} + 1 \big) \| w \|^{2}   \\ 
&+ E r^{-\alpha} ((1-\alpha) \ln (r) +1) \| w \|^{2} + 2r \ln (r)
(r^{-\alpha} \< r \>^{\alpha} +E r^{-\alpha}) (w',w). 
\end{align*}
Since $u$ is an eigenfunction of $H_{\alpha}$, we have
\begin{equation}  \label{equ}
w'' = -r^{-2} Bw + \frac{1}{4}(n-1)(n-3) r^{-2} w - \< r \>^{\alpha} w
+ Vw -Ew. 
\end{equation}
Then
\begin{align*}
\big( r \ln (r) F(r) \big)' =& \ln (r) \big( (1-\alpha )r^{-\alpha} \|
w' \|^{2} +  \| w \|^{2} \big) + r^{-\alpha} \| w' \|^{2} + \| w
\|^{2}  \\ 
&- r^{-2- \alpha } ( (1+\alpha ) \ln (r) -1) (w, Bw) +\frac{1}{2}
(n-1)(n-3) r^{-1-\alpha} \ln (r) (w',w)  \\ 
&+ E r^{-\alpha} ((1-\alpha) \ln (r) +1) \| w \|^{2} + {\mathcal O}
(r^{-1} \ln (r)) \| w \|^{2} + 2 r^{1-\alpha} \ln (r) (w' ,Vw). 
\end{align*}
Using $-B \geq 0$ and
\begin{gather*}
r^{-1-\alpha } \ln (r) (w',w) = o (1) r^{-\alpha } \| w' \|^{2} + o(1)
\| w \|^{2},  \\ 
r^{1-\alpha } \ln (r) (w' ,Vw) = o (1) r^{-\alpha /2} \| w ' \| \| w
\| = o (1) r^{-\alpha } \| w' \|^{2} + o(1) \| w \|^{2}, 
\end{gather*}
we get
\begin{equation}   \label{eee}
\big( r \ln (r) F(r) \big)' \geq \ln (r) \big( (1-\alpha )r^{-\alpha}
\| w' \|^{2} +  \| w \|^{2} \big) + \frac{r^{-\alpha}}{2} \| w' \|^{2}
+ \frac{1}{2}  \| w
\|^{2}. 
\end{equation}
for $r$ large enough. Here $o(1)$ denotes a function which tends to $0$
as $r$ tends to $+ \infty$. This computation is formal, but we can give,
as in \cite{ReedSimon4}, a rigorous proof of the integral version of
\eqref{eee}.

If $0 < \alpha < 1$, we get that $r \ln (r) F(r)$ is monotone
increasing for $r\geq R_1$ large enough. Integrate \eqref{eee} between $R_1$
and $r$: $ F(r) \geq \frac{R_1\ln R_1}{r\ln r}F(R_1)$. 
Since $F$ is integrable but
$(r \ln r)^{-1}$ is not, we have $F(r) \leq 0$ for $r > R_{1}$.
\esp

Now, we assume that $1 \leq \alpha \leq 2$. For $1 \ll a < b$, we have
\begin{equation}  \label{ee2}
b \ln b F(b) - a \ln a F(a) \geq \int_{a}^{b} \( \ln (r) \(
(1-\alpha )r^{-\alpha} \| w' \|^{2} +  \| w \|^{2} \) +
\frac{r^{-\alpha}}{2}  \| w' \|^{2} + \frac{1}{2} \| w \|^{2} \) dr. 
\end{equation}
Integration by parts yields:
\begin{align*}
\int_{a}^{b} r^{-\alpha}\ln r  (w'',w) dr =  r^{-\alpha} \ln r
(w',w) \Big|_{a}^{b} -\int_{a}^{b} r^{-\alpha}\ln r  \|
w' \|^{2} dr - \int_{a}^{b} (1-\alpha \ln r)r^{-1-\alpha } (w',w)
dr. 
\end{align*}
The eigenfunction relation \eqref{equ} yields:
\begin{align*}
\int_{a}^{b} r^{-\alpha}\ln r  (w'',w) dr =& - \int_{a}^{b} r^{-2-
\alpha }\ln r (w,Bw) dr + \int_{a}^{b} r^{- \alpha }\ln r  (w,Vw)
dr  \\ 
&+ \int_{a}^{b} r^{- \alpha }\ln r  \( \frac{1}{4}(n-1)(n-3)
r^{-2} -\< r \>^{\alpha} -E \)\| w \|^{2} dr \\ 
=& - \int_{a}^{b} ( \ln r + o(1)) \| w \|^{2} dr - \int_{a}^{b} 
r^{-2- \alpha }\ln r  (w,Bw) dr. 
\end{align*}
We infer:
\begin{align*}
\int_{a}^{b} r^{-\alpha}\ln r  \| w' \|^{2} dr =& \int_{a}^{b} \ln
r \| w \|^{2} dr + \int_{a}^{b} r^{-2- \alpha }\ln r  (w,Bw) dr
\\ 
&+ \big[ r^{-\alpha}\ln r  (w', w) \big]_{a}^{b} + \int_{a}^{b} o(1)
\( r^{-\alpha} \| w' \|^{2} + \| w \|^{2} \) dr, 
\end{align*}
and \eqref{ee2} becomes
\begin{align}
b \ln b \,F(b) - &a \ln a \,F(a) \geq \int_{a}^{b} \( (2- \alpha )
\ln r\, \| w \|^{2} + r^{-\alpha} \| w' \|^{2}/3 + \| w \|^{2}/3 \)
dr  \nonumber  \\ 
&\phantom{a \ln a \,F(a) \geq } + (1-\alpha ) \big[ r^{-\alpha}\ln r
\, (w', w) \big]_{a}^{b} + 
(1-\alpha) \int_{a}^{b} r^{-2- \alpha }\ln r \, (w,Bw) dr   \nonumber
\\  
\label{ee3}
\geq& \int_{a}^{b} \( (2- \alpha )  \ln r\, \| w \|^{2} +
r^{-\alpha} \| w' \|^{2}/3 + \| w \|^{2}/3 \) dr  
+ (1-\alpha ) \big[ r^{-\alpha}\ln r \, (w', w) \big]_{a}^{b}\ ,
\end{align}
since $(1-\alpha) B \geq 0$. Let $\widetilde{F}$ defined by
\begin{equation}
\widetilde{F} (r) = F(r) + (\alpha -1) r^{-\alpha } (w',w),
\end{equation}
which is integrable from \eqref{ii1} and \eqref{ii2}. Inequality
\eqref{ee3} implies that $r \ln r \widetilde{F} (r)$ is monotone
increasing,  and reasoning as before,  $\widetilde{F} (r) \leq 0$ for
$r >R_{1}$.
\esp

We now prove that $w(r)=0$ for $r>R_2$ large enough. For $m\in \N$,
let $w_{m} = r^{m} w$. It satisfies: 
\begin{equation}  \label{eqm}
w_{m}'' = 2mr^{-1} w'_{m} - r^{-2} B w_{m} - \( E+\< r \>^{\alpha}
- V + m(m+1) r^{-2} -\frac{1}{4} (n-1)(n-3) r^{-2} \) w_{m}. 
\end{equation}
We also define
\begin{equation}   \label{deg}
G(r) = r^{2} \| w_{m} ' \|^{2} + m(m+1) \| w_{m} \|^{2} +
(w_{m},Bw_{m}) + \( r^{2} \< r \>^{\alpha} + E r^{2} - r \) \|
w_{m} \|^{2}, 
\end{equation}
and we have:
\begin{align*}
G'(r) =& (4m+2) r \| w_{m}' \|^{2} + 2 \( \( r^{2} V +
\frac{1}{4} (n-1)(n-3) -r \) w_{m} , w_{m}' \) \\ 
&+ \( 2r\< r\>^{\alpha} + \alpha r^{3} \< r \>^{\alpha -2} + 2E r
-1 \) \| w_m \|^{2}. 
\end{align*}
Using \eqref{dec} and the Cauchy--Schwarz inequality, we get
\begin{align*}
G'(r) \geq (4m+1) r \| w_{m}' \|^{2} + r^{1+\alpha} \| w_{m} \|^{2},
\end{align*}
for $r > R_{2} > R_{1}$ independent of $m$. Then $G(r)$ is monotone
increasing on $]R_{2} , + \infty [$.

Suppose that $w(r_{ 0}) \neq 0$ for some $r_{0} > R_{2}$. Since we
have 
\begin{equation*}
G(r) = r^{2m} \( r^{2} \| w ' +mr^{-1} w\|^{2} + m(m+1) \| w \|^{2}
+ (w,Bw) + \( r^{2} \< r \>^{\alpha} + E r^{2} - r \) \| w
\|^{2} \), 
\end{equation*}
we get $G(r_{0}) >0$ for $m \geq 1$ large enough, and now fixed. So,
$G(r) >0$ for all $r > r_{0}$. On the other hand, if $r > r_{1} >
r_{0}$, with $r_{1}$ large enough, we have 
$m(2m+1) - r < 0$.
Since $\| w \|^{2}$ is integrable on $[r_{1}, + \infty[$, the function
$\| w \|^{2}$ is not monotone increasing; there exists $r >
r_{1}$ such that  
\begin{equation*}
\(\| w \|^{2}\)' (r) = 2 (w',w) (r) \leq 0.
\end{equation*}
Then
\begin{align*}
G(r) =& r^{2m} \Big( r^{2} \| w '\|^{2} + 2 mr (w' ,w) + m(2m+1) \| w
\|^{2} + (w,Bw) + \big( r^{2} \< r \>^{\alpha} + E r^{2} - r \big) \|
w \|^{2} \Big)  \\ 
\leq& r^{2m+2+\alpha} \Big( r^{-\alpha} \| w' \|^{2} + r^{-2-\alpha }
(w,Bw) + \big( r^{-\alpha} \< r \>^{\alpha} + E r^{-\alpha} \big) \| w
\|^{2} \Big) + 2 mr^{-1-\alpha} (w',w). 
\end{align*}
Therefore, we get
\begin{equation*}
\begin{aligned}
0& < G(r) \leq r^{2m+2+\alpha} F(r) \leq 0 \quad &&\text{ if } 0 < \alpha <1,   \\
0& < G(r) \leq r^{2m+2+\alpha} \widetilde{F}(r) \leq 0 &&\text{ if } 1 \leq \alpha \leq 2,
\end{aligned}
\end{equation*}
which is impossible, so $w (r)=0$ for $r > R_{2}$. Theorem
\ref{glm} follows from unique continuation. 
\end{proof}

\section{Mourre estimates}
\label{sec:mourre}

In the following, we use the Weyl calculus of L.~H\"{o}rmander, for
which  we refer to \cite[Section~XVIII]{Ho85}. More precisely, we
work with the simple metrics which are $\sigma$--temperate:  
\begin{align*}
g_{0} =& |dx|^{2} + |d \xi|^{2}\ , \\
g_{1} =& \frac{|dx|^{2}}{1+x^{2} + \xi^{2}} + \frac{|d
\xi|^{2}}{1+x^{2} + \xi^{2}}\ , \\ 
g_{1}^{\beta} =& \frac{|dx|^{2}}{(1+x^{2} + \xi^{2})^{\beta}} +
\frac{|d \xi|^{2}}{(1+x^{2} + \xi^{2})^{\beta}}\ , \\ 
g_{2} =& \frac{|dx|^{2}}{1+x^{2}} + \frac{|d \xi|^{2}}{1+\xi^{2}}\ ,
\end{align*}
for $\beta >0$. We refer to \cite{Ho85} for the definition of the
space of symbol $S(m,g)$ and we note $\Psi (m,g)$ the set of
pseudo-differential operators whose symbol is in a space $S(m,g)$. We
set $\Psi (g) = \bigcup_{m} \Psi (m,g)$.

The crucial point of the Mourre theory is the construction of the
conjugate operator. This is a self-adjoint operator $A_{\alpha}$ such
that $i[H_{\alpha} , A_{\alpha}]$ is positive on the energy level, and
$H_{\alpha}$--bounded. In our case, the generator of dilations
$ (xD+Dx)/2$ is not satisfactory for $H_{\alpha,0}$ since 
\begin{equation*}
i \left[H_{\alpha ,0} , (x D + D x) /2 \right] = -2 \Delta + \alpha x^{2} \< x
\>^{\alpha -2}, 
\end{equation*}
which is positive, but not $H_{\alpha,0}$--bounded for $\alpha >0$. So
we must find another conjugate operator. We look for $A_{\alpha}$ as a
pseudo-differential operator of symbol $a_{\alpha} (x,\xi)$. Consider
the case of dimension one, and start with $\alpha = 2$. Formally, we
want to solve:  
\begin{equation}  \label{mm1}
\left\{ \xi^{2} - x^{2} , a_{2} (x,\xi) \right\} \equiv 2 \xi
\partial_{x} a_{2} + 2 x 
\partial_{\xi} a_{2} = 1,\quad \text{on } \left\{(x,\xi)\ ;\ \xi^{2} -
x^{2} = E\right\}.  
\end{equation}
We saw in the introduction that a solution to this equation is given by:
\begin{equation}  \label{tx14}
a_{2} (x,\xi) = \frac{1}{4} \big( \ln (\xi +x) - \ln ( \xi -x) \big).
\end{equation}
Now consider the case $0 < \alpha < 2$. Replacing $\< x \>$ by $x>0$, we
try to solve  
\begin{equation}  \label{mm2}
\left\{ \xi^{2} - x^{\alpha} , a_{\alpha} (x,\xi) \right\} \equiv 2\xi
\partial_{x} 
a_{\alpha} + \alpha x^{\alpha -1} \partial_{\xi} a_{\alpha} = 2-\a,\quad
 \text{on } \left\{(x,\xi)\ ;\ \xi^{2} - x^{\alpha} = E\right\}. 
\end{equation}
As we saw in the introduction, $a_{\alpha} (x,\xi) = \xi x^{1-\alpha}$ should do the job, up to an error which is compact on the energy level (because it decays like $x^{-\alpha}$). In this section, we make this heuristic approach rigorous. The main results (Mourre estimates) are proved in Section~\ref{moumou}. Here, we can find again the short range condition: on the energy level, we have formally
\begin{align*}
\< \xi \> \approx \< x \>^{\alpha /2}
\text{ and then }
\left| a_{\alpha} (x, \xi ) \right| \approx \left\{
\begin{aligned}
&\ln \< x \> \quad && \text{if } \alpha =2,   \\
&\< x \>^{1-\alpha /2} \quad && \text{if } 0< \alpha <2 .
\end{aligned}
\right.
\end{align*}
By (\ref{mm1}) and (\ref{mm2}), we obtain that the position variable increases exactly like $t$ along the evolution. Then, in section \ref{mimi}, we will replace $a_{\alpha} (x, \xi )$ by $p_{\alpha} (x)$ and we will require that the potential decays as $p_{\alpha} (x)^{-1-\varepsilon}$.

\subsection{General framework}

We recall some results that we will use to prove 
regularity results on the groups generated by $H_{\alpha}$. A full
presentation of such issues can be found in the book of O.~Amrein, A.~Boutet de Monvel and V.~Georgescu \cite{AmBoGe90}. We start with the defintion of $C^{1} (A)$.

\begin{defin}
Let $A$ and $H$ be self-adjoint operators on a Hilbert space
${\mathcal H}$. We say that $H$ is of class $C^{r} (A)$ for $r>0$, if there is $z \in \C \setminus \sigma (H)$ such that
\begin{align*}
\R \ni t \to e^{it A} (H-z)^{-1} e^{-itA},
\end{align*}
is $C^{r}$ for the strong topology of ${\mathcal L} ({\mathcal H})$.
\end{defin}

We have the following useful characterization of the regularity $C^{1} (A)$. 

\begin{theo}[{\cite[Theorem~6.2.10]{AmBoGe90}}]  \label{ABG1}
Let $A$ and $H$ be self-adjoint operators on a Hilbert space
${\mathcal H}$. Then $H$ is of class $C^{1}(A)$ if and only if the
following conditions are satisfied:  
\begin{enumerate}
\item There exists $c< \infty$ such that for all $u \in D(A)
\cap D(H)$, 
\begin{equation*}
|(Au,Hu) - (Hu,Au) | \leq c \left( \| H u \|^{2} + \| u\|^{2} \right).
\end{equation*}

\item For some $z \in \C\backslash \sigma (H)$, the set $\{ u \in
D(A); \ (H-z)^{-1} u \in D(A) \text{ and } (H-\bar{z})^{-1} u \in D(A)
\}$ is a core for $A$. 
\end{enumerate}
If $H$ is of class $C^{1} (A)$, then the following holds:
\begin{enumerate}
\item The space $(H-z)^{-1} D(A)$ is independent of $z\in \C \backslash
\sigma (H)$, and contained in $D(A)$. It is a core for $H$, and a dense
subspace of $D(A) \cap D(H)$ for the intersection topology (i.e. the
topology associated to the norm $ \| H u\| + \| Au\| + \| u \|)$.

\item The space $D(A) \cap D(H)$ is a core for $H$, and the form
$[A,H]$ has a unique extension to a continuous sesquilinear form on
$D(H)$ (equipped with the graph topology). If this extension is
denoted by $[A,H]$, the following identity holds on ${\mathcal H}$ (in
the form sense): 
\begin{equation*}
\left[ A, (H-z)^{-1} \right] = - (H-z)^{-1} [A,H] (H-z)^{-1},
\end{equation*}
for $z \in \C \backslash \sigma (H)$.
\end{enumerate}
\end{theo}

We also have the following theorem from \cite[Theorem~6.3.4]{AmBoGe90}: 

\begin{theo}   \label{ABG2}
Let $A$ and $H$ be self-adjoint operators in a Hilbert space
${\mathcal H}$. Assume that the unitary one-parameter group $\{ \exp
(iA \tau) \}_{\tau \in \R}$ leaves the domain $D(H)$ of $H$
invariant. Then $H$ is of class $C^{1}(A)$ if and only if $[H,A]$ is
bounded from $D(H)$ to $D(H)^{*}$. 
\end{theo}

A criterion for  the above assumption to be satisfied is given by the
following result:   
\begin{lem}[{\cite[Lemma~2]{GeGe99}}]  \label{prm}
Let $A$ and $H$ be self-adjoint operators in a Hilbert space
${\mathcal H}$. Let $H \in C^{1} (A)$ and suppose that the commutator
$[iH,A]$ can be extended to a bounded operator from $D(H)$ to
${\mathcal H}$. Then $e^{itA}$ preserves $D(H)$.  
\end{lem}

\subsection{A technical result}
\label{ddd}

It is not clear that the energy cut-offs $\chi (H_{\alpha ,0})$, $\chi
\in C_{0}^{\infty} (\R)$, are pseudo-differential operators, since
$H_{\alpha ,0}$ is not elliptic, and $\chi \( \xi^{2} - \< x \>^{\alpha} \)$ is not a good symbol. The next proposition will allow us to use pseudo-differential calculus. Such techniques have been used by M. Dimassi and V. Petkov \cite{DiPe03}.

\begin{prop}   \label{pse}
Let $0 < \alpha \leq 2$, $1/2 < \beta \leq 1$, $z \in \C \setminus \R$
and $\psi \in C_{0}^{\infty} (\R)$ such that $\psi = 1$ near $0$. Then 
\begin{align*}
(H_{\alpha,0} -z)^{-1} = (H_{\alpha,0} -z)^{-1} {\rm Op} \left( \psi
\left( \frac{\xi ^{2} 
-\< x \>^{\alpha}}{( \xi^{2} + \< x \>^{\alpha} )^{\beta}} \right)
\right) + {\mathcal O} 
(1) {\rm Op} \left( r \right), 
\end{align*}
with $r (x,\xi ) \in S \( (\xi^{2} + \< x \>^{\alpha})^{-\beta}
,g_{0} \)$.  
\end{prop}

\begin{proof}
For $\gamma \geq 0$, let $g_{3}^{\alpha , \gamma}$ be the 
$\sigma$--temperate metric
\begin{align*}
g_{3}^{\alpha , \gamma} = |dx|^{2} + \frac{|d \xi|^{2}}{(\xi^{2} + \< x \>^{\alpha})^{\gamma}}.
\end{align*}
Let $B = {\rm Op} (b)$, with $b = ( \xi ^{2} - \< x \>^{\alpha} )/(
\xi^{2} + \< x \>^{\alpha} )^{\beta} \in S \( ( \xi^{2} + \< x
\>^{\alpha} )^{1-\beta},g_{3}^{\alpha ,1} \)$. It satisfies 
\begin{equation}  \label{mpp}
\partial_{x,\xi}^{\delta} b \in S \( (\xi^{2} + \< x
\>^{\alpha})^{1-\beta -\min (1, |\delta|/2)} ,g_{3}^{\alpha ,1} \)
, \quad \forall \delta \in \N^{2n}\,.
\end{equation}
We have
\begin{equation}  \label{aqq1}
(H_{\alpha,0} -z)^{-1} \left( 1- \psi (B) \right) = (H_{\alpha,0}
-z)^{-1} B B^{-1} \left( 1- \psi (B) \right). 
\end{equation}
Theorem~18.5.4 of \cite{Ho85} on the composition of
pseudo-differential operators in $\Psi (g_{3}^{\alpha ,1})$, and
 \eqref{mpp} imply that  
\begin{align*}
(H_{\alpha ,0} -z)^{-1} B = (H_{\alpha ,0} -z)^{-1} \big( H_{\alpha ,0} {\rm Op} \big( ( \xi^{2}
+ \< x \>^{\alpha} )^{-\beta} \big) + {\rm Op} (r) \big), 
\end{align*}
where $r \in S \( ( \xi^{2} + \< x \>^{\alpha} )^{-\beta} ,
g_{3}^{\alpha ,1} \)$. So we have 
\begin{equation}   \label{aqq2}
(H_{\alpha ,0} -z)^{-1} B = {\mathcal O} (1)  {\rm Op} (r),
\end{equation}
for some other $r \in S \( ( \xi^{2} + \< x \>^{\alpha} )^{-\beta} ,
g_{3}^{\alpha ,1} \)$.

Let $\varphi(y) = (1-\psi (y))/y$, and $\widetilde{\varphi}$ be an
almost analytic extension of $\varphi$ (see \cite{Ho68}, \cite{Da95}
and \cite[Appendix~C.2]{DG}). This is a $C^{\infty} (\C)$ function
which coincides with $\varphi$ on $\R$, whose support is contained in
a region like $|{\rm Im} \, z| < C \< {\rm Re} \, z \>$, and which
satisfies  
\begin{align*}
\left| \partial_{\overline{z}} \widetilde{\varphi} (z) \right| \leq
C_{k} \< z \>^{-2-k} |{\rm Im} \, z |^{k}\, ,\quad \forall k\in\N\, . 
\end{align*}
Using the Helffer--Sj\"{o}strand formula (see
\cite{HeSo89} or \cite{Da95}), we can write 
\begin{equation}
{\rm Op} (r) B^{-1} \left( 1- \psi (B) \right) = \frac{1}{\pi} \int
\partial_{\overline{z}} \widetilde{\varphi} (z) {\rm Op} (r)
(B-z)^{-1} L(dz).   \label{mal1}  
\end{equation} 
For ${\rm Im} \, z \neq 0$ and $1/2 < \beta \leq 1$, we have
$(b(x,\xi) -z)^{-1} \in S \( 1,g_{3}^{\alpha, 2 \beta -1} \)$
and 
\begin{equation}  \label{mpp2}
\partial_{x,\xi}^{\delta} (b(x,\xi) -z)^{-1} \in S \( (\xi^{2} + \<
x \>^{\alpha})^{-(\beta -1/2) \min (2, |\delta|)} ,g_{3}^{\alpha, 2
\beta -1} \), \quad \forall  \delta \in \N^{2n}\, .
\end{equation}
Using \cite[Theorem~18.5.5]{Ho85} on the composition of 
pseudo-differential operators in $\Psi \big( g_{3}^{\alpha, 1} \big)$,
and $\Psi \big( g_{3}^{\alpha, 2 \beta -1} \big)$, \eqref{mpp} and
\eqref{mpp2}, we have, for ${\rm Im} \, z \neq 0$,  
\begin{align*}
(B-z) {\rm Op} \left( (b(x,\xi) -z)^{-1} \right) = 1 + {\rm Op} (d(z)),
\end{align*}
where $d (z) \in S \big( (\xi^{2} + \< x \>^{\alpha} )^{1-3 \beta},
g_{3}^{\alpha, 2\beta -1} \big)$, and each semi-norm of $d(z)$ in this
space is bounded by some power of $1+ |{\rm Im} \, z |^{-1}$. On the
other hand, we have 
\begin{align*}
{\rm Op} \left( (\xi^{2} + \< x \>^{\alpha} +i \lambda)^{1-3 \beta}
\right) {\rm Op} \left( (\xi^{2} + \< x \>^{\alpha} +i \lambda)^{3
\beta -1} \right) = 1+{\mathcal O} (\lambda^{-1})\,, 
\end{align*}
and then, for $\lambda$ large enough,
\begin{align*}
{\rm Op} (d(z)) =&\,  {\rm Op} (d(z))  {\rm Op} \left( (\xi^{2} + \< x
\>^{\alpha} +i \lambda)^{3 \beta -1} \right) {\mathcal O} (1) {\rm Op} \left( (\xi^{2}
+ \< x \>^{\alpha} +i \lambda)^{3 \beta -1} \right)^{-1}     \\ 
=&\  {\rm Op} (\widetilde{d}(z)) {\mathcal O} (1) {\rm Op} \left(
(\xi^{2} + \< x \>^{\alpha} +i \lambda)^{1-3 \beta} \right), 
\end{align*}
with $\widetilde{d} (z) \in S \big(1, g_{3}^{\alpha, 2\beta -1} \big)$
and each semi-norm is bounded by some power of $1+ |{\rm Im} \, z
|^{-1}$. The continuity in $L^{2} (\R^{n})$ of 
pseudo-differential operators yields  
\begin{equation*}
(B-z) {\rm Op} \left( (b(x,\xi) -z)^{-1} \right) = 1 + {\mathcal O}
(1+|{\rm Im} \, z |^{-M}) {\rm Op} \left( (\xi^{2} + \< x \>^{\alpha}
+i \lambda)^{1-3 \beta} \right),  
\end{equation*}
for some $M>0$. We infer:
\begin{equation}       \label{mal2}
(B-z)^{-1} =  {\rm Op} \left( (b(x,\xi) -z)^{-1} \right) + {\mathcal
O}(1+|{\rm Im} \, z |^{-M-1}) {\rm Op} \left( (\xi^{2} + \< x
\>^{\alpha} +i \lambda)^{1-3 \beta} \right) .  
\end{equation}
Then, using the pseudo-differential calculus, \eqref{mal1} becomes
\begin{equation}\label{aqw3} 
\begin{aligned}
{\rm Op} (r) B^{-1} \left( 1- \psi (B) \right) =& \,\frac{1}{\pi} \int
\partial_{\overline{z}} \widetilde{\varphi} (z) {\rm Op} (r) {\rm Op}
\left( (b(x,\xi) -z )^{-1} \right) L(dz) \\
&+ {\mathcal O} (1) {\rm Op}
\left( (\xi^{2} + \< x \>^{\alpha} +i \lambda)^{1-3 \beta} \right)\\
=&\, {\mathcal O} (1) {\rm Op} (r),   
\end{aligned}
\end{equation}
for some other $r \in S \big( (\xi^{2} + \< x \>^{\alpha} )^{-\beta},
g_{3}^{\alpha, 2\beta -1} \big)$. From \eqref{aqq1}, \eqref{aqq2} and
\eqref{aqw3}, we have 
\begin{equation}   \label{fgh1}
(H_{\alpha ,0} -z)^{-1} = (H_{\alpha ,0} -z)^{-1} \psi (B) + {\mathcal
O} (1) {\rm Op} (r), 
\end{equation}
Using the Helffer--Sj\"{o}strand formula:
\begin{align*}
\psi (B) &= \frac{1}{\pi} \int \partial_{\overline{z}}
\widetilde{\psi} (z) (B-z)^{-1} \, L(dz),   \\
{\rm Op} (\psi (b(x,\xi)) &= \frac{1}{\pi} \int
\partial_{\overline{z}} \widetilde{\psi} (z) {\rm Op} \left(
(b(x,\xi) -z )^{-1} \right) \, L(dz) ,
\end{align*}
and \eqref{mal2}, we obtain:
\begin{equation}  \label{fgh2}
\psi (B) = {\rm Op} (\psi (b(x,\xi)) + {\mathcal O}(1) {\rm Op} (r),
\end{equation}
with $r \in S \big( (\xi^{2} + \< x \>^{\alpha} )^{-\beta},
g_{3}^{\alpha, 2\beta -1} \big)$. The proposition follows from
\eqref{fgh1} and \eqref{fgh2}. 
\end{proof}

\subsection{Conjugate operator}
\label{iii}

Following the discussion of the beginning of Section~\ref{sec:mourre},
we choose for the conjugate operator if $\alpha =2$,
\begin{equation}   \label{nnn}
A_2 = {\rm Op} \left( a_2 (x,\xi) \right),
\quad\text{with}\quad
a_2 = \big( \ln \< \xi + x\> - \ln \< \xi - x\> \big) .
\end{equation}
One can see that $a_{2} (x,\xi) \in S \left( \< \ln \< x \> \> , g_{0} \right)$. Indeed, we have, for $|\xi| < 2 |x|$,
\begin{equation}   \label{mm3}
\ln \< \xi + x \> - \ln \< \xi -x \> \leq \ln \< 3 x \> \leq \< \ln \< x \> \> + C,
\end{equation}
with $C>0$. On the other hand, we get for $|\xi| \ge 2 |x|$,
\begin{equation}  \label{mm4}
\ln \< \xi + x \> - \ln \< \xi -x \> = \ln \left( \frac{1+ (\xi +x)^{2}}{1+ (\xi -x)^{2}} \right) \leq \ln \left( \frac{1+ 9 \xi^{2} /4}{1+ \xi^{2} /4} \right) \leq C,
\end{equation}
with $C>0$. For computational reasons, it is better to have a another writing for $A_2$.

\begin{lem}  \label{xxx}
We have
\begin{equation*}
A_2 = U \big( \ln \big\< \sqrt{2} x \big\> - \ln \big\< \sqrt{2} D \big\> \big) U^{*}.
\end{equation*}
\end{lem}

\begin{proof}
Using the exact composition of pseudo-differential operators (Theorem~18.5.4 of \cite{Ho85}), we have  
\begin{gather*}
{\rm Op} \big( \ln \< \xi + x \> \big) = {\rm Op} \big( (\xi +x)^{2} +1 \big) {\rm Op} \big( \ln \< \xi + x \> ((\xi +x)^{2} +1)^{-1} \big), \\ 
{\rm Op} \big ((\xi +x)^{2} +z)^{-1} \big) = \big( (D+x)^{2} +z \big)^{-1}, 
\end{gather*}
for ${\rm Im} \, z \neq 0$. Let $\gamma \subset \C$ be a contour
enclosed $[0,+\infty [$ in the region where $\ln (z+1) (z+1)^{-1}$ is
holomorphic and coinciding with ${\rm Re} \, z = |{\rm Im} \, z|$ for
$z$ large enough. Using the Cauchy formula and \eqref{pmo}, we get 
\begin{align}
{\rm Op}  \big( \ln \< \xi + x \> \big) =& \frac{1}{2i \pi} \big( (D +x)^{2} +1 \big) \int_{\gamma} \ln (z+1) (z+1)^{-1}  {\rm Op} \big( ((\xi +x)^{2} -z \big)^{-1} \big) dz \nonumber \\ 
=& \frac{1}{2i \pi} \big( (D +x)^{2} +1 \big) \int_{\gamma} \ln (z+1)
 (z+1)^{-1} \big( (D +x)^{2} -z \big)^{-1} dz   \nonumber \\ 
=& U \frac{1}{2i \pi} \big( 2 x^{2} +1 \big) \int_{\gamma} \ln (z+1) (z+1)^{-1} \big( 2 x^{2} -z \big)^{-1} dz \, U^{*}  \nonumber \\ 
=& U \ln \< \sqrt{2} x \> U^{*} , 
\end{align}
and the lemma follows.
\end{proof}

In the case $0<\alpha<2$, we choose for the conjugate operator
$A_{\alpha} = {\rm Op} \left( a_{\alpha} (x,\xi) \right)$, where 
\begin{equation}  \label{nji}
a_{\alpha} (x,\xi) = x \cdot\xi \< x \>^{-\alpha} \psi \left( \frac{\xi^{2}
-\< x \>^{\alpha}}{\xi^{2} + \< x \>^{\alpha}} \right) \in S \( \<
x \>^{1-\alpha/2} , g_{1}^{\alpha /2} \) \cap S \( \< \xi \> \<
x \>^{1-\alpha} , g_{2} \) , 
\end{equation}
with $\psi \in C_{0}^{\infty} ([-1/2, 1/2])$, and  $\psi = 1$ near
$0$. Notice that on $\supp a_{\a}$, $| \xi|$ is like $\< x
\>^{\alpha /2}$.

\subsection{Regularity results}
\label{qqq}

The aim of this section is to prove some regularity results for $H_{\alpha}$. First, we give a common core for the operators $H_{\alpha}$ and $A_{\alpha}$. Using the results of the previous section, we get

\begin{lem}  \label{nel01}
Let $0<\alpha\leq 2$. The operator $A_{\alpha}$ is essentially
self-adjoint on $D(N_2)$. 
\end{lem}

\begin{proof}
As in the proof of Lemma \ref{lem:essadj0}, we use
Theorem~\ref{nel}. We distinguish the cases $\alpha=2$
and $0<\alpha<2$. First,  we suppose that $0<\alpha<2$; 
for $u \in D(N)$, we have  
\begin{align*}
\| A_{\alpha} u \| \lesssim \| \< x \>^{1-\alpha/2} u \| \lesssim \| N_2 u\|,
\end{align*}
and the composition rules for an operator in $\Psi (g_{1})$ by an
operator in $\Psi (g_{1}^{\alpha /2})$ yield 
\begin{align*}
[ A_{\alpha} , N_2 ] \in \Psi \( \< x, \xi \>^{1 - \alpha /2} \< x
\>^{1-\alpha/2}, g_{1}^{\alpha/2} \),  
\end{align*}
which implies \eqref{nel2}, from \cite[Theorem~18.6.3]{Ho85}. The
lemma follows for $0<\alpha<2$. When $\alpha=2$, 
\begin{gather*}
\| A_2 u \| \lesssim \| \< \ln \< x \> \> u \| \lesssim \|N_2u\|, 
\end{gather*}
which proves \eqref{nel1}. Moreover, since $a_{2} (x,\xi) \in S \left(
\< \ln \< x \> \> , g_{0} \right)$ and $N_{2} \in \Psi \left( \< x,
\xi \>^{2} ,g_{1} \right)$, we get $[ N_2 , A_2 ] \in \Psi \left( \<
\ln \< x \> \> \< x, \xi \> ,g_{0} \right)$ and then 
\begin{align*}
\big| \( [ N_2 , A_2 ] u ,u \) \big| \lesssim \left\| N_2^{1/2} u 
\right\|^{2},
\end{align*}
which yields \eqref{nel2} and the lemma.
\end{proof}

\subsubsection{Regularity for $H_{2,0}$}

\begin{lem}  \label{rere}
For $z \in \C \backslash \R$, $(H_{2,0} -z)^{-1}$ maps $D(N_2)$ into itself.
\end{lem}

\begin{proof} We use the notations of Section~\ref{sec:mehler}. For $u
\in D(N_2)$, we have 
\begin{align*}
\| x^{2} e^{-it H_{2,0}} u \| =&\, \| x^{2} {\mathcal M}_{t} {\mathcal
D}_{t} {\mathcal F} {\mathcal M}_{t} u \| = \| x^{2} {\mathcal D}_{t}
{\mathcal F} {\mathcal M}_{t} u \|  \\ 
=&\, \| ( \sinh 2t)^{2} x^{2} {\mathcal F} {\mathcal M}_{t} u \| = \| - (
\sinh 2t)^{2} \Delta {\mathcal M}_{t} u \|  \\ 
=&\, \| {\mathcal M}_{t} \big( - (\sinh 2t)^{2} \Delta + (\cosh 2t)^{2}
x^{2} - \tanh 2t (xD+Dx) \big) u \|  \\
\lesssim& \,  e^{4|t|} \| N_2u\|.
\end{align*}
We also have:
\begin{align*}
\| -\Delta e^{-it H_{2,0}} u \| = \| (H_{2,0} + x^{2} ) e^{-it
H_{2,0}} u \|  
\lesssim \| H_{2,0} u \| + e^{4|t|} \| N_2u\| 
\lesssim e^{4|t|} \| N_2u\|.
\end{align*}
So, for ${\rm Im} \, z >4$, we get
\begin{align*}
\| N_2 (H_{2,0} - z)^{-1} u \| =\left\| iN_2 \int_{0}^{+\infty} e^{itz}
e^{-it H_{2,0}} u \, dt \right\|  
\lesssim \int_{0}^{+ \infty} e^{- t {\rm Im} \, z} e^{4t} \, dt \,
\| N_2u\|
\lesssim \| N_2u \|,
\end{align*}
which shows that $(H_{2,0} - z)^{-1}$ maps $D(N_2)$ into itself for ${\rm
Im} \, z >4$. Then the lemma follows from \cite[Lemma~6.2.1]{AmBoGe90}. 
\end{proof}

\begin{lem}  \label{comutH0}
$H_{2,0}$ is in $C^{1} (A_2)$ and $[H_{2,0},A_2]$ is bounded on $L^{2} (\R^{n})$.
\end{lem}

\begin{proof}
From Theorem \ref{ABG1}, it is enough to estimate $[H_{2,0},A_2]$. Recall that from (\ref{eq:defN_Htilde}), $\widetilde{H}_{2,0} = D_x x + x D_x-1$. Using (\ref{mmm}) , Lemma \ref{xxx} and $\widetilde{H}_{2,0} +1 = -{\mathcal F} ( \widetilde{H}_{2,0} +1) {\mathcal
F}^{*}$, we have:
\begin{align}
[H_{2,0},A_2] =& U \big[\widetilde{H}_{2,0}, \ln \big\< \sqrt{2} x \big\> - \ln \big\< \sqrt{2} D \big\> \big] U^{*} \nonumber \\ 
=& U \big[ \widetilde{H}_{2,0}, \ln \big\<  \sqrt{2}x \big\> \big] U^{*} + U {\mathcal F} \big[ \widetilde{H}_{2,0}, \ln \big\< \sqrt{2} x \big\> \big] {\mathcal F}^{*} U^{*} \nonumber \\ 
=& -i U \frac{4x^{2}}{\big\< \sqrt{2} x \big\>^{2}} U^{*} -i U {\mathcal F} \frac{4 x^{2}}{\big\< \sqrt{2} x \big\>^{2}} {\mathcal F}^{*} U^{*} \nonumber \\
=& -i U \bigg( \frac{4x^{2}}{\big\< \sqrt{2} x \big\>^{2}} + \frac{4 D^{2}}{\big\< \sqrt{2} D \big\>^{2}} \bigg) U^{*}, \label{com1} 
\end{align}
which is bounded on $L^{2} (\R^{n})$.
\end{proof}

For the asymptotic completeness, we need more regularity. We begin with
\begin{lem}\label{doublecom}
$H_{2,0}$ is in $C^{2} (A_2)$ and $[[H_{2,0},A_2],A_2]$ is bounded on $L^{2} (\R^{n})$.
\end{lem}

\begin{proof}
Since we know that $H_{2,0}$ is in $C^{1} (A_{2})$, it is enough to prove
that $[[H_{2,0},A_{2}],A_{2}]$ is bounded. From (\ref{com1}), we can write 
\begin{equation*}
[[H_{2,0},A_2],A_2] = -i U \bigg[ \frac{4x^{2}}{\< \sqrt{2} x\>^{2}} + \frac{4 D^{2}}{\big\< \sqrt{2} D \big\>^{2}} , \ln \big\< \sqrt{2} x \big\> - \ln \big\< \sqrt{2} D \big\> \bigg] U^{*}. 
\end{equation*}
The symbols
\begin{align*}
f(x,\xi) =& \frac{2x^{2}}{\< x\>^{2}} + \frac{2 \xi^{2}}{\< \xi \>^{2}},  \\
g(x,\xi) =& \ln \< x\> - \ln \< \xi \>,
\end{align*}
satisfy $f \in S(1,g_{2})$ and $g \in S(\ln \< x \> + \ln \< \xi \>,
g_{2})$. Then, from Theorems~18.5.4 and 18.6.3 of \cite{Ho85},
we have $\left[ {\rm Op} (f), {\rm Op}(g) \right]= {\mathcal O} (1)$
which completes the proof. 
\end{proof}

\subsubsection{Regularity for $H_{\alpha ,0}$ for $0<\alpha<2$}

\begin{lem}  \label{tre}
The operator $[H_{\alpha ,0} , A_{\alpha} ]$ is in $\Psi (1, g_{2} )$,
and its symbol is supported inside the support of $a_{\alpha}
(x,\xi)$, modulo $S \big( \< x , \xi \>^{-\infty} ,
g_{2}  \big)$. 
\end{lem}

\begin{proof}
Since $a_{\alpha}  (x,\xi) \in S \big( \< x \>^{1-\alpha} \< \xi \> ,
g_{2} \big)$ and $ \xi^{2} - \< x \>^{\alpha} \in S \big( \< \xi
\>^{2} + \< x \>^{\alpha} , g_{2} \big)$, we get $[A_{\alpha} ,
H_{\alpha ,0}] \in \Psi \big( 1+ \< x \>^{-\alpha} \< \xi \>^{2} ,
g_{2} \big)$, and each term in the development of its symbol is
supported inside the support of $a_{\alpha}  (x, \xi)$. Since $\< \xi
\>$ is like $\< x \>^{\alpha/2}$ on the support of $a_{\alpha}
(x,\xi)$, we get the lemma.  
\end{proof}

\begin{lem}  \label{xac0}
$H_{\alpha ,0}$ is in $C^{1} (A_{\alpha})$, and $[H_{\alpha ,0},
A_{\alpha}]$ is bounded on $L^{2} (\R^{n})$.  
\end{lem}

\begin{proof}
As for Lemma \ref{tre}, we have
\begin{equation}   \label{ddt1}
[A_{\alpha} , N_{\alpha}] \in \Psi \big( 1 , g_{2} \big),
\end{equation}
and its symbol is supported inside the support of $a_{\alpha}
(x,\xi)$, modulo a term in $S \big( \< x , \xi
\>^{-\infty} , g_{2}  \big)$. Then $[A_{\alpha}  , N_{\alpha}]$ is
bounded on $L^{2} (\R^{n})$. On the other hand, from the
pseudo-differential calculus, one can show that $( N_{\alpha} +i
)^{-1}$ maps $D(N_2)$ into itself. Then, from Theorem \ref{ABG1},
$N_{\alpha}$ is $C^{1} ( A_{\alpha})$.

Since $[A _{\alpha} , N_{\alpha} ]$ is bounded, we get from the proof
of \cite[Lemma~2]{GeGe99}, that $e^{itA_{\alpha}}$ preserves
$D(N_{\alpha})$ and that 
\begin{equation}   \label{ddt2}
N_{\alpha} e^{itA_{\alpha}} = e^{itA_{\alpha}} N_{\alpha} + i
\int_{0}^{t} e^{i(t-s)A_{\alpha}} [N_{\alpha}, A_{\alpha} ]
e^{isA_{\alpha}} ds, 
\end{equation}
on $D(N_{\alpha})$. From \eqref{ddt1}, we infer that $[[A_{\alpha} ,
N_{\alpha}] , N_{\alpha}] \in \Psi \big( \< x \>^{-1} \< \xi \> + \< x
\>^{\alpha -1} \< \xi \>^{-1} , g_{2} \big)$ and that each term in 
the development of its symbol is supported inside the support of
$a_{\alpha}  (x, \xi)$. Then 
$[[A_{\alpha} , N_{\alpha}] , N_{\alpha}] \in \Psi \big( \< x
\>^{\alpha /2 -1} , g_{2} \big) \subset \Psi \big( 1 , g_{2} \big)$
because $0< \alpha <2$. By induction, we obtain that 
\begin{equation}   \label{ddt3}
[ \ldots [[[ A_{\alpha} , N_{\alpha} ] , N_{\alpha} ], N_{\alpha} ],
\ldots N_{\alpha} ] \in \Psi \big( 1, g_{2} \big). 
\end{equation} 
Using \eqref{ddt2} and \eqref{ddt3}, we obtain that $e^{itA_{\alpha}}$
preserves $D(N_{\alpha}^{k})$ for all $k \in \N$ and that 
\begin{equation}
\big\| N_{\alpha}^{k} e^{itA_{\alpha}} u \big\| \lesssim \left\|
N_{\alpha}^{k} u \right\| + \| u \|, 
\end{equation}
for all $u \in D(N_{\alpha}^{k})$. Since $\alpha \neq 0$, there is $k
\in \N$ such that $e^{itA_{\alpha}}$ maps continuously $D(N^{k}_\a)$ into
$D (N^{2}_\a)$. Then 
\begin{align*}
t \mapsto e^{i(t-s)A_{\alpha}} H_{\alpha ,0} e^{isA_{\alpha}} u,
\end{align*}
is well-defined and $C^{1}$ for $u \in D(N_\a^{k})$. It follows that
\begin{equation}   \label{ddt4}
H_{\alpha ,0} e^{itA_{\alpha}} = e^{itA_{\alpha}} H_{\alpha ,0} + i
\int_{0}^{t} e^{i(t-s)A_{\alpha}} [H_{\alpha ,0}, A_{\alpha} ]
e^{isA_{\alpha}} ds, 
\end{equation}
on $D(N_\a^{k})$. Using Lemma~\ref{tre}, $[H_{\alpha ,0}, A_{\alpha}]$
can be extended as a bounded operator on $L^{2} (\R^{n})$. On the other
hand $H_{\alpha ,0}$ satisfies the Nelson's Theorem~\ref{nel} with
$N^{k}_\a$ as reference operator. Then \eqref{ddt4} can be extended on
$D(H_{\alpha ,0})$ and $e^{itA_{\alpha}}$ preserves $D( H_{\alpha
,0})$.

Since  $[H_{\alpha ,0}, A_{\alpha}]$ is bounded on $L^{2}
(\R^{n})$, Theorem~\ref{ABG2} shows that  $H_{\alpha ,0}$ is in
$C^{1} (A_{\alpha})$.  
\end{proof}

\begin{lem}  \label{xac1}
$H_{\alpha ,0}$ is in $C^{2} (A_{\alpha})$ and $[ [H_{\alpha ,0},
A_{\alpha}] , A_{\alpha}]$ is bounded on $L^{2} (\R^{n})$.  
\end{lem}

\begin{proof}
As in Lemma~\ref{doublecom}, it is enough to estimate $[ [H_{\alpha
,0}, A_{\alpha}] , A_{\alpha}]$. Since $[H_{\alpha ,0}, A_{\alpha}]
\in \Psi (1,g_{2})$ and $A_{\alpha} \in \Psi \left( \< x\>^{1-\alpha}
\< \xi \> ,g_{2} \right)$, we get $[ [H_{\alpha ,0}, A_{\alpha}] ,
A_{\alpha}] \in \Psi \left( \< x\>^{-\alpha} ,g_{2} \right)$ which
implies the lemma. 
\end{proof}

\subsubsection{Regularity for $H_{\alpha}$}

\begin{prop}  \label{reg}
Assume that $V_\a$ satisfies the assumptions of Theorem~\ref{comple}. Then
$H_{\alpha}$ is of class $C^{1 + \delta} (A_{\alpha})$ for some
$\delta >0$. Moreover $[H_{\alpha} , A_{\alpha}]$ is bounded from
$D(H_{\alpha})$ to $L^{2}(\R^{n})$. 
\end{prop}

\begin{proof}
We use an interpolation argument as in
\cite[Proposition~3.7.5]{GeLa02}. We begin by proving that 
$H_{\alpha}$ is in $C^{1} (A_{\alpha})$ if $V_{\alpha}^{2}$ satisfies
\eqref{va2} with $\varepsilon \geq 0$. Since $H_{\alpha ,0}$ is $C^{1}
(A_{\alpha})$ and $[ H_{\alpha} , A_{\alpha}]$ is bounded from
$D(H_{\alpha})$ to $L^{2}(\R^{n})$, $e^{it A_{\alpha}}$ preserves
$D(H_{\alpha ,0}) = D (H_{\alpha})$, from Lemma~\ref{prm}. Then, from
Theorem~\ref{ABG2}, it is enough to show that $[H_{\alpha} ,
A_{\alpha}]$ is bounded from $D(H_{\alpha})$ to $L^{2}(\R^{n})$.

Since we know from Lemma~\ref{comutH0} and Lemma \ref{xac0} that
$[H_{\alpha ,0} , A_{\alpha}]$ is bounded from $D(H_{\alpha})$ to
$L^{2}(\R^{n})$, it is enough to show that: 
\begin{equation}  \label{cov}
[ V_{\alpha} , A_{\alpha} ] (H_{\alpha ,0} +i)^{-1} \text{ is compact
(resp. continuous) on } L^{2}(\R^{n}) \text{ if } \varepsilon >0
\text{ (resp. } \varepsilon =0 \text{)}, 
\end{equation}
where $\varepsilon$ is the constant in \eqref{va2}. We can write:
\begin{equation}  \label{ttt}
[ V_{\alpha} , A_{\alpha} ] (H_{\alpha ,0} +i)^{-1} =  V_{\alpha}^{1}
A_{\alpha} (H_{\alpha ,0} +i)^{-1} - A_{\alpha} V_{\alpha}^{1}
(H_{\alpha ,0} +i)^{-1} + [ V_{\alpha}^{2} , A_{\alpha} ] (H_{\alpha
,0} +i)^{-1}. 
\end{equation}
Let $\chi \in C^{\infty}_{0} (\R^{n})$ be equal to $1$ near the
support of $V_{\alpha}^{1}$. Since $A_{\alpha} \in \Psi \left( \<
p_{\alpha} (x) \> , g_{0} \right)$, $A_{\alpha} \chi$ and $\chi
A_{\alpha}$ are bounded. So 
\begin{equation}  \label{tx4}
A_{\alpha} V_{\alpha}^{1} (H_{\alpha ,0} +i)^{-1} = A_{\alpha} \chi
V_{\alpha}^{1} (H_{\alpha ,0} +i)^{-1} = {\mathcal O}(1)
V_{\alpha}^{1} (H_{\alpha ,0} +i)^{-1}, 
\end{equation}
which is compact because $V_{\alpha}^{1}$ is $H_{\alpha
,0}$--compact. Since $\chi A_{\alpha}$ is bounded, we can write:
\begin{align}
\chi A_{\alpha} (H_{\alpha ,0} +i)^{-1} =& (H_{\alpha ,0} +i)^{-1}
\chi A_{\alpha} + (H_{\alpha ,0} +i)^{-1} \left[ H_{\alpha ,0} , \chi
A_{\alpha} \right] (H_{\alpha ,0} +i)^{-1}    \nonumber  \\ 
=&(H_{\alpha ,0} +i)^{-1} {\mathcal O} (1) + (H_{\alpha ,0} +i)^{-1}
\left[ H_{\alpha ,0} , \chi \right] A_{\alpha} (H_{\alpha ,0} +i)^{-1}
\nonumber   \\ 
&+ (H_{\alpha ,0} +i)^{-1} \chi \left[ H_{\alpha ,0} , A_{\alpha}
\right] (H_{\alpha ,0} +i)^{-1}.  \label{tx1}  
\end{align}
We have
\begin{align*}
{\rm Op} \left( c (x,\xi ) \right) := \left[ H_{\alpha ,0} , \chi \right] A_{\alpha} \in \Psi \left( \< \xi \> \< x \>^{-\infty}, g_{0} \right) .
\end{align*}
Then Proposition \ref{pse} and the pseudo-differential calculus imply
\begin{align*}
\left[ H_{\alpha ,0} , \chi \right] A_{\alpha} (H_{\alpha ,0} +i)^{-1} = \left( {\rm Op} \left( c(x,\xi ) \psi \left( \frac{\xi ^{2} -\< x \>^{\alpha}}{ \xi^{2} + \< x \>^{\alpha}} \right)  \right) + R \right) {\mathcal O} (1), 
\end{align*}
with $R \in \Psi \left( \< x \>^{-\infty} , g_{0} \right)$. Since we also have $c(x,\xi ) \psi \left( \frac{\xi ^{2} -\< x \>^{\alpha}}{ \xi^{2} + \< x \>^{\alpha}} \right) \in S \left( \< x \>^{-\infty} , g_{0} \right)$, we get
\begin{equation}
\left[ H_{\alpha ,0} , \chi \right] A_{\alpha} (H_{\alpha ,0} +i)^{-1} = {\mathcal O} (1).
\end{equation}
As $\left[ H_{\alpha ,0} , A_{\alpha} \right]
(H_{\alpha ,0} -z)^{-1}$ is bounded, \eqref{tx1} becomes 
\begin{align*}
\chi A_{\alpha} (H_{\alpha ,0} +i)^{-1} = (H_{\alpha ,0} +i)^{-1}
{\mathcal O} (1), 
\end{align*}
and then
\begin{equation}   \label{tx3}
V_{\alpha}^{1} A_{\alpha} (H_{\alpha ,0} +i)^{-1} = V_{\alpha}^{1}
(H_{\alpha ,0} +i)^{-1} {\mathcal O}(1), 
\end{equation}
which is compact. Since $A_{\alpha} \in \Psi \left( \< p_{\alpha} (x)
\> ,g_{0} \right)$ and $V_{\alpha}^{2}$ satisfies \eqref{va2}, we have 
\begin{align}
\left[ V_{\alpha}^{2} , A_{\alpha} \right] (H_{\alpha ,0} +i)^{-1} =&
\,\left( V_{\alpha}^{2} A_{\alpha} - A_{\alpha} V_{\alpha}^{2} \right)
(H_{\alpha ,0} +i)^{-1}  \nonumber \\ 
=&\, {\mathcal O} (1) \< p_{\alpha} (x) \>^{-\varepsilon} (H_{\alpha ,0}
+i)^{-1},  \label{tx2} 
\end{align}
which is compact (resp. bounded) if $\varepsilon >0$
(resp. $\varepsilon =0$) from Lemma~\ref{lem:compH}. So we have
\eqref{cov}, $H_{\alpha}$ is of class $C^{1} (A_{\alpha})$ and
$[H_{\alpha} , A_{\alpha}]$ is bounded from $D(H_{\alpha})$ to
$L^{2}(\R^{n})$.

To have $H_{\alpha}$ in $C^{1+\delta} (A_{\alpha})$, it remains to show that
\begin{align*}
T(V_{\alpha}^{2}) := \left[ (H_{\alpha} +i)^{-1} , A_{\alpha} \right] =
(H_{\alpha} +i)^{-1} \left[ H_{\alpha} , A_{\alpha} \right]
(H_{\alpha} +i)^{-1}, 
\end{align*}
is of class $C^{\delta} (A_{\alpha})$. We use an interpolation
argument as in \cite{GeLa02}. For $\rho>0$, we set: 
\begin{align*}
S^{\rho}_{\alpha} = \{ W \in L^{\infty} (\R^{n};\R); \ |W (x)| \lesssim
\< p_{\alpha} (x) \>^{-\rho} \text{ a.e. } x \in \R^{n} \}\, . 
\end{align*}
Then $S^{\rho}_{\alpha} $ is a Banach space, equipped with the norm $\| W
\|_{\rho ,\alpha} = \| \< p_{\alpha} (x) \>^{\rho} W(x) \|_{L^{\infty}
(\R^{n})}$. We already have proved that $T(\cdot)$ maps $S^{1}_{\alpha}$
into $C^{0} (A_{\alpha})$. From \eqref{tx2}, we get 
\begin{align*}
\big\| T(W) - T(\widetilde{W}) \big\| \lesssim \big\| W -
\widetilde{W} \big\|_{1, \alpha}, 
\end{align*} 
and then $T$ is continuous. We now show that $T(V_{\alpha}^{2})$ is of
class $C^{1} (A_{\alpha})$ for $V_{\alpha}^{2} \in
S^{2}_{\alpha}$. Using $H_{\alpha} \in C^{1} (A_{\alpha})$,
$[H_{\alpha} , A_{\alpha} ] (H_{\alpha} +i)^{-1} = {\mathcal O} (1)$
and Lemma 6.2.9 of \cite{AmBoGe90}, it is enough to show that $[
[H_{\alpha} , A_{\alpha} ] , A_{\alpha} ] (H_{\alpha} +i)^{-1}$ is
bounded. From Lemmas~\ref{doublecom} and \ref{xac1}, it is enough to
show that $[ [V_{\alpha} , A_{\alpha} ] , A_{\alpha} ] (H_{\alpha}
+i)^{-1}$ is bounded. As for \eqref{tx4} and \eqref{tx3}, we have 
\begin{align*}
[[V_{\alpha}^{1} , A_{\alpha} ] , A_{\alpha} ] (H_{\alpha} +i)^{-1} =&
\left( V_{\alpha}^{1} A_{\alpha}^{2} - 2 A_{\alpha} V_{\alpha}^{1}
A_{\alpha} + A_{\alpha}^{2} V_{\alpha}^{1} \right) (H_{\alpha}
+i)^{-1}   \\ 
=& V_{\alpha}^{1} (H_{\alpha} +i)^{-1} {\mathcal O}(1) - 2 {\mathcal
O}(1)  V_{\alpha}^{1} (H_{\alpha} +i)^{-1} {\mathcal O}(1) + {\mathcal
O}(1) V_{\alpha}^{1} (H_{\alpha} +i)^{-1}, 
\end{align*}
which is bounded. On the other hand,
\begin{align*}
[[V_{\alpha}^{2} , A_{\alpha} ] , A_{\alpha} ] (H_{\alpha} +i)^{-1} =&
\left( V_{\alpha}^{2} A_{\alpha}^{2} - 2 A_{\alpha} V_{\alpha}^{2}
A_{\alpha} + A_{\alpha}^{2} V_{\alpha}^{2} \right) (H_{\alpha}
+i)^{-1}, 
\end{align*}
is bounded since $V_{\alpha}^{2} \in S^{2}_{\alpha}$ and $A_{\alpha}
\in \Psi \left( \< p_{\alpha} (x) \> , g_{0} \right)$. Then
$T(V_{\alpha}^{2}) \in C^{1} (A_{\alpha})$ for $V_{\alpha}^{2} \in
S_{\alpha}^{2}$.

Moreover, we know that for $0< \varepsilon' <1$, $C^{\varepsilon '}
(A_{\alpha})$ is a real interpolation space between $C^{0}
(A_{\alpha})$ and $C^{1} (A_{\alpha})$. Using the notation of
\cite{AmBoGe90}, \cite[Equation (5.2.22)]{AmBoGe90} implies 
\begin{align*}
C^{\varepsilon '} (A_{\alpha}) = \left( C^{0} (A_{\alpha}) , C^{1}
(A_{\alpha}) \right)_{1-\varepsilon ' , \infty}. 
\end{align*}
On the other hand, mimicking the proof of Lemma A.3 of \cite{GeLa02}
with $\chi_{R} = \chi ( p_{\alpha} (x) /R )$, we prove that for $\rho
\in ]1,2[$, 
\begin{align*}
S_{\alpha}^{\rho} \subset \left( S_{\alpha}^{1} , S_{\alpha}^{2}
\right)_{\rho -1, \infty}. 
\end{align*}
By interpolation (see e.g. \cite[Theorem 2.6.1]{AmBoGe90}),
there is $\delta >0$ such that $H_{\alpha}$ is of class $C^{1+\delta}
(A_{\alpha})$ for $V_{\alpha} \in S_{\alpha}^{1+\varepsilon}$ with
$\varepsilon >0$. 
\end{proof}

\subsection{Mourre estimates}  \label{moumou}

First, we prove a Mourre estimate for $H_{2,0}$.

\begin{lem} \label{m20}
Let $\eta >0$ and $\chi \in C^{\infty}_{0} (\R)$. There exists 
 $K$ compact on $L^{2} (\R^{n})$ such that 
\begin{equation}
\chi (H_{2,0}) \left[ i H_{2,0} , A_{2} \right] \chi (H_{2,0}) \geq
(2-\eta ) \chi^{2} (H_{2,0}) + \chi (H_{2,0}) K \chi (H_{2,0}). 
\end{equation}
\end{lem}

\begin{proof}
Using \eqref{com1}, we have
\begin{align}
\chi (H_{2,0}) \left[ i H_{2,0} , A \right] \chi (H_{2,0})=& 2 \chi (H_{2,0})
U \left( \frac{2 x^{2}}{\< \sqrt{2} x \>^{2}} + \frac{2 D^{2}}{\<
\sqrt{2} D \>^{2}} \right) U^{*} \chi (H_{2,0}). 
\end{align}
For $x^{2}+ \xi^{2} > C$ with $C \gg1$, the symbol
\begin{equation*}
f(x,\xi) = \frac{2 x^{2}}{\< \sqrt{2} x \>^{2}} + \frac{2 D^{2}}{\<
\sqrt{2} D \>^{2}} \in S(1,g_{2}) 
\end{equation*}
satisfies $f \geq 1- \eta /2$. Then G\aa rding inequality
(Theorem~18.6.7 of \cite{Ho85}) yields 
\begin{equation*}
{\rm Op} (f) \geq (1- \eta /2) - \widetilde{C} {\rm Op} \big(
\chi (x,\xi) \big) - R , 
\end{equation*}
with $\widetilde{C} >0$, $\chi \in C_{0}^{\infty} (\R^{2n})$ and $R
\in \Psi \left( \< x \>^{-1} \< \xi \>^{-1} , g_{2} \right)$. Then 
\begin{equation*}
\chi (H_{2,0}) \left[ i H_{2,0} , A \right] \chi (H_{2,0}) \geq (2 -
\eta) \chi (H_{2,0})^{2} + \chi (H_{2,0}) K \chi (H_{2,0}), 
\end{equation*}
where $K$ is compact.
\end{proof}

We have also a Mourre estimate for $H_{\alpha ,0}$ with $0< \alpha <2$.

\begin{lem} \label{ma0}
Let $\eta >0$. If the support of $\psi$ in \eqref{nji} is close enough
to $0$ and $\chi \in C^{\infty}_{0} (\R)$, there exists a compact operator
$K$ on $L^{2} (\R^{n})$ such that 
\begin{equation}
\chi (H_{\alpha,0}) \left[ i H_{\alpha,0} , A_{\alpha} \right] \chi
(H_{\alpha,0}) \geq (2 -\alpha -\eta ) \chi^{2} (H_{\alpha,0}) + \chi
(H_{\alpha,0}) K \chi (H_{\alpha,0}). 
\end{equation}
\end{lem}

\begin{proof}
Since $a_{\alpha} (x,\xi ) \in S \left( \< x \>^{1-\alpha} \< \xi \>
,g_{2} \right)$, $\xi^{2} - \< x \>^{\alpha} \in S \left( \xi^{2} + \<
x \>^{\alpha} ,g_{2} \right)$ and $\< \xi \>$ is like $\< x \>^{\alpha
/2}$ on the support of $a_{\alpha}$, we have  
\begin{align*}
[iH_{\alpha ,0} ,A_{\alpha} ] = {\rm Op} (b_{1}) + {\rm Op} (b_{2}) +
K_{1}, 
\end{align*}
where
\begin{align*}
b_{1} (x,\xi ) =& \big( 2\xi^{2} \< x\>^{-\alpha} - 2 \alpha (x.\xi
)^{2} \< x \>^{-\alpha -2} + \alpha x^{2} \< x \>^{-2} \big) \psi
\left( \frac{\xi^{2} -\< x \>^{\alpha}}{\xi^{2} + \< x \>^{\alpha}}
\right)  \in S \left( 1,g_{2} \right), 
\end{align*}
$b_{2} (x,\xi ) \in S (1, g_{2})$ with support inside the support of
$\psi ' \Big( \frac{\xi^{2} -\< x \>^{\alpha}}{\xi^{2} + \< x
\>^{\alpha}} \Big)$ and $K_{1} \in \Psi \left( \< x \>^{-1} \< \xi
\>^{-1} , g_{2} \right)$. If the support of $\psi$ is close enough to
$0$, we have  
\begin{align*}
b_{1} (x,\xi) \geq (2-\alpha - \eta ) \psi \left( \frac{\xi^{2} -\< x
\>^{\alpha}}{\xi^{2} + \< x \>^{\alpha}} \right),  
\end{align*}
for $(x,\xi)$ large enough, and the G\aa rding inequality implies
\begin{align*}
{\rm Op} (b_{1}) \geq (2- \alpha - \eta ) {\rm Op} \left(
\psi \left( \frac{\xi^{2} -\< x \>^{\alpha}}{\xi^{2} + \< x
\>^{\alpha}} \right) \right) +K_{2}, 
\end{align*}
with $K_{2} \in \Psi \left( \< x \>^{-1} \< \xi \>^{-1} ,g_{2}
\right)$. Therefore  
\begin{equation}   \label{tx6}
\begin{aligned}
\chi (H_{\alpha,0}) [iH_{\alpha ,0} ,A_{\alpha} ] \chi (H_{\alpha,0})
\geq&\, (2-\alpha - \eta ) \chi (H_{\alpha,0})  {\rm Op} \left( \psi
\left( \frac{\xi^{2} -\< x \>^{\alpha}}{\xi^{2} + \< x \>^{\alpha}}
\right) \right) \chi (H_{\alpha,0})   \\  
&+ \chi (H_{\alpha,0}) \left( {\rm Op} (b_{2}) + K_{1} + K_{2} \right)
\chi (H_{\alpha,0}) . 
\end{aligned}
\end{equation}
Let $\widetilde{\chi} \in C_{0}^{\infty} (\R)$ be equal to $1$ near the support of $\chi$. Using Proposition (\ref{pse}), we get
\begin{align*}
{\rm Op} \left( \psi \left( \frac{\xi^{2} -\< x \>^{\alpha}}{\xi^{2} + \< x \>^{\alpha}} \right) \right) \chi (H_{\alpha,0}) =& {\rm Op} \left( \psi \left( \frac{\xi^{2} -\< x \>^{\alpha}}{\xi^{2} + \< x \>^{\alpha}} \right) \right) (H_{\alpha,0} +i)^{-1} (H_{\alpha,0} +i) \widetilde{\chi} (H_{\alpha,0}) \chi (H_{\alpha,0})  \\
=& \chi (H_{\alpha,0}) - {\mathcal O} (1) {\rm Op} \left( r \right) (H_{\alpha,0} +i) \widetilde{\chi} (H_{\alpha,0}) \chi (H_{\alpha,0}) \, .
\end{align*}
Then we have
\begin{align*}
\chi (H_{\alpha,0})  {\rm Op} \left( \psi \left( \frac{\xi^{2} -\< x
\>^{\alpha}}{\xi^{2} + \< x \>^{\alpha}} \right) \right) \chi
(H_{\alpha,0}) = \chi^{2} (H_{\alpha,0}) + \chi (H_{\alpha,0}) K_{3}
\chi (H_{\alpha,0})\,, 
\end{align*}
where $K_{3}$ is compact. Let $\widetilde{\psi} \in C^{\infty}_{0} ([-1/2 , 1/2])$ such that $\widetilde{\psi} = 1$ near $0$ and $\psi =1$ on the support of $\widetilde{\psi}$. Using Proposition~\ref{pse} with $\widetilde{\psi}$, we get: 
\begin{align*}
{\rm Op} (b_{2}) (H_{\alpha ,0}+i)^{-1} =& {\rm Op} \left( b_{2} (x,\xi) \widetilde{\psi} \left( \frac{\xi^{2} -\< x \>^{\alpha}}{\xi^{2} + \< x \>^{\alpha}} \right) \right) {\mathcal O} (1) + {\rm Op} (s) {\mathcal O} (1)  \\
=& 0 + {\rm Op} (s)  {\mathcal O} (1),
\end{align*}
with $s(x,\xi) \in S \left( \< x\>^{-1} \< \xi \>^{-1}, g_{0}  \right)$. Here we have used the fact that $b_{2}=0$ on the support of $\widetilde{\psi} \left( \frac{\xi^{2} -\< x \>^{\alpha}}{\xi^{2} + \< x \>^{\alpha}} \right)$. Then
\begin{align*}
\chi (H_{\alpha,0}) {\rm Op} (b_{2}) \chi (H_{\alpha,0} ) = \chi
(H_{\alpha,0}) K_{4} \chi (H_{\alpha,0} )\,, 
\end{align*}
with $K_{4}$ compact. Then \eqref{tx6} becomes
\begin{align*}
\chi (H_{\alpha,0}) [iH_{\alpha ,0} ,A_{\alpha} ] \chi (H_{\alpha,0})
\geq\, & (2-\alpha - \eta ) \chi^{2} (H_{\alpha,0}) \\
&+ \chi (H_{\alpha,0})
\left( K_{1} + K_{2} + (2-\alpha - \eta ) K_{3} + K_{4} \right) \chi
(H_{\alpha,0}), 
\end{align*}
which implies the lemma.
\end{proof}

Finally, we obtain a Mourre estimate for $H_{\alpha}$ for $0< \alpha
\leq 2$. 

\begin{prop}  \label{maa}
Let $\eta >0$ and $0< \alpha \leq 2$. If the support of $\psi$ in
\eqref{nji} is close enough to $0$ and $\chi \in C^{\infty}_{0} (\R)$,
there exists a compact operator $K$ on $L^{2} (\R^{n})$ such that 
\begin{equation}
\chi (H_{\alpha}) \left[ i H_{\alpha} , A_{\alpha} \right] \chi
(H_{\alpha}) \geq ( \sigma_{\alpha} -\eta ) \chi^{2} (H_{\alpha}) +
\chi (H_{\alpha}) K \chi (H_{\alpha}). 
\end{equation}
\end{prop}

\begin{proof}
Let $\widetilde{\chi} \in C^{\infty}_{0} (\R)$ with $\widetilde{\chi}
=1$ near the support of $\chi$. $\widetilde{\chi} (H_{\alpha}) -
\widetilde{\chi} (H_{\alpha ,0})$ is compact because $V_{\alpha}$ is
$H_{\alpha ,0}$--compact. Since $[H_{\alpha} , A_{\alpha} ]$ is
$H_{\alpha}$--bounded from Proposition~\ref{reg}, 
\begin{equation*}   
\chi (H_{\alpha}) \left[ i H_{\alpha} , A_{\alpha} \right] \chi
(H_{\alpha}) = \chi (H_{\alpha}) \widetilde{\chi} (H_{\alpha ,0})
\left( \left[ i H_{\alpha ,0} , A_{\alpha} \right] + i \left[
V_{\alpha} , A_{\alpha} \right] \right) \widetilde{\chi} (H_{\alpha
,0}) \chi (H_{\alpha}) + \chi (H_{\alpha}) K \chi (H_{\alpha}), 
\end{equation*}
with $K$ compact. From \eqref{cov}, Lemmas~\ref{m20} and \ref{ma0}, we have
\begin{align*}
\chi (H_{\alpha}) \left[ i H_{\alpha} , A_{\alpha} \right] \chi
(H_{\alpha}) \geq \sigma_{\alpha} \chi (H_{\alpha})
\widetilde{\chi}^{2} (H_{\alpha ,0}) \chi (H_{\alpha}) + \chi
(H_{\alpha}) K \chi (H_{\alpha}), 
\end{align*}
with another $K$ compact. Therefore,
\begin{align*}
\chi (H_{\alpha}) \left[ i H_{\alpha} , A_{\alpha} \right] \chi
(H_{\alpha}) \geq \sigma_{\alpha} \chi^{2} (H_{\alpha}) + \chi
(H_{\alpha}) K \chi (H_{\alpha}), 
\end{align*}
which implies the lemma.
\end{proof}

\section{Asymptotic completeness}
\label{sec:complet}

\subsection{Limiting absorption principle}

From Proposition~\ref{reg}, $H_{\alpha}$ is of class $C^{1+\delta}
(A_{\alpha})$, and $[H_{\alpha} , A_{\alpha}]$ maps $D(H_{\alpha})$
into $L^{2} (\R^{n})$. Using Proposition~\ref{maa},
Theorem 1.1 and 4.13 of \cite{BoGeMa93} yield:

\begin{theo}[Limiting Absorption Principle]  \label{lap}
Let  $0< \alpha \leq 2$.
The singular continuous spectrum of $H_{\alpha}$ is empty, and his
point spectrum is locally finite. For $\Lambda \subset \R \backslash
\sigma_{pp} (H)$ and $\nu > 1/2$, we have, for some $\eta>0$,  
\begin{equation}
\sup_{z \in \Lambda + i[-\eta , \eta]} \left\| \< A_{\alpha} \>^{-\nu}
(H_{\alpha} - z )^{-1}\< A_{\alpha} \>^{-\nu} \right\| < \infty.  
\end{equation}
\end{theo}

As a corollary, we have:

\begin{prop}  \label{mou}
Let $0< \alpha \leq 2$ and $\eta >0$. Assume that the support of $\psi$ in
\eqref{nji} is close enough to $0$. For $\lambda \in \R \setminus
\sigma_{pp} (H_{\alpha})$, there is a real open interval $\Lambda$
containing $\lambda$ such that 
\begin{equation}
{\bf 1}_{\Lambda} (H_{\alpha}) \left[ i H_{\alpha} , A_{\alpha}
\right] {\bf 1}_{\Lambda} (H_{\alpha}) \geq ( \sigma_{\alpha} -\eta )
{\bf 1}_{\Lambda} (H_{\alpha}). 
\end{equation}
\end{prop}

\begin{proof}
Since the spectrum of $H_{\alpha}$ is absolutely continuous near
$\lambda$, we have  
\begin{equation*}
{\rm s-} \lim_{\delta \to 0} {\mathbf 1}_{[\lambda-\delta, \lambda +
\delta]} (H_{\alpha}) = 0. 
\end{equation*}
Using Proposition \ref{maa}, we infer that  $K {\mathbf 1}_{[\lambda-\delta,
\lambda + \delta]} (H_{2,0})$ goes to $0$ in norm when $\delta \to
0$, since $K$ is compact. So we can find $\Lambda$ such that  
\begin{equation*}
{\bf 1}_{\Lambda} (H_{\alpha}) \left[ i H_{\alpha} , A_{\alpha}
\right] {\bf 1}_{\Lambda} (H_{\alpha}) \geq ( \sigma_{\alpha} - 2 \eta
) {\bf 1}_{\Lambda} (H_{\alpha}). 
\end{equation*}
Since $\eta>0$ is arbitrary, this yields the proposition. 
\end{proof}

\subsection{Minimal velocity estimate}   \label{mimi}

The following proposition is a simplified version of Proposition~A.1
of \cite{GeNi98}, due to ideas of I.~M.~Sigal and A.~Soffer.

\begin{prop}[\cite{GeNi98}]  \label{prop:minvel}
Let $H$ and $A$ be two self-adjoint operators on a separable Hilbert
space $\mathcal{H}$. We suppose that   
\begin{enumerate}
\item[i)] $H$ is in $C^{1+\delta}(A)$ for some $\delta>0$.
\item [ii)] There exists an interval  $\Lambda$ such that: 
${\bf 1}_{\Lambda}(H)[H,iA] {\bf 1}_{\Lambda}(H) \geq c {\bf
1}_{\Lambda}(H)$, with $c>0$.
\end{enumerate}
Then for any $g\in C_0^{\infty} (]-\infty ,c[)$ and any $f \in
C_0^{\infty}(\Lambda)$, we have  
\begin{equation} 
\int_1^{+\infty} \Big\| g \Big( \frac A t \Big) e^{-itH} f (H)u
\Big\|^2 \frac{dt}{t} \lesssim \|u\|^2, 
\end{equation}
for $u \in {\mathcal H}$, and
\begin{equation}
s-\lim_{t\rightarrow +\infty} g \Big( \frac{A}{t} \Big) e^{-itH} f
(H)=0. 
\end{equation}
\end{prop}

From Propositions~\ref{reg}, \ref{mou} and \ref{prop:minvel}, we get:

\begin{prop}  \label{mim}
For any $g \in C_0^{\infty} (]- \infty, \sigma_{\alpha} [)$ and any $f
\in C_0^{\infty}( \R )$, we have  
\begin{equation} 
\int_1^{+\infty} \left\| g \( \frac{A_{\alpha}}{t} \)
e^{-itH_{\alpha}} f (H_{\alpha}) u \right\|^2 \frac {dt}{t} \lesssim
\|u\|^2, 
\end{equation} 
for $u \in L^{2} (\R^{n})$, and
\begin{equation}
s-\lim_{t\rightarrow +\infty} g \( \frac{A_{\alpha}}{t} \)
e^{-itH_{\alpha}} f (H_{\alpha}) = 0. 
\end{equation}
\end{prop}

We want to replace $A_{\alpha}$ by $p_{\alpha} (x)$ in
Proposition~\ref{mim}. For that, we use a slight modification of
Lemma~A.3 of C.~G\'erard and F.~Nier \cite{GeNi98}.

\begin{lem}[\cite{GeNi98}]  \label{GeNi3}
Let ${\mathcal A}$ and ${\mathcal B}$ be two self-adjoint operators on
a separable Hilbert space ${\mathcal H}$ such that, for each $\mu >0$,
we have 
\begin{gather}
D ({\mathcal B}) \subset D ({\mathcal A}) \text{ and } 1 \leq
{\mathcal B}\,,    \label{GN1}\\  
{\mathcal A} \leq (1+\mu ) {\mathcal B} + C_{\mu}\,,
\label{GN2} 
\end{gather}
with $C_{\mu} \geq 0$, and
\begin{equation}
[ {\mathcal A} , {\mathcal B} ] {\mathcal B}^{-1} \in {\mathcal L}
({\mathcal H})\,.    \label{GN3} 
\end{equation}
Then for each $\lambda \in \R$, let $\varphi \in C^{\infty} (\R)$ with $\supp (\varphi ) \subset ]-\infty , \lambda [$, $\varphi = 1$ near $- \infty$ and $\psi \in C^{\infty}
(\R )$ with $\supp ( \psi) \subset ]\lambda , + \infty[$, $\psi = 1$ near $+ \infty$. We have  
\begin{equation}
\left\| \varphi \left( {\mathcal B} / t \right) \psi \left(
{\mathcal A} / t \right) \right\| = {\mathcal O}
\(t^{-1}\)\quad\text{as} \quad t\to +\infty\,.
\label{GN4} 
\end{equation}
\end{lem}

\begin{proof} We follow the proof of \cite[Lemma~A.3]{GeNi98}. Let
$\varphi_{1} \in C^{\infty} ( \R )$ with $\supp ( \varphi_{1} ) \subset]-\infty , \lambda[$ and $\varphi_{1} =1$ near the support of $\varphi$. We have by \eqref{GN2},  
\begin{equation*}
\varphi_{1} ({\mathcal B} /t) {\mathcal A} \varphi_{1} ({\mathcal B}
/t) \leq (1+\mu ) \varphi_{1} ({\mathcal B} /t) {\mathcal B}
\varphi_{1} ({\mathcal B} /t) + {\mathcal O} (1)    
\leq (1+\mu ) \lambda t + {\mathcal O} (1).
\end{equation*}
So, if $\mu$ is small enough and $t$ large enough, we get
\begin{equation*}
\psi \left( \frac{\varphi_{1} ({\mathcal B} /t) {\mathcal A}
\varphi_{1} ({\mathcal B} /t)}{t} \right) = 0, 
\end{equation*}
and it remains to show that
\begin{equation}
\phi ({\mathcal B}/t) \left( \psi ({\mathcal A} / t) - \psi \left(
\frac{\varphi_{1} ({\mathcal B} /t) {\mathcal A} \varphi_{1} 
({\mathcal B} /t)}{t} \right) \right) = {\mathcal O} (t^{-1}).
\label{GN5} 
\end{equation}
The proof of \eqref{GN5} is the same as in \cite{GeNi98}.
\end{proof}

To apply the above result, we distinguish the cases $\a=2$ and $0<\a<2$.

\begin{lem}    \label{rem2}
The pairs $( {\mathcal A} , {\mathcal B}) = (A_{2} , \<
\ln \< x \> \>)$ and $( {\mathcal A} , {\mathcal B}) = (-A_{2} , \<
\ln \< x \> \>)$  satisfy the assumptions of Lemma~\ref{GeNi3},
provided that the
support of $\psi$ in \eqref{nnn} is small enough according to $\mu$. 
\end{lem}

\begin{proof}
We prove the lemma only for $( {\mathcal A} , {\mathcal B}) = (A_{2} ,
\< \ln \< x \> \>)$; the proof is the same in the other case. Since
${\mathcal A} \in \Psi \left( \< \ln \< x \> \> , g_{0} \right)$,
${\mathcal A}$ is well-defined and symmetric on $D ({\mathcal B}) = \{
u \in L^{2} (\R^{n}) ; \ \< \ln \< x \> \> u \in L^{2} (\R^{n})
\}$. So the assumption \eqref{GN1} of Lemma~\ref{GeNi3} is true from
Theorem \ref{nel}. Moreover ${\mathcal B} \in \Psi \left( \< \ln \< x
\> \> , g_{1} \right)$ implies $[ {\mathcal A}, {\mathcal B}]
{\mathcal B}^{-1} \in \Psi \big( \< \ln \< x \> \> \< x \>^{-1}, g_{0}
\big)$, from  \cite[Theorem 18.5.5]{Ho85}. Then the assumption
(\ref{GN3}) is also true.

Let $f$, $g \in C_{0}^{\infty} ([0,1]; [0,1])$ be equal to $1$ near
$0$. For $\delta$, $M >0$, we can write  
\begin{equation}
{\mathcal A} = {\rm Op} (s_{1}) + {\rm Op} (s_{2}) + {\rm Op} (s_{3}),
\end{equation} 
with
\begin{align*}
s_{1} &= \( \ln \< \xi + x\> f \( \frac{\< \xi + x \>}{\< x
\>^{\delta}} \) - \ln \< \xi - x\> f \( \frac{\< \xi - x
\>}{\< x \>^{\delta}} \) \) g \( \< x \> / M \)
,  \\  
s_{2} &= \( \ln \< \xi + x\> f \( \frac{\< \xi + x \>}{\< x
\>^{\delta}} \) - \ln \< \xi - x\> f \( \frac{\< \xi - x
\>}{\< x \>^{\delta}} \) \) (1-g) \( \< x \> / M
\) ,  \\  
s_{3} &= \( \ln \< \xi + x\> (1-f) \( \frac{\< \xi + x \>}{\<
x \>^{\delta}} \) - \ln \< \xi - x\> (1-f) \( \frac{\< \xi - x
\>}{\< x \>^{\delta}} \) \) . 
\end{align*}
Then
\begin{equation}     \label{mp1}
|( {\rm Op} (s_{1} ) u,u)| \leq C \| u \|^{2} , 
\end{equation}
because $s_{1} \in S \left( \< x \>^{-\infty} , g_{0} \right)$. On the other hand,
since 
$$s_{2} \in S \left( \< \ln \< x \> \>, g_{0} \right)\quad\text{ and
}\quad
\<\ln \< x \> \>^{-1/2} \in S \( \< \ln \< x \> \>^{-1/2} , \frac{|dx|^{2}}{\< x \>^{2}} + \frac{|d \xi |^{2}}{\< x \>^{2}} \),$$
we get: 
\begin{equation}  \label{tx10}
\< \ln \< x \> \>^{-1/2} {\rm Op} (s_{2}) \< \ln \< x \> \>^{-1/2} =
{\rm Op} \big( \< \ln \< x \> \>^{-1} s_{2} (x,\xi) \big) + R, 
\end{equation}
with $R \in \Psi \left( \< x \>^{-1}, g_{0} \right)$. Since $\supp f
\subset [0,1]$, we get: 
\begin{align*}
\left| \< \ln \< x \> \>^{-1} s_{2} (x,\xi) \right| \leq 2 \delta \,.
\end{align*}
We also have
\begin{align*}
\big| \partial_{x}^{\alpha} \partial_{\xi}^{\beta} \< \ln \< x \> \>^{-1} s_{2} (x,\xi) \big|
\leq C_{\alpha ,\beta ,\delta} \ln (M)^{-1},
\end{align*} 
where $C_{\alpha ,\beta ,\delta}$ depends on $\delta$. Fix $\delta$
small enough, and then, $M$ large enough. Theorem~18.6.3 of \cite{Ho85}
yields:  
\begin{align*}
\big\| {\rm Op} \big( \< \ln \< x \> \>^{-1} s_{2} (x,\xi) \big)
\big\| < \eta /2. 
\end{align*}
Then \eqref{tx10} implies
\begin{align*}
\big( \< \ln \< x \> \>^{-1/2} {\rm Op} (s_{2}) \< \ln \< x \>
\>^{-1/2} u, u\big ) \leq \eta /2 \| u \|^{2} + ( R u,u), 
\end{align*}
and since $R \in \Psi \left( \< x \>^{-1}, g_{0} \right)$,
\begin{equation}   \label{mp2}
\big( {\rm Op} (s_{2}) u, u \big ) \leq \eta /2 ({\mathcal B} u,u) + C
\| u \|^{2}. 
\end{equation}
So it remains to study $s_{3} (x,\xi)$. Using (\ref{mm3}) and (\ref{mm4}), we get
\begin{align}
s_{3} (x, \xi ) \leq& \ln \< \xi + x \> - \ln \< \xi -x \> + \ln \< \xi - x\> f \( \frac{\< \xi - x \>}{\< x \>^{\delta}} \)    \nonumber  \\
\leq& (1+\delta) \< \ln \< x \> \> + C .
\end{align}
We also have
\begin{align*}
s_{3} (x, \xi ) \in S \left( \< \ln \< x \> \>, \frac{|dx|^{2}}{\< x \>^{2 \delta}} + \frac{|d \xi |^{2}}{\< x \>^{2 \delta}} \right).
\end{align*}
If we assume $\delta < \eta /2$, G\aa rding inequality implies: 
\begin{align*}
\( {\rm Op} (s_{3}) u, u \) \leq (1+\eta /2) ({\mathcal B} u,u)
+ C \| u \|^{2} + (Ru,u),  
\end{align*}
with $R \in \Psi \left( \< \ln \< x \> \> \< x \>^{-2 \delta } ,
g_{0} \right)$. So, 
\begin{equation} \label{mp3}
\( {\rm Op} (s_{3}) u, u \) \leq (1+\eta /2) ({\mathcal B} u,u)
+ C \| u \|^{2}. 
\end{equation}
Combining with \eqref{mp1}, \eqref{mp2} and \eqref{mp3}, we get 
\begin{align*}
( {\mathcal A} u,u)  \leq (1+\eta) ({\mathcal B} u ,u) + C \| u \|^{2},
\end{align*}
which is \eqref{GN2}.
\end{proof}

\begin{lem}  \label{rem1}
Let $0< \alpha <2$. The pairs of operators $( {\mathcal A} , {\mathcal
B}) = (A_{\alpha} , \< x \>^{1-\alpha /2})$ and $( {\mathcal A} ,
{\mathcal B}) = (-A_{\alpha} , \< x \>^{1-\alpha /2})$ satisfy the
assumptions of Lemma~\ref{GeNi3}, provided that the support of $\psi$ in
\eqref{nji} is small enough according to $\mu$.  
\end{lem}

\begin{proof}
We prove the lemma only for $( {\mathcal A} , {\mathcal B}) =
(A_{\alpha} , \< x \>^{1-\alpha /2})$; the proof is the same in the
other case. Since ${\mathcal A} \in \Psi \big( \< x \>^{1-\alpha /2} ,
g_{1}^{\alpha /2} \big)$, ${\mathcal A}$ is well-defined and symmetric
on $D ({\mathcal B}) = \{ u \in L^{2} (\R^{n}) ; \ \< x \>^{1-\alpha
/2 } u \in L^{2} (\R^{n}) \}$. So the hypothesis \eqref{GN1} of
Lemma~\ref{GeNi3} is true from Theorem~\ref{nel}. Moreover, ${\mathcal
B} \in \Psi \big( \< x \>^{1-\alpha /2} , g_{1} \big)$ implies $[ {\mathcal
A}, {\mathcal B}] {\mathcal B}^{-1} \in \Psi \big( \< x \>^{- \alpha}
, g_{1}^{\alpha /2} \big)$ from Theorem 18.5.5 of \cite{Ho85}. Then
the assumption (\ref{GN3}) is also true.

On the support of $a_{\alpha} (x,\xi)$, with $(x,\xi)$ large enough, we
have $| \xi | = \< x \>^{\alpha /2} (1+o(1))$, where $o(1)$ stands for an
arbitrary small function as $\supp \psi \rightarrow \{ 0 \}$. Then, 
\begin{align*}
a_{\alpha} (x,\xi ) \leq (1+\eta /2) \< x \>^{1-\alpha /2} + C,
\end{align*}
with $C>0$. The G\aa rding inequality in $\Psi \big( g_{1}^{\alpha /2}
\big)$ implies that
\begin{equation}     \label{gh1}
( {\mathcal A} u,u)  \leq (1+\eta /2) ({\mathcal B} u ,u) + C \| u \|
^{2} + \left( {\rm Op} (r) u,u \right),  
\end{equation}
with $r \in S \big( \< x \>^{1- 3 \alpha /2} , g_{1}^{\alpha /2}
\big)$. We have 
\begin{align*}
\left| \left( {\rm Op} (r) u,u \right) \right| =&\, \left| \( \< x
\>^{-1/2 + \alpha /4} {\rm Op} (r) \< x \>^{-1/2 + 5 \alpha /4} \< x
\>^{- \alpha } \< x 
\>^{1/2 - \alpha /4} u, \< x \>^{1/2 - \alpha /4} u \) \right|    \\   
\lesssim&\,  \left\| \< x \>^{- \alpha } \< x \>^{1/2 - \alpha /4} u \right\|
\left\| \< x \>^{1/2 - \alpha /4} u \right\|    \\ 
\leq& \, \( \eta /4 \left\| \< x \>^{1/2 - \alpha /4} u \right\| + C
\| u\| \) \left\| \< x \>^{1/2 - \alpha /4} u \right\|  \\  
\leq& \, \eta /2 \left\| \< x \>^{1/2 - \alpha /4} u \right\|^{2} + C
\| u\|^{2}.
\end{align*}
So, \eqref{gh1} becomes
\begin{equation}
( {\mathcal A} u,u)  \leq (1+\eta) ({\mathcal B} u ,u) + C \| u \| ^{2},
\end{equation}
which proves \eqref{GN2} and the lemma.
\end{proof}

From Proposition~\ref{mim}, Lemmas~\ref{GeNi3}, \ref{rem2} and
\ref{rem1}, we obtain:

\begin{prop}[Minimal Velocity Estimate] \label{P.P.6}
For any $\chi \in C_0^{\infty} (\R)$ with $\supp  \chi  \cap
\sigma_{pp} (H_{\alpha}) = \emptyset$, $0< \theta < \sigma_{\alpha}$,
and $u \in L^{2} (\R^{n})$ we have: 
\begin{gather*}
\int_1^{\infty} \left\| {\bf 1}_{[0, \theta ]} \(
\frac{p_{\alpha}(x)}{t} \) e^{itH_{\alpha}} \chi (H_{\alpha}) u
\right\|^2 \frac{dt}{t} \lesssim \| u\|^2,   \\ 
s-\lim_{t \rightarrow + \infty} {\bf 1}_{[ 0, \theta ]} \(
\frac{p_{\alpha}(x)}{t} \) e^{-itH_{\alpha}} {\bf 1}^{c} (H_{\alpha}) =0. 
\end{gather*}
\end{prop}

\subsection{Proof of Theorem \ref{comple}}
\label{subsec:demo}

As mentioned in the introduction, we prove Theorem \ref{comple}, with
 $\H_{\a,0}$ (resp. $\H_\a$)  replaced by $H_{\a,0}$ (resp. $H_\a$),
 since this substitution will turn out to yield a short range
 perturbation. 
We prove \eqref{W2}; it will be clear from the
proof that \eqref{W1} follows the same way.
By a density argument
and using that $\sigma_{pp} (H_{\alpha ,0}) = \emptyset$, and
$\sigma_{pp} (H_{\alpha})$ has no accumulating point, it is enough to
show the existence of  
\begin{align*}
s-\lim_{t\to + \infty} e^{itH_{\alpha ,0}} e^{-itH_{\alpha}} \chi^{2}
(H_{\alpha}),  
\end{align*}
with $\supp \chi  \cap \sigma_{pp}(H_{\alpha}) = \emptyset$. We have 
\begin{equation}    \label{rtt2}
e^{itH_{\alpha ,0}} e^{-itH_{\alpha}} \chi^{2} (H_{\alpha}) = \chi
(H_{\alpha ,0}) e^{itH_{\alpha ,0}} e^{-itH_{\alpha}} \chi
(H_{\alpha})  + e^{itH_{\alpha ,0}} \big( \chi (H_{\alpha ,0}) - \chi
(H_{\alpha}) \big) e^{-itH_{\alpha}} \chi (H_{\alpha}). 
\end{equation}
As the spectrum of $H_{\alpha}$ is absolutely continuous on $\supp (
\chi )$, $e^{-itH_{\alpha}} \chi (H_{\alpha}) \to 0$ weakly. Since
$\chi (H_{\alpha ,0}) - \chi (H_{\alpha})$ is compact, the second term
in \eqref{rtt2} converges strongly to $0$.

Let $g_{1} \in C^{\infty}_{0} (]- \infty , \sigma_{\alpha}
[)$ such that $g_{1} = 1$ near $0$. From Proposition \ref{P.P.6}, we
deduce that 
\begin{equation}  \label{tx11}
s-\lim_{t\rightarrow +\infty} \chi (H_{\alpha ,0}) e^{itH_{\alpha,0}}
g_{1} \bigg( \frac{p_{\alpha}(x)}{t} \bigg) e^{-itH_{\alpha}} \chi
(H_{\alpha}) =0.  
\end{equation}
Now, let us consider
\begin{align*}
G(t) = \chi (H_{\alpha ,0}) e^{itH_{\alpha ,0}} (1 - g_{1}) \bigg(
\frac{p_{\alpha}(x)}{t} \bigg) e^{-itH_{\alpha}} \chi (H_{\alpha}) u.  
\end{align*}
The function $G(t)$ is differentiable and
\begin{equation}
\begin{split}  \label{gro}
G'(t) =& \chi (H_{\alpha ,0}) e^{itH_{\alpha ,0}} \bigg[ g_{1} \bigg(
\frac{p_{\alpha}(x)}{t} \bigg) , iH_{\alpha ,0}
\bigg]e^{-itH_{\alpha}} \chi (H_{\alpha}) u \\ 
&+ \chi (H_{\alpha ,0}) e^{itH_{\alpha ,0}}
\frac{p_{\alpha}(x)}{t^{2}} g_{1}' \bigg( \frac{p_{\alpha}(x)}{t}
\bigg) 
e^{-itH_{\alpha}} \chi (H_{\alpha}) u  \\ 
&+ \chi (H_{\alpha ,0}) e^{itH_{\alpha ,0}} (1 - g_{1}) \bigg(
\frac{p_{\alpha}(x)}{t} \bigg) V_{\alpha} (x) e^{-itH_{\alpha}} \chi
(H_{\alpha})u \, .
\end{split}
\end{equation}
The first term of \eqref{gro} is equal to
\begin{align*}
I (t) = \frac{1}{t} \chi (H_{\alpha ,0}) e^{itH_{\alpha ,0}} \left(
{\rm Op} (f) g_{1}' \bigg( \frac{p_{\alpha}(x)}{t} \bigg) -i \frac{|
\nabla_x p_{\alpha} (x) |^{2}}{t} g_{1}'' \bigg(
\frac{p_{\alpha}(x)}{t} \bigg)\right) e^{-itH_{\alpha}} \chi
(H_{\alpha}) u , 
\end{align*}
with 
\begin{align*}
f(x,\xi) = -2 \nabla_x p_{\alpha}(x) \xi \in S \big( \< x \>^{-\alpha
/2} \< \xi \> , g_{2} \big) . 
\end{align*}
Using Proposition~\ref{pse} and the fact that $|\nabla_x
p_{\alpha}(x)|$ is bounded, we get 
\begin{align*}
I(t) =& \frac{1}{t} \chi (H_{\alpha ,0}) e^{itH_{\alpha ,0}} {\rm Op}
\bigg( \psi \bigg( \frac{\xi^{2} - \< x \>^{\alpha}}{\xi^{2} + \< x
\>^{\alpha}} \bigg) \bigg) {\rm Op} (f) g_{1}' \bigg(
\frac{p_{\alpha}(x)}{t} \bigg) e^{-itH_{\alpha}} \chi (H_{\alpha}) u
\\ 
&+ {\mathcal O}(t^{-1}) \bigg\| {\rm Op} (r) {\rm Op} (f) g_{1}'
\bigg( \frac{p_{\alpha}(x)}{t} \bigg) \bigg\| + {\mathcal O} (t^{-2}), 
\end{align*}
with $r \in S \left( \< \xi \>^{-2} , g_{0} \right)$. Using the
pseudo-differential calculus and the fact that $\< x \>$ is like
$t^{\frac{1}{1-\alpha /2}}$ (resp. $e^t$) if $0< \alpha < 2$
(resp. $\alpha =2$) on the support of $g_{1}' \left( p_{\alpha}(x)/ t
\right)$, we get 
\begin{equation*}
\bigg\| {\rm Op} (r) {\rm Op} (f) g_{1}' \bigg(
\frac{p_{\alpha}(x)}{t} \bigg) \bigg\|  
= \left\{
\begin{aligned}
\O& \big( t^{- \frac{\alpha /2}{1 - \alpha /2}} \big) &&\quad
\text{for } 0<\alpha<2 ,     \\ 
\O& \big( e^{-\alpha t/2} \big) &&\quad \text{for } \alpha =2 .
\end{aligned}  \right.
\end{equation*}
On the other hand, the pseudo-differential calculus in $\Psi (g_{1})$
implies 
\begin{align*}
{\rm Op} \bigg( \psi \bigg( \frac{\xi^{2} - \< x \>^{\alpha}}{\xi^{2}
+ \< x \>^{\alpha}} \bigg) \bigg) {\rm Op} (f) = {\rm Op} (m), 
\end{align*}
with $m(x,\xi ) \in S (1, g_{2})$. Let $g_{2} \in C^{\infty}_{0}
(]-\infty , \sigma_{\alpha} [)$ such that $g_{2} =0$ near $0$ and
$g_{2} =1$ near the support of $g_{1} '$. Using the
pseudo-differential calculus in $\Psi \big( | dx |^{2} \<
x\>^{-\alpha} + | d\xi |^{2} \big)$, we get 
\begin{align*}
{\rm Op} (m) g_{1}' \bigg( \frac{p_{\alpha} (x)}{t} \bigg) =& {\rm Op}
\bigg( m (x,\xi )  g_{1}' \bigg( \frac{p_{\alpha} (x)}{t} \bigg)
\bigg) g_{2} \bigg( \frac{p_{\alpha} (x)}{t} \bigg) + {\mathcal O} (1)
\< x \>^{-\alpha /2} g_{2} \bigg( \frac{p_{\alpha} (x)}{t}
\bigg)  \\ 
=& g_{2} \bigg( \frac{p_{\alpha} (x)}{t} \bigg)  {\mathcal O} (1)
g_{2} \bigg( \frac{p_{\alpha} (x)}{t} \bigg) + {\mathcal O} \( t^{-
\delta } \), 
\end{align*}
with $\delta >0$. Then
\begin{align*}
I(t) =& \frac{1}{t} \chi (H_{\alpha ,0}) e^{itH_{\alpha ,0}} g_{2}
\bigg( \frac{p_{\alpha} (x)}{t} \bigg) {\mathcal O} (1) g_{2} \bigg(
\frac{p_{\alpha} (x)}{t} \bigg) e^{-itH_{\alpha}} \chi (H_{\alpha}) u
+ {\mathcal O} \( t^{-1-\delta} \). 
\end{align*}
Proposition~\ref{P.P.6} and a duality argument imply that $I(t)$ is integrable.

The second term in \eqref{gro} can be written
\begin{equation}
\begin{split}     \label{gro2}
\chi ( H_{\alpha ,0}) e^{itH_{\alpha ,0}} \frac{p_{\alpha}(x)}{t^{2}}
& g_{1}' \bigg( \frac{p_{\alpha}(x)}{t} \bigg) e^{-itH_{\alpha}} \chi
(H_{\alpha}) u = \\ 
&= \frac{1}{t} \chi (H_{\alpha ,0}) e^{itH_{\alpha ,0}} g_{2} \bigg(
\frac{p_{\alpha}(x)}{t} \bigg) {\mathcal O} (1) g_{2} \bigg(
\frac{p_{\alpha}(x)}{t} \bigg)e^{-itH_{\alpha}} \chi (H_{\alpha}) u. 
\end{split}
\end{equation}
Like for $I(t)$, we get that \eqref{gro2} is integrable.

Finally, using assumptions \eqref{va1} and \eqref{va2}, we get, for
$t\gg 1$, 
\begin{equation*}
\bigg| (1 - g_{1}) \bigg( \frac{p_{\alpha}(x)}{t} \bigg) V_{\alpha}
(x) \bigg| = \bigg| (1 - g_{1}) \bigg( \frac{p_{\alpha}(x)}{t} \bigg)
V_{\alpha}^{2} (x) \bigg|  ={\mathcal O}\(t^{-1-\varepsilon}\) ,
\end{equation*}
which is integrable. This implies that the third term in \eqref{gro},
and then $G'(t)$, is integrable. So $G(t)$ has a limit when $t
\rightarrow + \infty$ and the theorem follows from \eqref{rtt2} and
\eqref{tx11}.

\qed  \def\H{{\mathbb H}}

\section{Asymptotic velocity}
\label{sP}

In this section, we construct the asymptotic velocity and describe
its spectrum. 
In \eqref{tx34}, we defined the position variable so that
it increases like $t$ along 
the evolution. We define the local velocity as: 
\begin{equation*}
\V_{\alpha}:=\left[iH_{\alpha},p_{\alpha}\(x\)\right].
\end{equation*}
(We use typewriter style letters to avoid any
confusion with previous notations.)
We denote $N=N_2$ the harmonic oscillator. 
The observable $\V_{\alpha}$ is defined as a quadratic form on
$D\(N\)$. By a direct calculation and an application of
Theorem~\ref{nel}, we obtain that $\V_{\alpha}$ is (well defined as an
operator 
and) essentially self-adjoint with domain $D(N)$; we note again
$\V_{\alpha}$ the self-adjoint extension. Thanks to
Theorem~\ref{comple}, it is sufficient to construct the asymptotic 
velocity and to describe its spectrum in the free case. Nevertheless,
the local velocity does not commute with the free evolution, in
particular the asymptotic velocity is different from the local
velocity even in the free case. In particular we have to establish propagation estimates even in the free case.
Before stating the main result of this section, we introduce some
notations. For an observable $\Theta\(t\)$, we denote ${\bf D}\Theta\(t\)$
its Heisenberg derivative with respect to $H_{\a,0}$, i.e. 
\begin{equation*}
{\bf D}\Theta\(t\):=\frac{d}{dt}\Theta\(t\)+{[}iH_{\a,0},\Theta\(t\){]}.
\end{equation*}
The main result we prove in this section is Theorem~\ref{Th.P.1}. 
It can be proved for all $0<\alpha\le 2$ in the same
way. Like for the asymptotic completeness, we give some
generalizations of this result in the case $\alpha=2$, see
Section~\ref{subsP.4}.

\subsection{Local velocity}
This section is devoted to the study of the local velocity and the local
acceleration. A direct calculation yields $2\nabla p_\a(x) = \sigma_\a
\<x\>^{-1-\a/2}x$, and: 
\begin{equation}
\label{P.1.1}
\V_{\alpha} = \frac{\sigma_{\alpha}}{2} \left( \frac{x}{\<x\>^{1+\alpha/2}}D + hc \right) ,
\end{equation}
where $hc$ stands for the adjoint of the first term. We set
\begin{equation}
\label{P.1.2}
{\v}_{\alpha}\(x,\xi\)=\sigma_{\alpha}\frac{x\cdot
\xi}{\<x\>^{1+\alpha/2}}\, .
\end{equation}
\begin{lem}  \label{L.P.1}
The operator $\(\V_{\alpha}, D\(N\)\)$ is well defined as an operator, and
essentially self-adjoint. We have $\V_{\alpha}\in
\Psi\(\<\xi\>\<x\>^{-\alpha/2},g_2\)$, and the symbol of $\V_{\alpha}$ is
${\v}_{\alpha}$. 
\end{lem}
\begin{proof}
We clearly have:
$\|\V_{\alpha}u\|\lesssim \|Nu\|$, for $u\in D\(N\)$.
We also have
\begin{equation*}
\{\xi^2+x^2,{\v}_{\alpha}\(x,\xi\)\}=
\sigma_{\alpha}\(\frac{2\xi^2}{\<x\>^{1+\alpha/2}}-\(2+\alpha\)
\frac{\(x\cdot\xi\)^2}{\<x\>^{3+\alpha/2}}-
2\frac{x^2}{\<x\>^{1+\alpha/2}}\),  
\end{equation*}
and the lemma follows from Theorem~\ref{nel}. 
\end{proof}

Define the acceleration:
\begin{align}
\tA_{\alpha}:=\left[iH_{\a,0},\V_{\alpha}\right] \quad ;\quad
\ta_{\alpha}\(x,\xi\):=\sigma_{\alpha}
\(\frac{2\xi^2}{\<x\>^{1+\alpha/2}}-
\(2+\alpha\)\frac{\(x\cdot\xi\)^2}{\<x\>^{3+\alpha/2}}+
\alpha\frac{x^2}{\<x\>^{3-\alpha/2}}\). \label{P.1.3}
\end{align}
The operator $\tA_{\alpha}$ is a pseudo-differential operator, with
principal symbol 
$$\ta_{\alpha}\(x,\xi\)\in
S\(\<\xi\>^2\<x\>^{-1+\alpha/2},g_2\).$$ 
We will often use the
decomposition $\ta_{\alpha} (x,\xi) = \ta_{\alpha}^1 (x,\xi) +\ta_{\alpha}^2 (x,\xi)$, with
$\ta_{\alpha}^2 \(x,\xi\) =\sigma_{\alpha} \alpha
\frac{x^2}{\<x\>^{3-\alpha/2}}\in S \(\<x\>^{-1+\alpha/2},g_2\)$, and
$\ta_{\alpha}^1 \(x,\xi\) \in S \(\<\xi\>^2\<x\>^{-1-\alpha/2},g_2\)$. 
\begin{lem}
\label{L.P.2}
The operators $\V_{\alpha}\(i+H_{\a,0}\)^{-1}$ and
$\tA_{\alpha}\(i+H_{\a,0}\)^{-1}$, defined on $D\(N\)$, can be
extended to bounded operators. 
\end{lem}
\begin{proof}
We prove a slightly more general result. Let $c\in
S\(\<\xi\>^m\<x\>^{-k},g_2\)$, with $\frac{\alpha m}{2}-k\le 0$, $0\le
m\le 2$. We prove:  
\begin{equation*}
\text{The operator ${\rm Op}\(c\)\(i+H_{\a,0}\)^{-1}$, defined on 
 $D\(N\)$, can be extended to a bounded operator.}
\end{equation*}
The lemma then follows using the decomposition  $\ta_{\alpha} = \ta_{\alpha}^1 + \ta_{\alpha}^2$.
Recall from Proposition~\ref{pse} that for all
$1\ge\beta>1/2$, 
\begin{equation}
\label{P.1.7}
 \(i+H_{\a,0}\)^{-1}={\rm
 Op}\(\psi\(
 \frac{\xi^2-\<x\>^{\alpha}}{\(\xi^2+\<x\>^{\alpha}\)^{\beta}}\)\)
 \(i+H_{\a,0}\)^{-1}+R_{\beta}, 
\end{equation}
with $N_{\alpha}^{\beta}R_{\beta}$ bounded, and $\psi\in
C_0^{\infty}(\R),\psi=1$ in a neighborhood of zero. Since ${\rm
Op}\(c\)N_{\alpha}^{-1}$ is bounded, it is sufficient to prove: 
\begin{equation}
\label{P.1.8}
{\rm Op}\(c\){\rm
Op}\(\psi\(\frac{\xi^2- \<x\>^{\alpha}}{\(\xi^2+
\<x\>^{\alpha}\)^{\beta}}\)\)\,\mbox{is 
bounded}. 
\end{equation}
This is a pseudo-differential operator, with principal symbol
\begin{equation*}
\left|c\(x,\xi\)\psi\(\frac{\xi^2- \<x\>^{\alpha}}{\(\xi^2+
 \<x\>^{\alpha}\)^{\beta}}\)\right|\lesssim\<x\>^{\frac{\alpha
 m}{2}-k}\psi\(\frac{\xi^2-
 \<x\>^{\alpha}}{\(\xi^2+\<x\>^{\alpha}\)^{\beta}}\)\lesssim 1 ,
\end{equation*}
where we have used $\<\xi\>\lesssim\<x\>^{\alpha/2}$ on
$\supp\psi\(\frac{\xi^2-\<x\>^{\alpha}}{\(\xi^2+\<x\>^{\alpha}\)^{\beta}}\)$.
This yields \eqref{P.1.8}, and the lemma.
\end{proof}
\begin{lem}
\label{L.P.3}
Let $f,\chi\in C_0^{\infty}(\R)$. Then, as $t\to\infty$:
\begin{itemize}
\item[(i)]$\displaystyle \left[\chi\(H_{\alpha}\),f\(
\frac{p_{\alpha}\(x\)}{t}\)\right]= {\mathcal O}\(t^{-1}\)$.
\item[(ii)] If $f$ is constant in a  neighborhood of $0$, then there
exists $\varepsilon >0$ such that:
\begin{equation*}
\left[f\(\frac{p_{\alpha}\(x\)}{t}\), \V_{\alpha}\right]=
\left\{
\begin{aligned}
&{\mathcal O}\(t^{-\frac{2+\alpha}{2-\alpha}}\) \quad  &&\text{if }\
0<\alpha<2,\\  
&{\mathcal O}\(e^{-\varepsilon t}\) \quad  &&\text{if }\ 
\alpha=2.
\end{aligned}
\right. 
\end{equation*}
\item[(iii)] If $f$ is constant in a  neighborhood of $0$, then there
exists $\varepsilon >0$ such that:
\begin{equation*}
\chi\(H_{\a,0}\)\left[f\(\frac{p_{\alpha}\(x\)}{t}\),
\tA_{\alpha}\right]\chi\(H_{\a,0}\)
=
\left\{
\begin{aligned}
&{\mathcal O}\(t^{\frac{-4}{2-\alpha}}\) \quad  &&\text{if }\ 0<\alpha<2,\\ 
&{\mathcal O}\(e^{-\varepsilon t}\)\quad &&\text{if }\ \alpha =2 .
\end{aligned}
\right. 
\end{equation*}
\end{itemize}
\end{lem}

\begin{proof}
$\(i\)$ Using Helffer--Sj\"ostrand formula, it is sufficient to show:
\begin{equation*}
\(z-H_{\alpha}\)^{-1}\frac{1}{2t}\(\V_{\alpha}f'
\(\frac{p_{\alpha}\(x\)}{t}\)+hc\)\(z-H_{\alpha}\)^{-1}={\mathcal
O}\(t^{-1}\),\quad \forall z\in \C\setminus\R .
\end{equation*}
The above relation follows from Lemma~\ref{L.P.2}, and the fact that
$D\(H_{\alpha}\)=D\(H_{\a,0}\)$. 

$(ii)$ We have:
\begin{equation*}
\left[i\V_{\alpha},f\(\frac{p_{\alpha}\(x\)}{t}\)\right]=
\frac{\sigma_{\alpha}^2}{2t}f'\(\frac{p_{\alpha}\(x\)}{t}\)
\frac{x^2}{\<x\>^{2+\alpha}}=
\left\{
\begin{aligned}
&{\mathcal O}\(t^{-\frac{2+\alpha}{2-\alpha}}\) \quad && \text{if }0<\alpha<2, \\
&{\mathcal O}\(e^{-\varepsilon t} \) \quad && \text{if }\alpha=2 ,
\end{aligned} \right.
\end{equation*}
for some $\varepsilon >0$.

$\(iii\)$ First notice that 
\begin{equation*}
\left[if\(\frac{p_{\alpha}\(x\)}{t}\),\tA_{\alpha}\right]=
\frac{1}{t}f'\(\frac{p_{\alpha}\(x\)}{t}\){\rm Op}\(c\) \, ,
\end{equation*}
with $c\in S\(\<\xi\>\<x\>^{-1-\alpha},g_2\)$. We now use
Proposition~\ref{pse}:  
\begin{equation*}
\chi\(H_{\a,0}\)={\rm Op}\(\psi\(\frac{\xi^2- \<x\>^{\alpha}}{\(\xi^2+
\<x\>^{\alpha}\)^{\beta}}\)\)\chi\(H_{\a,0}\)+ R_{\beta},\quad \forall
1\ge\beta>1/2, 
\end{equation*}
with $N_{\alpha}^{\beta}R_{\beta}$ bounded, and $\psi\in
C_0^{\infty}(\R),\psi=1$ in a small neighborhood of zero. Let
$0<\alpha<2$. We have: 
\begin{equation*}
{\rm Op}\(c\)N_{\alpha}^{-\beta}\in
S\(\<x\>^{-1-\alpha-\alpha\beta/2},g_0\), 
\end{equation*}
and thus
\begin{equation*}
\frac{1}{t}f'\( \frac{p_{\alpha}\(x\)}{t}\){\rm Op}\(c\)R_{\beta}= {\mathcal
O}\(t^{-\frac{4+ \alpha\beta+\alpha}{2-\alpha}}\).
\end{equation*}
Furthermore
\begin{equation*}
{\rm Op}\(c\){\rm Op}\(\psi\(\frac{\xi^2- \<x\>^{\alpha}}{\(\xi^2+
\<x\>^{\alpha}\)^{\beta}}\)\)\in  S\(\<x\>^{-1-\alpha/2},g_2\)
\end{equation*}
and thus
\begin{equation*}
\frac{1}{t}f'\( \frac{p_{\alpha}\(x\)}{t}\){\rm Op}\(c\){\rm Op}\(
\psi\(\frac{\xi^2-\<x\>^{\alpha}}{\(\xi^2+ \<x\>^{\alpha}\)^{\beta}}\)\)=
{\mathcal O}\(t^{-\frac{4}{2-\alpha}}\). 
\end{equation*}
For $\alpha=2$ all these terms are in ${\mathcal O}\(e^{-\varepsilon t}\)$
for some $\varepsilon>0$. This completes the proof of $\(iii\)$. 
\end{proof}

We now set $\displaystyle g_{\gamma,\eta}=\frac{|dx|^2}{\<x\>^{2\gamma}}+
\frac{|d\xi|^2}{\<\xi\>^{2\eta}}$.
\begin{lem}
\label{L.P.4}
Let $J\in C^{\infty}_b(\R),\,J=0$ on
$\left[-\varepsilon,\varepsilon\right]$ for some $\varepsilon>0$, and $a\in
S\(\<\xi\>^{\eta 
m},g_{\gamma,\eta}\)$ for some real $m$. Then for all $\delta$ with
$\varepsilon>\delta>0$, there exists $\widetilde{J}\in
C^{\infty}_b(\R)$, $\widetilde{J}=0$ on $\left[-\delta,\delta\right]$, with
$\widetilde{J}J=J$, such that:
\begin{equation*}
{\rm
Op}\(a\)J\(\frac{p_{\alpha}\(x\)}{t}\)=
\widetilde{J}\(\frac{p_{\alpha}\(x\)}{t}\){\rm
Op}\(a\)J\(\frac{p_{\alpha}\(x\)}{t}\)+
\left\{\begin{aligned} 
&{\mathcal
O}\(t^{-\infty}\) \quad && \text{if }0<\alpha<2, \\ 
&{\mathcal O}\(e^{-\infty t}\) \quad && \text{if }\alpha=2. 
\end{aligned} \right. 
\end{equation*}
\end{lem}

\begin{proof}
Let $\varepsilon>\delta'>\delta$, $\widetilde{J}\in
C^{\infty}_b(\R)$, with
$\supp\widetilde{J}\subset\R\setminus\left[-\delta,\delta\right]$, 
and $\widetilde{J}=1$ on $\R\setminus\left[-\delta',\delta'\right]$. Let
$\hat{J}:=1-\widetilde{J}.$ We have to estimate: 
\begin{equation*}
R(t,x):=
\hat{J}\(\frac{p_{\alpha}\(x\)}{t}\)\({\rm Op}\(a\)J\(\frac{p_{\alpha}\(\cdot
\)}{t}\)\Phi\)\(x\)
\end{equation*}
By definition,
\begin{equation*}
R(t,x)= 
\iint e^{i(x-y)\cdot\xi}a\(\frac{x+y}{2},\xi\)\hat{J}\(
\frac{p_{\alpha}\(x\)}{t}\)J\(\frac{p_{\alpha}\( y\)}{t}\)\Phi\(y\)dyd\xi.
\end{equation*}
Introduce the operator ${}^tL_{\xi}=\frac{(
x-y)\cdot\partial_{\xi}}{i|x-y|^2}$. We have for any $k\in \N$: 
\begin{equation*}
R(t,x)= 
\iint e^{i(x- y)\cdot \xi}L^k_{\xi}a\(\frac{x+y}{2},\xi\)\hat{J}\(
\frac{p_{\alpha}\(x\)}{t}\)J\(\frac{p_{\alpha}\(y\)}{t}\)\Phi\(y\)dyd\xi. 
\end{equation*}
We treat the case $0<\alpha<2$, the other case is analogous. Notice that: 
\begin{equation*}
\left.
\begin{aligned} 
y\in \supp J\(\frac{p_{\alpha}\(
\cdot\)}{t}\)&\Rightarrow |y|\ge \varepsilon^{\frac{2}{2-\alpha}}
t^{\frac{2}{2-\alpha}}-1\\ 
x\in \supp \hat{J}\(\frac{p_{\alpha}\( \cdot\)}{t}\)&\Rightarrow |x|\le
\delta'^{\frac{2}{2-\alpha}} t^{\frac{2}{2-\alpha}}
\end{aligned}\right\}\Rightarrow|x-y|
\ge\(\varepsilon^{\frac{2}{2-\alpha}}-
\delta'^{\frac{2}{2-\alpha}}\)t^{\frac{2}{2-\alpha}}-1.  
\end{equation*}
We infer:
\begin{align*}
\left| R (t,x) \right|
&\lesssim\left|\hat{J}\(\frac{p_{\alpha}\(x\)}{t}\)\right | 
\iint|x-y|^{-k}\<\xi\>^{\eta m-k}|\Phi\(y\)|dyd\xi\\
&\lesssim\left|\hat{J}\(\frac{p_{\alpha}\(x\)}{t}\)\right|
t^{-\frac{k}{2-\alpha}}\|\Phi\|_{L^2}\, ,
\end{align*}
for any $k$, and $t$ sufficiently large. Thus,
\begin{equation*}
\left\|\hat{J}\(\frac{p_{\alpha}\(x\)}{t}\)\( {\rm Op}\(a\)J\(
  \frac{p_{\alpha}\(.\)}{t}\)\Phi\)\(x\)\right\|\lesssim
  t^{-\frac{k}{2\(2-\alpha\)}}\|\Phi\|_{L^2}\, ,
\end{equation*}
for any $k$, and $t$ sufficiently large.
\end{proof}

\begin{lem}
\label{L.P.5}
Let $\chi\in C_0^{\infty}(\R;\R_+)$, $J\in C_b^{\infty}(\R;\R_+)$,  with
$J=0$ in a  neighborhood of zero. Then we have for some
$\varepsilon>0$: 
\begin{itemize}
\item[(i)] Denote $\Theta (t) =\displaystyle
\sigma_{\alpha}\chi\(H_{\a,0}\)J^2\(
\frac{p_{\alpha}\(x\)}{t}\)\chi\(H_{\a,0}\)$.  Then:
\begin{equation*}
-\Theta (t) +{\mathcal
O}\(t^{-\varepsilon}\) \leq \chi\(H_{\a,0}\)J\(\frac{p_{\alpha}\(x\)}{t}\) 
\V_{\alpha}J\(\frac{p_{\alpha}\(x\)}{t}\)\chi\(H_{\a,0}\)\leq 
\Theta (t) +{\mathcal
O}\(t^{-\varepsilon}\) .
\end{equation*}
\item[(ii)]$\displaystyle \chi\(H_{\a,0}\)J\(\frac{p_{\alpha}\(x\)}{t}\)
\tA_{\alpha}J\(\frac{p_{\alpha}\(x\)}{t}\)\chi\(H_{\a,0}\) \ge {\mathcal
O}\(t^{-1-\varepsilon}\).$ 
\item[(iii)]$\displaystyle \chi\(H_{\a,0}\)
J\(\frac{p_{\alpha}\(x\)}{t}\)
\tA_{\alpha}\(\sigma_{\alpha}-\V_\a\)
J\(\frac{p_{\alpha}\(x\)}{t}\)\chi\(H_{\a,0}\)  \ge {\mathcal
O}\(t^{-1-\varepsilon}\).$ 
\end{itemize}
\end{lem}

\begin{proof}
We start with an inequality which we will use in the following. Let $\varepsilon>0$ such
that $J=0$ on $\left[-\varepsilon',\varepsilon'\right]$, for some
$\varepsilon'>\varepsilon$. For $1\leq j\leq d$, let $c_j\in
S\(\<\xi\>^{m_j}\<x\>^{-k_j},g_2\)$, 
with $0\le m_j\le 3$, and let
$l:=\min_j\{k_j-\alpha m_j/2\}\ge 0$. We suppose that for
$\psi\in C_0^{\infty}(\R)$, $ \psi=1$ in a small neighborhood of 0, we
have: 
\begin{equation}
\label{P.1.9} 
c\(x,\xi\)=\sum_{j=1}^d
c_j\(x,\xi\)\ge-C\<x\>^{-q\(\beta\)},\quad\,\text{on }\,
\supp\psi\(\frac{\xi^2-\<x\>^{\alpha}}{\(\xi^2+
\<x\>^{\alpha}\)^{\beta}}\),  
\end{equation}
with $q(\beta) >0$. For $1\ge\beta>3/4$, let
$\gamma:=\min\{1,q\(\beta\),\alpha\beta+l\}$. We prove that: 
\begin{equation}
\label{P.1.10}
\chi\(H_{\a,0}\)J\(\frac{p_{\alpha}\(x\)}{t}\){\rm Op}\(
c\)J\(\frac{p_{\alpha}\(x\)}{t}\)\chi\(H_{\a,0}\)\ge
\left\{\begin{aligned} 
&{\mathcal O}\(t^{-\gamma\frac{2}{2-\alpha}}\) \quad && \text{if }0<\alpha<2,\\
&{\mathcal O}\(e^{-\varepsilon \gamma t}\) && \text{if }\alpha=2.
\end{aligned} 
\right.
\end{equation}
Before proving \eqref{P.1.10}, we show that it implies the
lemma. First, on 
$\supp\psi\(\frac{\xi^2-\<x\>^{\alpha}}{\(\xi^2+
\<x\>^{\alpha}\)^{\beta}}\)$, we have: 
\begin{equation}
\label{P.1.15}
|\xi^2-\<x\>^{\alpha}|\lesssim\<x\>^{\alpha\beta}.
\end{equation}
We start with proving $\(i\)$. We have on
$\supp\psi\(\frac{\xi^2-\<\xi\>^{\alpha}}{\(\xi^2+
\<x\>^{\alpha}\)^{\beta}}\)$:  
\begin{equation}
\label{C.11.a}
|{\v}_{\alpha}\(x,\xi\)|\le \sigma_{\alpha}
\frac{|x|\(\<x\>^{\alpha/2}+C\<x\>^{\alpha\beta/2}\)}{\<x\>^{1+\alpha/2}} 
\le\sigma_{\alpha}+C\<x\>^{\(\beta-1\)\alpha/2}.
\end{equation}
We have $q\(\beta\)=\(1-\beta\)\alpha/2>0$ and
$\alpha\beta+l=\alpha\beta>0$,  for all $1>\beta>3/4$. This yields
$(i)$. In order to prove $(ii)$, we decompose
$\ta_{\alpha}=\ta^1_{\alpha}+\ta^2_{\alpha}$. We then use that on
$\supp\psi\(\frac{\xi^2-\<x\>^{\alpha}}{\(\xi^2+
\<x\>^{\alpha}\)^{\beta}}\)$:  
\begin{equation}
\label{C.11.b}
\ta_{\alpha}\(x,\xi\)\ge2\sigma_{\alpha}\xi^2
\frac{\<x\>^2-x^2}{\<x\>^{3+\alpha/2}}+\alpha\sigma_{\alpha}x^2
\frac{\<x\>^{\alpha}-\xi^2}{\<x\>^{3+\alpha/2}}\
\ge-C\<x\>^{\alpha\beta-1-\alpha/2}.
\end{equation}
Now observe that $\displaystyle \frac{2}{2-\alpha}>1$, and:
\begin{align*}
 \(1+\alpha\(\frac{1}{2}-\beta\)\)\frac{2}{2-\alpha}>1 \quad
 \text{and} \quad
 \(\alpha\beta+1-\frac{\a}{2}\)\frac{2}{2-\alpha}>1,&\quad \text{for }
\frac{3}{4}<\beta<1. 
\end{align*}
This yields $(ii)$.
Let us prove $(iii)$. We have:
\begin{equation*}
\tA_{\alpha}\(\sigma_{\alpha}-\V_\a\)= {\rm
Op}\(\ta_{\alpha}\(\sigma_{\alpha}-\v_{\alpha}\)\)+{\rm
Op}\({r}_1+{r}_2\),
\end{equation*}
with:
\begin{equation*}
{r}_1\in S\(\<x\>^{-2},g_2\),\,{r}_2 \in
S\(\<\xi\>^2\<x\>^{-2-\alpha},g_2\). 
\end{equation*}
In particular we have, on $\supp\psi\(\frac{\xi^2-
\<x\>^{\alpha}}{\(\xi^2+\<x\>^{\alpha}\)^{\beta}}\)$: 
\begin{equation*}
 {r}_j\ge -C\<x\>^{-2}\quad \(j=1,2\).
\end{equation*}
We apply \eqref{P.1.10} to ${r}_1+{r}_2$, and find
$l=2=q\(\beta\)$. Therefore $\gamma\frac{2}{2-\alpha}>1$.  

Let us consider ${\rm Op}\(\ta_{\alpha}
\(\sigma_{\alpha}-\v_{\alpha}\)\)$. We have 
\begin{align*}
\ta_{\alpha}&=\frac{1}{p_{\alpha}\(x\)}
\frac{2+\alpha}{\sigma_{\alpha}}\(\sigma_{\alpha}^2-
\v_{\alpha}^2\)+{r}_3\ ;\quad  
{r}_3=2\sigma_{\alpha}\frac{\xi^2-
\<x\>^{\alpha}}{\<x\>^{1+\alpha/2}}-\alpha\sigma_{\alpha}
\frac{1}{\<x\>^{3-\alpha/2}}={r}_3^1+ {r}^2_3\, ,\\
{r}^1_3&\in S\(\<\xi\>^2\<x\>^{-1- \alpha/2},g_2\),\quad
{r}^2_3\in S\(\<x\>^{-1+\alpha/2},g_2\). 
\end{align*}
We find on $\supp\psi\(\frac{\xi^2-\<x\>^{\alpha}}{\(\xi^2+
\<x\>^{\alpha}\)^{\beta}}\)$:  
\begin{equation*}
{r}_3\ge -C\<x\>^{-1-\alpha/2+\alpha\beta}.
\end{equation*}
Using that $|\v_{\alpha}|\lesssim 1$ on
$\supp\psi\(\frac{\xi^2-\<x\>^{\alpha}}{\(\xi^2+
\<x\>^{\alpha}\)^{\beta}}\)$, we find:
\begin{equation*}
{e}_3:={r}_3\(\sigma_{\alpha}- \v_{\alpha}\)\ge -
C\<x\>^{-1-\alpha/2+\alpha\beta}\quad\text{on}\,
\supp\psi\(\frac{\xi^2-\<x\>^{\alpha}}{\(\xi^2+\<x\>^{\alpha}\)^{\beta}}\). 
\end{equation*}
We decompose
${e}_3={e}^1_3+{e}^2_3+{e}^3_3+{e}^4_3$ with:  
\begin{align*}
{e}^1_3\in S\(\<\xi\>^2\<x\>^{-1-\alpha/2},g_2\)&\quad ;\quad {e}^2_3\in
S\(\<x\>^{-1+\alpha/2},g_2\),\\ 
{e}^3_3\in S\(\<\xi\>^{3}\<x\>^{-1-\alpha},g_2\)&\quad ;\quad {e}^4_3\in
S\(\<\xi\>\<x\>^{-1},g_2\). 
\end{align*}
We apply \eqref{P.1.10}, and find
$l=1-\alpha/2$, $q(\beta)=1+\alpha/2-\alpha\beta$. In particular,
$\gamma\frac{2}{2-\alpha}>1$ if $\beta<1$. It remains to consider: 
\begin{equation*}
b_{\alpha} :=\frac{2+\alpha}{\sigma_{\alpha}}
\frac{1}{p_{\alpha}}\(\sigma_{\alpha}-\v_{\alpha}\)^2\(\sigma_{\alpha}+
\v_{\alpha}\).
\end{equation*}
We use \eqref{C.11.a} and find 
\begin{equation*}
b_{\alpha}\ge -C\<x\>^{-1+\alpha\beta/2}\quad
\text{on }\,\supp\psi\(\frac{\xi^2-\<x\>^{\alpha}}{\(\xi^2+
\<x\>^{\alpha}\)^{\beta}}\). 
\end{equation*}
Notice that $b_\a=\sum_j b_{\alpha}^j$, with
$b_{\alpha}^j\in
S\(\<x\>^{-1+\alpha/2}\<\xi\>^{m_j}\<x\>^{-\frac{\alpha
m_j}{2}},g_2\)$, $0\le m_j\le 3$. We have $l=1-\alpha/2$ and
$q\(\beta\)=1-\alpha\beta/2$. In particular, $\gamma\frac{2}{2-\alpha}>1$
if $\beta<1$. This yields $(iii)$.

It remains to show \eqref{P.1.10}. Recall from Proposition~\ref{pse} that:
\begin{equation*}
\chi\(H_{\a,0}\)=\chi\(H_{\a,0}\){\rm Op}
\(\psi\(\frac{\xi^2-\<x\>^{\alpha}}{\(\xi^2+\<x\>^{\alpha}\)^{\beta}}\)\)+
R_{\beta}\, , 
\end{equation*}
with $N_{\alpha}^{\beta}R_{\beta}$ bounded. Let
\begin{equation*}
g_4=\frac{|dx|^2}{\<x\>^{2\(\alpha\beta+1-\alpha\)}}+
\frac{|d\xi|^2}{\<x\>^{2\alpha\(\beta-1/2\)}}\,. 
\end{equation*}
Notice that
$\psi\(\frac{\xi^2-\<x\>^{\alpha}}{\(\xi^2+ \<x\>^{\alpha}\)^{\beta}}\)\in  
S\(1,g_4\)\cap  S\(1,g_{\alpha\beta+1-\alpha,2\beta-1}\)$. 
We have:
\begin{align*}
\chi\(H_{\a,0}\) J \bigg( \frac{p_{\alpha}\(x\)}{t} & \bigg) {\rm Op}\(c\)
J\(\frac{p_{\alpha}\(x\)}{t}\)\chi\(H_{\a,0}\)\\ 
=&\chi\(H_{\a,0}\){\rm Op} \left( \psi\(
\frac{\xi^2-\<x\>^{\alpha}}{\(\xi^2+
\<x\>^{\alpha}\)^{\beta}} \) \right) J\(\frac{p_{\alpha}\(x\)}{t}\){\rm Op}\(c\) \times \\
&\times J\(\frac{p_{\alpha}\(x\)}{t}\){\rm Op} \left( \psi \( \frac{\xi^2-\<x\>^{\alpha}}{\(\xi^2+ \<x\>^{\alpha}\)^{\beta}} \) \right) \chi\(H_{\a,0}\)\\  
&+R_{\beta}J\(\frac{p_{\alpha}\(x\)}{t}\){\rm Op}\(c\)
J\(\frac{p_{\alpha}\(x\)}{t}\){\rm Op}\(\psi\(\frac{\xi^2-
\<x\>^{\alpha}}{\(\xi^2+ \<x\>^{\alpha}\)^{\beta}}\)\)
\chi\(H_{\a,0}\)\\
&+R_{\beta}J\(\frac{p_{\alpha}\(x\)}{t}\){\rm Op}\(c\)
J\(\frac{p_{\alpha}\(x\)}{t}\)R_{\beta}\\ 
=:& F_1 + F_2 + F_3.
\end{align*}
Let us start with estimating $F_2$. Using Lemma~\ref{L.P.4}, we find
$\widetilde{J}\in C_b^{\infty}(\R)$, $\widetilde{J}J=J$,
$\widetilde{J}=0$ on $\left[-\varepsilon,\varepsilon\right]$, such that: 
\begin{align*}
F_2 =&R_{\beta}N_{\alpha}^{\beta}\widetilde{J}\(\frac{p_{\alpha}\(x\)}{t}\)
N_{\alpha}^{-\beta} J\(\frac{p_{\alpha}\(x\)}{t}\){\rm Op}\(c\)
J\(\frac{p_{\alpha}\(x\)}{t}\) {\rm
Op} \left( \psi\(\frac{\xi^2-\<x\>^{\alpha}}{\(\xi^2
+\<x\>^{\alpha}\)^{\beta}} \right) \)\chi\(H_{\a,0}\)\\
&+ {\mathcal
O}\(t^{-\infty}\)\\ 
=&R_{\beta}N_{\alpha}^{\beta} \widetilde{J}\(\frac{p_{\alpha}\(x\)}{t}\) {\rm
Op}\(\frac{c\(x,\xi\)}{\(\xi^2+\<x\>^{\alpha}\)^{\beta}} 
\psi\(\frac{\xi^2-\<x\>^{\alpha}}{\(\xi^2+ \<x\>^{\alpha}\)^{\beta}}\)
J^2\(\frac{p_{\alpha}\(x\)}{t}\)\) \widetilde{J}\(\frac{p_{\alpha}\(x\)}{t}\)
\chi\(H_{\a,0}\)\\ 
&+R_{\beta}N_{\alpha}^{\beta}\widetilde{J} \(\frac{p_{\alpha}\(x\)}{t}\) {\rm
Op}\(e_2^t\) \widetilde{J}\(\frac{p_{\alpha}\(x\)}{t}\) \chi\(H_{\a,0}\)+
{\mathcal O}\(t^{-\infty}\). 
\end{align*}
We estimate
\begin{align*}
|f_2\(x,\xi\)|&=\left|\frac{c\(x,\xi\)}{\(\xi^2+
 \<x\>^{\alpha}\)^{\beta}} 
\psi\(\frac{\xi^2- \<x\>^{\alpha}}{\(\xi^2+
\<x\>^{\alpha}\)^{\beta}}\)J^2\(\frac{p_{\alpha}\(x\)}{t}\)\right|\\ 
&\lesssim \sum_{j=1}^d\<x\>^{-\alpha\beta}\<\xi\>^{m_j}\<x\>^{-k_j}
\psi\(\frac{\xi^2-\<x\>^{\alpha}}{\(\xi^2+\<x\>^{\alpha}\)^{\beta}}\)\\ 
&\lesssim\sum_{j=1}^d\<x\>^{-\alpha\beta} \<x\>^{\frac{\alpha
m_j}{2}-k_j} 
\lesssim\<x\>^{-\alpha\beta-l}\, ,\quad\text{ uniformly in $t$}.
\end{align*}
Thus,
\begin{equation*}
F_2=\left\{
\begin{aligned}
&{\mathcal O}\(t^{-\(\alpha\beta+l\)
\frac{2}{2-\alpha}}\) \quad  &&\text{if } 0<\alpha<2,\\ 
&{\mathcal O}\(e^{-\varepsilon\(\alpha\beta+l\)t}\) \quad  &&\text{if } \alpha=2.
\end{aligned}
\right.
\end{equation*}
The estimate for $F_3$ is analogous; we have to estimate
\begin{align*}
|f_3\(x,\xi\)|&\lesssim\sum_{j=1}^d\frac{1}{\(\xi^2+
\<x\>^{\alpha}\)^{2\beta}}\<\xi\>^{m_j}\<x\>^{-k_j}
\lesssim\<x\>^{-2\beta\alpha-l},\quad\text{if $1\ge\beta>3/4$.}
\end{align*}
Using the same arguments as in the estimates for $F_2$, $F_3$, we find:
\begin{align*}
F_1=&\chi\(H_{\a,0}\)\widetilde{J}\(\frac{p_{\alpha}\(x\)}{t}\){\rm Op}
\(c\(x,\xi\)\psi^2\(\frac{\xi^2-\<x\>^{\alpha}}{\(\xi^2+
\<x\>^{\alpha}\)^{\beta}}\) 
J^2\(\frac{p_{\alpha}\(x\)}{t}\)\)\widetilde{J}
\(\frac{p_{\alpha}\(x\)}{t}\)\chi\(H_{\a,0}\)\\
&+\chi\(H_{\a,0}\)\widetilde{J}
\(\frac{p_{\alpha}\(x\)}{t}\){\rm Op}\(e_1^t\)\widetilde{J}
\(\frac{p_{\alpha}\(x\)}{t}\)\chi\(H_{\a,0}\)+{\mathcal
O}\(t^{-\infty}\).
\end{align*}
Using \cite[Theorem 18.5.5]{Ho85} on the composition of
pseudo-differential operators in
$\Psi\(g_{\alpha\beta+1-\alpha,2\beta-1}\)$ and $\Psi\(g_2\)$, as well as
the fact that $\<\xi\>\lesssim\<x\>^{\alpha/2}$ on $\supp\psi$, we find:
$$e_1^t\in S\(\<x\>^{\alpha/2-\alpha\beta-1},1/2\(g_2+g_4\)\)\subset
S\(\<x\>^{-1},1/2\(g_2+g_4\)\),$$
because $\beta>1/2$. 
Using \eqref{P.1.9} and the G\aa rding inequality in $\Psi\(g_4\)$ we get:
\begin{equation*}
F_1\ge \chi\(H_{\a,0}\)\widetilde{J}\(\frac{p_{\alpha}\(x\)}{t}\) {\rm
Op}\(\hat{e}_1^t\)\widetilde{J}\(\frac{p_{\alpha}\(x\)}{t}\)
\chi\(H_{\a,0}\)+{\mathcal O}\(t^{-\frac{2}{2-\alpha}}\), 
\end{equation*}
with $\hat{e}^t_1\in
S\(\<x\>^{-\min\{q\(\beta\),2\alpha\beta+1-3\alpha/2\}},g_4\)\subset
S\(\<x\>^{-\min\{q\(\beta\),1\}},g_4\)$ ($\beta>3/4$), uniformly in
$t\ge 1$. Therefore:
\begin{equation*}
F_1\ge
\left\{
\begin{aligned} 
&{\mathcal O}\(t^{-\min\{q\(\beta\),1\}\frac{2}{2-\alpha}}\) \quad  &&\text{if } 0<\alpha<2,\\ 
&{\mathcal O}\(e^{-\min\{q\(\beta\),1\}\varepsilon t}\) \quad  &&\text{if } \alpha=2. 
\end{aligned}
\right. 
\end{equation*}
The estimates for $F_1$, $F_2$ and $F_3$ yield \eqref{P.1.10}.
\end{proof}

\subsection{Propagation estimates}
\label{subsP.2}

\begin{lem}   \label{L.P.Z}
Let $H$ be a self-adjoint operator on a Hilbert space ${\mathcal H}$,
and let ${\bf D}_H$ be the associated Heisenberg derivative. Let
$\Theta\(t\)$ be a uniformly bounded observable. We suppose that
\begin{equation}
\label{P.C.1}
{\bf D}_H\Theta\(t\)\ge f\(t\)\in L^1\(\R^+, dt\).
\end{equation}
Then the limit 
\begin{equation*}
\lim_{t\rightarrow\infty}\(\Theta\(t\)e^{-itH}\Phi|e^{-itH}\Phi\)
\end{equation*}
exists for all $\Phi\in {\mathcal H}$.
\end{lem}
\begin{proof}
We set $\Phi_t=e^{-itH}\Phi$. By \eqref{P.C.1}, we have:
\begin{equation*}
\frac{d}{dt}\(\(\Theta\(t\)\Phi_t|\Phi_t\)-F\(t\)\)\ge
0,\quad\mbox{with }\,
F(t)=\int_1^{t}\(f(s)\Phi|\Phi\)ds. 
\end{equation*}
Thus $\(\Theta (t)\Phi_t|\Phi_t\)-F(t)$ is increasing and
bounded. Therefore,  the limit 
\begin{equation*}
\lim_{t\rightarrow\infty}\(\(\Theta\(t\)\Phi_t|\Phi_t\)-F\(t\)\)
\end{equation*}
exists. Since
$\displaystyle \lim_{t\rightarrow\infty}F(t)= 
\int_1^{\infty}\(f(s)\Phi|\Phi\)ds$ 
exists, this gives the lemma. 
\end{proof}

Proposition~\ref{P.P.6} yields a minimal velocity
estimate. There is also a maximal velocity estimate: 
\begin{prop}[Maximal Velocity Estimate]
\label{P.P.7}
Let $\chi\in C_0^{\infty}(\R)$,
$\sigma_{\alpha}<\theta_2<\theta_3<\infty$. Then:  
\begin{itemize}
\item[(i)]$\displaystyle \int_1^{\infty}\left\|{\bf
1}_{\left[\theta_2,\theta_3\right]} 
\(\frac{p_{\alpha}(x)}{t}\)e^{-itH_{\a,0}}\chi\(H_{\a,0}\)
\Phi\right\|^2\frac{dt}{t}\lesssim \|\Phi\|^2$.
\item[(ii)]Let  $F\in C^{\infty}(\R)$,
with $ F'\in
C_0^{\infty}(\R)$ and $ {\rm supp} F\subset
]\theta_2,\infty [$. Then:
\begin{equation*}
  s-\lim_{t\rightarrow\infty}F\(
\frac{p_{\alpha}(x)}{t}\)e^{-itH_{\a,0}}=0. 
\end{equation*}
\end{itemize}
\end{prop}

\begin{proof}
$\(i\)$ Let $f\in C_0^{\infty}(\R)$, with $ f=1$ on
$\left[\theta_2,\theta_3\right]$, 
and ${\rm supp} f\subset ]\sigma_{\alpha},\infty[$. Let 
\begin{equation*}
F (s):=\int_{-\infty}^sf^2(s)ds\quad ;\quad \Theta(t):=
\chi\(H_{\a,0}\) F \(\frac{p_{\alpha}(x)}{t}\)\chi\(H_{\a,0}\). 
\end{equation*}
We compute:
\begin{align*}
-{\bf D}\Theta\(t\)&=\chi\(H_{\a,0}\)
f^2\(\frac{p_{\alpha}(x)}{t}\)\frac{p_{\alpha}(x)}{t^2}
\chi\(H_{\a,0}\) 
-\frac{1}{t}\chi\(H_{\a,0}\) f\(\frac{p_{\alpha}(x)}{t}\)
\V_\a f\(\frac{p_{\alpha}(x)}{t}\)\chi\(H_{\a,0}\)\\ 
&\ge\frac{\mu}{t}\chi\(H_{\a,0}\) f^2\(\frac{p_{\alpha}(x)}{t}\)
\chi\(H_{\a,0}\)+{\mathcal O}\(t^{-1-\varepsilon}\) ,
\end{align*}
for some $\mu,\varepsilon>0$. We have used Lemma~\ref{L.P.5}. This proves $(i)$, using \cite[Lemma B.4.1]{DG}. 
\esp

$(ii)$ Let the function $F$ satisfy the conditions in $(ii)$. Clearly
we can assume that $F\ge0$, and $F\(s\)=1$ for $s\ge
R_0$. Let $\widetilde{f}\in C_0^{\infty}(\R)$ be such that $\widetilde{f}=1$
on ${\rm supp} F'$, and ${\rm
supp}\widetilde{f}\subset ]\theta_2,\infty [$. Then:
\begin{equation}
\label{P.2.1}
{\bf D}\Theta\(t\)=\frac{1}{t}\chi\(H_{\a,0}\)
\widetilde{f}\(\frac{p_{\alpha}(x)}{t}\)B\(t\)\widetilde{f}
\(\frac{p_{\alpha}(x)}{t}\)\chi\(H_{\a,0}\)+
{\mathcal O}\(t^{-1-\varepsilon}\),
\end{equation}
with $B\(t\)$ uniformly bounded in $t$. To see that \eqref{P.2.1} is
true, we introduce $\widetilde{\chi}\in C_0^{\infty}(\R)$ with
$\widetilde{\chi}\chi=\chi$. Using Lemma~\ref{L.P.3} to estimate the
commutator
$\left[\widetilde{\chi}\(H_{\a,0}\),\widetilde{f}
\(\frac{p_{\alpha}(x)}{t}\)\right]$, 
and arguments similar to the arguments used in the proof of
Lemma~\ref{L.P.5}, we check that $B(t)$ is uniformly bounded in $t\ge
1$. From $(i)$, 
there exists  
\begin{equation} 
\label{P.2.2}
s-\lim_{t\rightarrow\infty}e^{itH_{\a,0}}
\Theta\(t\)e^{-itH_{\a,0}}.    
\end{equation}
If, in addition, $F$ is compactly supported, then by $(i)$ we have:
\begin{equation*}
\int_1^{\infty}\(\Theta\(t\)e^{-itH_{\a,0}}
\Phi|e^{-itH_{\a,0}}\Phi\)\frac{dt}{t}\lesssim \|\Phi\|^2. 
\end{equation*}
Thus if $F$ satisfies the conditions in $(ii)$, and is compactly
supported, then the limit \eqref{P.2.2} is zero. Now take 
 $F_1\in C^{\infty}(\R)$, $ f\in C_0^{\infty}(\R)$ such that
${\rm supp} F_1\subset ]\theta_2,\infty [$ with $F_1=1$ in a neighborhood
of $\infty$, and $F_1'=f^2$. Set 
\begin{equation*}
\Theta_R\(t\):=\chi\(H_{\a,0}\)F_1\(\frac{p_{\alpha}(x)}{Rt}\)
\chi\(H_{\a,0}\) .
\end{equation*}
From the previous discussion, we know that for $R> 0$, the limit 
$\displaystyle s-\lim_{t\rightarrow\infty}e^{itH_{\a,0}}
\Theta_R\(t\)e^{-itH_{\a,0}}$ exists. 
Repeating the computations of the proof of $(i)$, and keeping track of
$R$, we obtain: 
\begin{align*}
-{\bf D}\Theta_R\(t\)=&\, \frac{1}{t}\chi\(H_{\a,0}\)
f^2\(\frac{p_{\alpha}(x)}{Rt}\)\frac{p_{\alpha}(x)}{Rt}\chi\(H_{\a,0}\) \\
&-\frac{1}{t}\chi\(H_{\a,0}\)f\(
\frac{p_{\alpha}(x)}{Rt}\)\frac{\V_\a}{R}
f\(\frac{p_{\alpha}(x)}{Rt}\)\chi\(H_{\a,0}\)\\ 
\ge&\, \frac{1}{t}\(\theta_2-\sigma_{\alpha} /R\)
\chi\(H_{\a,0}\)f^2\(\frac{p_{\alpha}(x)}{Rt}\)\chi\(H_{\a,0}\)
+{\mathcal O}\(\(tR\)^{-1-\varepsilon}\) .
\end{align*}
Hence for $R\ge 1$, $ - {\bf D}\Theta_R\(t\)\ge{\mathcal
O}\(\(tR\)^{-1-\varepsilon}\)$. Therefore, for $t_0\ge 1$, we have:
\begin{equation*}
\begin{aligned}
s-\lim_{t\rightarrow\infty} e^{itH_{\a,0}}
\Theta_R\(t\)e^{-itH_{\a,0}}
&=e^{it_0H_{\a,0}}\Theta_R\(t_0\) e^{-it_0H_{\a,0}}+
\int_{t_0}^{\infty}e^{isH_{\a,0}}{\bf D}
\Theta_R\(s\)e^{-isH_{\a,0}}ds\\ 
&\le e^{it_0H_{\a,0}}\Theta_R\(t_0\) e^{-it_0H_{\a,0}}+
{\mathcal O}\(t_0^{-\varepsilon}R^{-1-\varepsilon}\).
\end{aligned}
\end{equation*}
For a fixed $t_0$, the terms on the right hand side  go
strongly to $0$ as $R\rightarrow\infty$, hence: 
\begin{equation}
\label{P.2.4}
s-\lim_{R\rightarrow\infty}\(s-\lim_{t\rightarrow\infty}
e^{itH_{\a,0}}\Theta_R\(t\)e^{-itH_{\a,0}}\)=0. 
\end{equation}
We remark now that, for $R\ge1$, the function
$F_1\(p_{\alpha}(x)\)-F_1\(\frac{p_{\alpha}(x)}{R}\)$ has a compact
support included in $[\theta_2,\infty [$. So, 
\begin{equation}
\label{P.2.5}
s-\lim_{t\rightarrow\infty}e^{itH_{\a,0}}
\(\Theta_1\(t\)-\Theta_R\(t\)\)e^{-itH_{\a,0}}=0. 
\end{equation}
Letting $R$ go to infinity in \eqref{P.2.5} and using \eqref{P.2.4},
we obtain: 
\begin{eqnarray*}
s-\lim_{t\rightarrow\infty}e^{itH_{\a,0}}\
\Theta_1\(t\)e^{-itH_{\a,0}}=0. 
\end{eqnarray*}
This completes the proof of $(ii)$.
\end{proof}

The next estimate is a weak propagation estimate:
\begin{prop}
\label{P.P.8}
Let $\chi\in C_0^{\infty}(\R),\, 0<\theta_1<\theta_2.$ Then:
\begin{itemize}
\item[(i)] $\displaystyle \int_1^{\infty}\left\|{\bf
1}_{\left[\theta_1,\theta_2\right]} 
\(\frac{p_{\alpha}(x)}{t}\) 
\(\frac{p_{\alpha}(x)}{t}-\V_\a\)e^{-itH_{\a,0}}
\chi\(H_{\a,0}\)\Phi\right\|^2\frac{dt}{t}\lesssim \|\Phi\|^2$.
\esp

\item[(ii)] $\displaystyle s-\lim_{t\rightarrow\infty}  {\bf
1}_{\left[\theta_1,\theta_2\right]}\(\frac{p_{\alpha}(x)}{t}\)
\(\frac{p_{\alpha}(x)}{t}-\V_\a\)e^{-itH_{\a,0}}=0$.   
\end{itemize}
\end{prop}
\begin{proof}
Let $R\in C^{\infty}(\R)$, with $R''\ge 0$, $R'=0$ on
$\left[-\varepsilon,\varepsilon\right]$ for some $\varepsilon>0$, and 
$R(x)=x^2/2$ for $|x|\ge\theta_1$. We take $ \theta_3>\max\{\theta_2,
\sigma_{\alpha}\}$. Let $J\in 
C_0^{\infty}(\R)$, with $J=1 $ on $\left[0,\theta_3\right]$. We set: 
\begin{align*}
M\(t\)&:=\frac{1}{2}\(\V_\a-\frac{p_{\alpha}(x)}{t}\)
R'\(\frac{p_{\alpha}(x)}{t}\)+ \frac{1}{2}R'\(\frac{p_{\alpha}(x)}{t}\)
\(\V_\a-\frac{p_{\alpha}(x)}{t}\)+R\(\frac{p_{\alpha}(x)}{t}\),\\ 
\Theta\(t\)&:=
\chi\(H_{\a,0}\)
J\(\frac{p_{\alpha}(x)}{t}\)M\(t\)J\(\frac{p_{\alpha}(x)}{t}\) 
\chi\(H_{\a,0}\).
\end{align*}
The observable $\Theta\(t\)$ is bounded uniformly in $t$ (see
Lemmas~\ref{L.P.4} and \ref{L.P.5}). We have: 
\begin{align*}
{\bf D}\Theta\(t\)=&\, -\chi\(H_{\a,0}\)
J'\(\frac{p_{\alpha}(x)}{t}\)\frac{p_{\alpha}(x)}{t^2}M\(t\)
J\(\frac{p_{\alpha}(x)}{t}\)\chi\(H_{\a,0}\)\\ 
&+\frac{1}{2}\chi\(H_{\a,0}\)\(J'\(\frac{p_{\alpha}(x)}{t}\)
\frac{\V_\a}{t}+hc\)M\(t\)J\(\frac{p_{\alpha}(x)}{t}\)\chi\(H_{\a,0}\)
+hc\\
&+\chi\(H_{\a,0}\)J\(\frac{p_{\alpha}(x)}{t}\){\bf D}M\(t\)
J\(\frac{p_{\alpha}(x)}{t}\)\chi\(H_{\a,0}\)\\ 
=:&\, R_1+R_2+R_3\,.
\end{align*}
The first two terms are of the form
\begin{equation*}
\frac{1}{t}\chi\(H_{\a,0}\) \widetilde{j}\(\frac{p_{\alpha}(x)}{t}\)
B\(t\)\widetilde{j}\(\frac{p_{\alpha}(x)}{t}\)\chi\(H_{\a,0}\)+ {\mathcal
O}\(t^{-1-\varepsilon}\), 
\end{equation*}
with $\widetilde{j}\in C_0^{\infty}(\R)$, $ {\rm
supp}\widetilde{j}\subset ]\sigma_{\alpha},\infty [$ (see
Lemmas~\ref{L.P.3} and \ref{L.P.4}). They are integrable, from
Proposition~\ref{P.P.7}. For the last term, we have:

\begin{align*}
R_3=&\, \frac{1}{2}\chi\(H_{\a,0}\)J\(\frac{p_{\alpha}(x)}{t}\)
\left\{\tA_{\alpha}R'\(\frac{p_{\alpha}(x)}{t}\)+
hc\right\}J\(\frac{p_{\alpha}(x)}{t}\)\chi\(H_{\a,0}\)\\ 
&-\frac{1}{2}\chi\(H_{\a,0}\)
J\(\frac{p_{\alpha}(x)}{t}\)\left\{\frac{\V_\a}{t}R'
\(\frac{p_{\alpha}(x)}{t}\)+
hc\right\}J\(\frac{p_{\alpha}(x)}{t}\)\chi\(H_{\a,0}\)\\
&+\frac{1}{2}\chi\(H_{\a,0}\) J\(\frac{p_{\alpha}(x)}{t}\)
\left\{\(\V_\a- \frac{p_{\alpha}(x)}{t}\)
R''\(\frac{p_{\alpha}(x)}{t}\)\frac{\V_\a}{t}+hc\right\}
J\(\frac{p_{\alpha}(x)}{t}\)\chi\(H_{\a,0}\)\\ 
&+\chi\(H_{\a,0}\)
J\(\frac{p_{\alpha}(x)}{t}\)\frac{p_{\alpha}(x)}{t^2}
R'\(\frac{p_{\alpha}(x)}{t}\)J\(\frac{p_{\alpha}(x)}{t}\) 
\chi\(H_{\a,0}\)\\
&-\frac{1}{2}\chi\(H_{\a,0}\)
J\(\frac{p_{\alpha}(x)}{t}\)\left\{\(\V_\a-\frac{p_{\alpha}(x)}{t}\)
R''\(\frac{p_{\alpha}(x)}{t}\)\frac{p_{\alpha}(x)}{t^2}+hc\right\}
J\(\frac{p_{\alpha}(x)}{t}\)\chi\(H_{\a,0}\)\\ 
&+\frac{1}{2}\chi\(H_{\a,0}\)
J\(\frac{p_{\alpha}(x)}{t}\)\left\{\frac{\V_\a}{t} 
R'\(\frac{p_{\alpha}(x)}{t}\)+hc\right\}
J\(\frac{p_{\alpha}(x)}{t}\)\chi\(H_{\a,0}\)\\
&-\chi\(H_{\a,0}\)J\(\frac{p_{\alpha}(x)}{t}\)
\frac{p_{\alpha}(x)}{t^2}R'\(\frac{p_{\alpha}(x)}{t}\)
J\(\frac{p_{\alpha}(x)}{t}\)\chi\(H_{\a,0}\)+ {\mathcal O}\(t^{-2}\)\\
=&\, \frac{1}{2t}\chi\(H_{\a,0}\) J\(\frac{p_{\alpha}(x)}{t}\)
\(\V_\a-\frac{p_{\alpha}(x)}{t}\) R''\(\frac{p_{\alpha}(x)}{t}\)
\(\V_\a-\frac{p_{\alpha}(x)}{t}\)
J\(\frac{p_{\alpha}(x)}{t}\)\chi\(H_{\a,0}\)\\ 
&+\frac{1}{2}\chi\(H_{\a,0}\)J\(\frac{p_{\alpha}(x)}{t}\)
\{\tA_{\alpha}R'\(\frac{p_{\alpha}(x)}{t}\)+hc\}J\(
\frac{p_{\alpha}(x)}{t}\)\chi\(H_{\a,0}\)+{\mathcal O}\(t^{-2}\).
\end{align*}
Consider the second term. Let $\widetilde{J}\in
C_0^{\infty}(\R; \R^{+})$, with $\widetilde{J}J=J$ and
$R'\widetilde{J}=\hat{J}^2$. We have:  
\begin{align*}
\frac{1}{2}\(\tA_{\alpha}\hat{J}^2\(\frac{p_{\alpha}(x)}{t}\)+
\hat{J}^2\(\frac{p_{\alpha}(x)}{t}\) 
\tA_{\alpha}\)=&\, \hat{J}\(\frac{p_{\alpha}(x)}{t}\)\tA_{\alpha}
\hat{J}\(\frac{p_{\alpha}(x)}{t}\)\\
&+ \left[\hat{J}
\(\frac{p_{\alpha}(x)}{t}\),
\left[\hat{J}\(\frac{p_{\alpha}(x)}{t}\),\tA_{\alpha}\right]\right]. 
\end{align*}
We estimate the double commutator. From Lemma~\ref{L.P.4}, there
exists $\hat{J}_1$, with $\hat{J}_1=0$ in a neighborhood of zero and
$\hat{J}_1\hat{J}=\hat{J}$, such that: 
\begin{align*}
\left[\hat{J}\(\frac{p_{\alpha}(x)}{t}\),\left[\hat{J}
\(\frac{p_{\alpha}(x)}{t}\),
\tA_{\alpha}\right]\right]=&\,  \hat{J}_1\(\frac{p_{\alpha}(x)}{t}\)\left[
\hat{J}\(\frac{p_{\alpha}(x)}{t}\),\left[\hat{J}\(\frac{p_{\alpha}
(x)}{t}\),
\tA_{\alpha}\right]\right]\hat{J}_1\(\frac{p_{\alpha}(x)}{t}\)\\ 
&+{\mathcal O}\(t^{-\infty}\).
\end{align*}
Using
$\left[\hat{J}\(\frac{p_{\alpha}(x)}{t}\),\left[\hat{J}
\(\frac{p_{\alpha}(x)}{t}\),\tA_{\alpha}\right]\right]\in \Psi \(\<x\>^{-3-\alpha/2},g_2\)$ uniformly in $t$ we get: 
\begin{equation*}
\left[\hat{J}\(\frac{p_{\alpha}(x)}{t}\),\left[\hat{J}
\(\frac{p_{\alpha}(x)}{t}\),\tA_{\alpha}\right]\right]=
\left\{
\begin{aligned}
&{\mathcal O} \(t^{\(-3-\alpha/2\)\frac{2}{2-\alpha}}\) \quad && \text{if } 0<\alpha<2,\\ 
&{\mathcal O}\(e^{-\varepsilon t}\) \quad && \text{if } \alpha=2,
\end{aligned} 
\right.
\end{equation*}
for some $\varepsilon>0$. Putting all together and using Lemmas~\ref{L.P.3}, \ref{L.P.5} we
obtain:
\begin{equation*}
R_3\ge\frac{1}{2t}\chi\(H_{\a,0}\) \(\V_\a-
\frac{p_{\alpha}(x)}{t}\){\bf 1}_{\left[\theta_1, \theta_2\right]}\(
\frac{p_{\alpha}(x)}{t}\)\(\V_\a- \frac{p_{\alpha}(x)}{t}\)
\chi\(H_{\a,0}\)+{\mathcal O}\(t^{-1-\varepsilon}\). 
\end{equation*}
This yields the desired estimate, thanks to \cite[Lemma B.4.1]{DG}.
\esp

$(ii)$ We can suppose $\sigma_{\alpha}\in \left[\theta_1,\theta_2\right]$. In
the other case, $(ii)$ follows from Propositions~\ref{P.P.6},
\ref{P.P.7}, and Lemma~\ref{L.P.2}. Let us first observe that:
\begin{align*}
s-\lim_{t\rightarrow\infty}{\bf 1}_{\left[\theta_1,
\theta_2\right]}\(\frac{p_{\alpha}(x)}{t}\) \(\V_\a-
\frac{p_{\alpha}(x)}{t}\)e^{-itH_{\a,0}}
&=s-\lim_{t\rightarrow\infty}{\bf 1}_{\left[\theta_1,
\theta_2\right]}\(\frac{p_{\alpha}(x)}{t}\)\(\V_\a-
\sigma_{\alpha}\)e^{-itH_{\a,0}},
\end{align*}
that is:
\begin{equation}
\label{C.A}
s-\lim_{t\rightarrow\infty}{\bf 1}_{\left[\theta_1,
\theta_2\right]}\(\frac{p_{\alpha}(x)}{t}\)\(\frac{p_{\alpha}(x)}{t}-
\sigma_{\alpha}\)e^{-itH_{\a,0}}=0. 
\end{equation}
Indeed, let $\varepsilon>0$. Then for $\Phi\in {\mathcal H}$, we have:
\begin{align*}
\left\|{\bf 1}_{\left[\theta_1,\theta_2\right]}\(
\frac{p_{\alpha}}{t}\)\(\frac{p_{\alpha}}{t}-
\sigma_{\alpha}\)e^{-itH_{\a,0}}\Phi\right\| 
&\le\varepsilon+\left\|{\bf 1}_{\left[\theta_1,
\theta_2\right]\setminus\left[\sigma_{\alpha}-\varepsilon,\sigma_{\alpha}+
\varepsilon\right]}\(\frac{p_{\alpha}}{t}\)\(\frac{p_{\alpha}}{t}-
\sigma_{\alpha}\)e^{-itH_{\a,0}}\Phi\right\|\\
&\le 2\varepsilon\, , 
\end{align*}
for $t$ sufficiently large, using Propositions~\ref{P.P.6} and
 \ref{P.P.7}. This yields \eqref{C.A}. Let $J\in
 C_0^{\infty}\(\R^*\)$, with $J\ge 0$ and $J(x)=1$ in a neighborhood of
 $\left[\theta_1,\theta_2\right]$.  
 Set 
\begin{equation*}
\Theta\(t\):=\chi\(H_{\a,0}\)\(\V_\a- \sigma_{\alpha}\)
 J^2\(\frac{p_{\alpha}(x)}{t}\)\(\V_\a-\sigma_{\alpha}\)
 \chi\(H_{\a,0}\). 
\end{equation*}
We compute:
\begin{align*}
-{\bf D}\Theta\(t\)=&\, \chi\(H_{\a,0}\)\tA_{\alpha}
J^2\(\frac{p_{\alpha}(x)}{t}\)\(\sigma_{\alpha}-\V_\a\)
\chi\(H_{\a,0}\)\\ 
&-\frac{1}{t}\chi\(H_{\a,0}\)\(\V_\a- \sigma_{\alpha}\)
\(J^2\)'\(\frac{p_{\alpha}(x)}{t}\)\(\V_\a-
\frac{p_{\alpha}(x)}{t}\)\(\V_\a-\sigma_{\alpha}\)\chi\(H_{\a,0}\)+hc. 
\end{align*}
The second term is integrable along the evolution, by
Propositions~\ref{P.P.6}, \ref{P.P.7} and Lemma~\ref{L.P.2} 
(the derivative of $J^2$ is zero in a neighborhood of
$\sigma_{\alpha}$). Using Lemmas~\ref{L.P.3} and \ref{L.P.5},  we find: 
\begin{align*}
\chi\(H_{\a,0}\)\tA_{\alpha}
J^2&\(\frac{p_{\alpha}(x)}{t}\)\(\V_\a-\sigma_{\alpha}\)
\chi\(H_{\a,0}\) \\  
&=\chi\(H_{\a,0}\)J\(\frac{p_{\alpha}(x)}{t}\) \tA_{\alpha}\(\V_\a-
\sigma_{\alpha}\)J\(\frac{p_{\alpha}(x)}{t}\) \chi\(H_{\a,0}\)+
{\mathcal O}\(t^{-1-\varepsilon}\)\\ 
&\ge{\mathcal O}\(t^{-1-\varepsilon}\),
\end{align*}
for some $\varepsilon>0$. Set $\Phi_t=e^{-itH_{\a,0}}\Phi$.
By Lemma~\ref{L.P.Z}, the limit $\displaystyle
\lim_{t\rightarrow\infty}\(\Theta\(t\)\Phi_t|\Phi_t\)$
exists. Let
\begin{equation*}
\widetilde{\Theta}\(t\)=\chi\(H_{\a,0}\)\(\V_\a-
\frac{p_{\alpha}(x)}{t}\)J^2\(\frac{p_{\alpha}(x)}{t}\)\(\V_\a-
\frac{p_{\alpha}(x)}{t}\)\chi\(H_{\a,0}\). 
\end{equation*}
We have
\begin{align*}
\widetilde{\Theta}\(t\)-\Theta\(t\)=&\, \chi\(H_{\a,0}\) \bigg( \(\sigma_{\alpha}-
\frac{p_{\alpha}(x)}{t}\)J^2\(\frac{p_{\alpha}(x)}{t}\)\(\V_\a-
\frac{p_{\alpha}(x)}{t}\) \\
&+\(\V_\a- \sigma_{\alpha}\) J^2\(\frac{p_{\alpha}(x)}{t}\)\(\sigma_{\alpha}- \frac{p_{\alpha}(x)}{t}\) \bigg) \chi\(H_{\a,0}\) .
\end{align*}
Using \eqref{C.A}, we obtain:
\begin{equation*}
\lim_{t\rightarrow\infty}\(\Theta\(t\)\Phi_t|\Phi_t\)=
\lim_{t\rightarrow\infty}\(\widetilde{\Theta}\(t\)\Phi_t|\Phi_t\). 
\end{equation*}
But by $(i)$, we have $\displaystyle
\int_1^{\infty}\(\Phi_t|\widetilde{\Theta}\(t\)\Phi_t\) \frac{dt}{t}\lesssim
\|\Phi\|^2$. 
Hence the limit is zero.
\end{proof}

\subsection{Asymptotic velocity}\label{subsP.3}
We now have all the technical tools to prove Theorem~\ref{Th.P.1}. 
\begin{prop}
\label{P.P.9}
Let $J\in \CI(\R)$. Then there exists
\begin{equation*}
s-\lim_{t\rightarrow\infty}e^{itH_{\alpha}}
J\(\frac{p_{\alpha}(x)}{t}\)e^{-itH_{\alpha}}. 
\end{equation*}
Moreover, if $J\(0\)=1$, then 
$\displaystyle s-\lim_{R\rightarrow\infty}\(s-\lim_{t\rightarrow\infty}
e^{itH_{\alpha}}J\(\frac{p_{\alpha}(x)}{Rt}\)e^{-itH_{\alpha}}\)={\bf
1}$.\\
If we define 
\begin{equation*}
P_{\alpha}^+=s-\CI-\lim_{t\rightarrow\infty}
e^{itH_{\alpha}}\frac{p_{\alpha}(x)}{t}e^{-itH_{\alpha}}, 
\end{equation*}
then $P_{\alpha}^+$ is a self-adjoint operator, which commutes with
$H_{\alpha}$. 
\end{prop}
\begin{proof}
We prove the proposition in two steps:
\esp

\noindent{\bf First step.} We assume $V_\a\equiv 0$: $H_{\alpha}=H_{\a,0}$.

By density, we may assume that $J\in C_0^{\infty}(\R)$, and that $J$ is
constant in a neighborhood of $0$ and in a neighborhood of
$\sigma_{\alpha}$. It also suffices to prove the existence of  
\begin{align*}
s-\lim_{t\rightarrow\infty}&e^{itH_{\a,0}}
J\(\frac{p_{\alpha}(x)}{t}\)e^{-itH_{\a,0}}\chi^2\(H_{\a,0}\)\\
=s-\lim_{t\rightarrow\infty}& \chi\(H_{\a,0}\)e^{itH_{\a,0}}
J\(\frac{p_{\alpha}(x)}{t}\)e^{-itH_{\a,0}} \chi\(H_{\a,0}\),
\end{align*}
for any $\chi\in C_0^{\infty}(\R)$, using Lemma~\ref{L.P.3}. Let 
$ \displaystyle \Theta\(t\):=\chi\(H_{\a,0}\)J\(\frac{p_{\alpha}(x)}{t}\)
\chi\(H_{\a,0}\)$.
We compute:
\begin{equation*}
{\bf D}\Theta\(t\)=\chi\(H_{\a,0}\)\frac{1}{2}
\(J'\(\frac{p_{\alpha}(x)}{t}\)\frac{\V_\a}{t}+hc\)\chi\(H_{\a,0}\)-
\chi\(H_{\a,0}\)J'\(\frac{p_{\alpha}(x)}{t}\)
\frac{p_{\alpha}(x)}{t^2}\chi\(H_{\a,0}\). 
\end{equation*}
This is integrable along the evolution by
Lemmas~\ref{L.P.2}, \ref{L.P.3} and Propositions~\ref{P.P.6},
\ref{P.P.7}.
\esp

\noindent{\bf Second step.} General
case. 

Let $P_{\a,0}^+$ be the asymptotic velocity associated with
$H_{\a,0}$. Using Theorem~\ref{comple}, we obtain the
existence of $P_{\alpha}^+$ by the formula: 
\begin{equation}
\label{P.3.2}
J\(P_{\alpha}^+\)=\Omega^+J\(P_{\a,0}^+\)\(\Omega^+\)^*+ J\(0\){\bf
1}^{pp}\(H_{\alpha}\). 
\end{equation}
The fact that $H_{\alpha}$ commutes with $P_{\alpha}^+$ follows from
Lemma~\ref{L.P.3}.
\end{proof}

\begin{prop}
\label{P.P.10}
We have
$\displaystyle \sigma\(P_{\alpha}^+\)=
\left\{
\begin{aligned} 
&\{0,\sigma_{\alpha}\} \quad &&\text{if} \quad \sigma_{\rm pp}\(H_{\alpha}\)\neq\emptyset, \\
&\{\sigma_{\alpha}\} \quad &&\text{if} \quad \sigma_{\rm pp}\(H_{\alpha}\)=\emptyset. 
\end{aligned}\right.$ 
\end{prop}

\begin{proof}
We first observe that
$\sigma\(P_{\a,0}^+\)\subset\{\sigma_{\alpha}\}$, by
Propositions~\ref{P.P.6} and \ref{P.P.7}. But the spectrum of
$P_{\a,0}^+$ cannot be empty, thus
$\sigma\(P_{\a,0}^+\)=\{\sigma_{\alpha}\}$. If
$\sigma_{\rm pp}\(H_{\alpha}\)=\emptyset$, then by \eqref{P.3.2},
$\sigma\(P_{\alpha}^+\)=\sigma\(P_{\a,0}^+\)=\{\sigma_{\alpha}\}$. We
suppose in the following that\footnote{ From Proposition~\ref{glm} and
Remark~\ref{rema:unique}, this does not seem to be the generic case.}
$\sigma_{\rm pp}\(H_{\alpha}\)\neq\emptyset$. By 
\eqref{P.3.2}, we have
$\sigma\(P_{\alpha}^+\)\subset\{0,\sigma_{\alpha}\}$. Let $J\(0\)\neq 0$,
and $\Phi\neq 0$ be an eigenfunction of $H_{\alpha}$. Then
$J\(P_{\alpha}^+\) \Phi =J\(0\)\Phi\neq0$, thus $0\in
\sigma\(P_{\alpha}^+\)$. Let 
now $J\(\sigma_{\alpha}\)\neq 0$, $J\(P_{\a,0}^+\)\Phi\neq 0$, and
$\psi=\Omega^+\Phi$. Since ${\rm Im}\Omega^+\subset{\bf
1}^c\(H_{\alpha}\)$,  we obtain by \eqref{P.3.2}:
$J\(P_{\alpha}^+\)\psi=\Omega^+J\(P^+_{\a,0}\)\Phi\neq 0$, in
particular $\sigma_{\alpha}\in \sigma\(P_{\alpha}^+\)$. 
\end{proof}

\begin{prop}
\label{P.P.11}
For $0<\a\leq 2$, we have:
$\displaystyle {\bf 1}_{\{0\}}\(P_{\alpha}^+\)={\bf
1}^{pp}\(H_{\alpha}\)$. 
\end{prop} 
\begin{proof}
Take (an approximation of) $J={\bf 1}_{\{0\}}$ in \eqref{P.3.2}, and
use ${\bf 1}_{\{0\}}\(P_{\a,0}^+\)=0$. 
\end{proof}
\begin{prop}
\label{P.P.12}
Let $g\in \CI(\R)$. Then
\begin{equation*}
s-\lim_{t\rightarrow\infty}e^{itH_{\alpha}}g\(\V_\a\)
e^{-itH_{\alpha}}{\bf 1}_{\R\setminus\{0\}}\(P_{\alpha}^+\)=
g\(P_{\alpha}^+\){\bf 1}_{\R\setminus\{0\}}\(P_{\alpha}^+\). 
\end{equation*}
\end{prop}
\begin{proof}
We first treat the case $H_{\alpha}=H_{\a,0}$. It is enough to
assume $g\in C_0^{\infty}(\R)$ and to prove that 
\begin{equation*}
s-\lim_{t\rightarrow\infty}e^{itH_{\a,0}}\( g\(\V_\a\)-
g\(\frac{p_{\alpha}(x)}{t}\)\)J\(\frac{p_{\alpha}(x)}{t}\)
\chi\(H_{\a,0}\)e^{-itH_{\a,0}}=0, 
\end{equation*}
for any $J\in C_0^{\infty}\(\R^*\)$ and $\chi\in C_0^{\infty}(\R)$. By
the Helffer--Sj\"ostrand formula, it is sufficient to show that for all
$z\in \C\setminus\(\sigma\(\V_\a\)\cup\R^+\)$: 
\begin{equation*}
s-\lim_{t\rightarrow\infty}\(z-\V_\a\)^{-1}\(\V_\a-
\frac{p_{\alpha}(x)}{t}\)\(z-\frac{p_{\alpha}(x)}{t}\)^{-1}
J\(\frac{p_{\alpha}(x)}{t}\)\chi\(H_{\a,0}\)e^{-itH_{\a,0}}=0. 
\end{equation*}
We have 
\begin{align*}
\(z-\V_\a\)^{-1} & \(\V_\a- \frac{p_{\alpha}(x)}{t}\)
\(z-\frac{p_{\alpha}(x)}{t}\)^{-1}J\(\frac{p_{\alpha}(x)}{t}\)
\chi\(H_{\a,0}\) \\  
=&\,\(z-\V_\a\)^{-1}\(z-\frac{p_{\alpha}(x)}{t}\)^{-1}
\(\V_\a-\frac{p_{\alpha}(x)}{t}\) J\(\frac{p_{\alpha}(x)}{t}\)
\chi\(H_{\a,0}\)\\ 
&+\(z-\V_\a\)^{-1}\(z- \frac{p_{\alpha}(x)}{t}\)^{-1}
\frac{\left[p_{\alpha},\V_\a\right]}{t}\(z- \frac{p_{\alpha}(x)}{t}\)^{-1}
J\(\frac{p_{\alpha}(x)}{t}\)\chi\(H_{\a,0}\) \\
=&:R_1+R_2\,.
\end{align*}
A direct computation shows that the commutator
$\left[\V_\a,p_{\alpha}\right]$ is bounded. Thus
$\displaystyle s-\lim_{t\rightarrow\infty} R_2 e^{-itH_{\a,0}}=0$. We have
$\displaystyle s-\lim_{t\rightarrow\infty} R_1 e^{-itH_{\a,0}}=0$ by
Proposition~\ref{P.P.8}.  

For the general case, we notice that
\begin{align*}
g\(P_{\alpha}^+\){\bf 1}_{\R\setminus\{0\}}\(P_{\alpha}^+\)& =
g\(P_{\alpha}^+\) {\bf 1}^c\(H_{\alpha}\)= \Omega^+g\(P_{\a,0}^+\)
\(\Omega^+\)^*\\
&=\Omega^+s-\lim_{t\rightarrow\infty} e^{itH_{\a,0}}
g\(\V_\a\)e^{-itH_{\a,0}}\(\Omega^+\)^*\\ 
&=s-\lim_{t\rightarrow\infty} e^{itH_{\alpha}} g\(\V_\a\)
e^{-itH_{\alpha}}{\bf 1}^c\(H_{\alpha}\)\\ 
&=s-\lim_{t\rightarrow\infty} e^{itH_{\alpha}} g\(\V_\a\)
e^{-itH_{\alpha}}{\bf 1}_{\R\setminus\{0\}}\(P_{\alpha}^+\), 
\end{align*}
which implies the proposition.
\end{proof}

Propositions~\ref{P.P.9}, \ref{P.P.10}, \ref{P.P.11} and \ref{P.P.12}
correspond to the four points in Theorem~\ref{Th.P.1}.

\section{Generalizations in the case $\alpha =2$}
\label{www}
\def\H{H}

In this section, we generalize our results, in the special case where
the reference Hamiltonian is exactly the Laplace operator plus a
second order polynomial. Resuming the notations of
Section~\ref{sec:mehler}, we may assume, as in \eqref{tx21} that 
\begin{equation}  \label{tx22}
H_0 = -\Delta - \sum_{k=1}^{n_{-}} \omega_{k}^{2} x_{k}^{2} +
\sum_{k=n_{-} + 1}^{n_{-} + n_{+}} \omega_{k}^{2} x_{k}^{2} +
\sum_{k=n_{-} + n_{+} +1}^{n_{-} + n_{+} + n_{E}} E_{k} x_{k}=:
-\Delta +U(x), 
\end{equation}
on $L^{2} (\R^{n})$, with $n_{-} + n_{+} + n_{E} \leq n$, $\omega_{k}
>0$ and $E_{k} \neq 0$. By convention $\sum_{j=a}^{b} =0$ if $b<a$. We
study the scattering theory for the Hamiltonians $(H_{0} , H = H_{0} +
V (x))$. This setting includes the presence of a Stark potential
($n_{E} \neq 0$). Notice that, if $U(x)$ is a general second order
polynomial with real coefficients, the operator 
\begin{align*}
-\Delta + U(x),
\end{align*}
can always be written as \eqref{tx22}, modulo a constant term, after a
change of orthonormal basis and origin (which leaves the Laplace
operator invariant). With that approach, one could even demand
$n_E\leq 1$. The reason why we do not reduce to that case is that the
pointwise decay estimates required for the potential (see
\eqref{hyp:W} below) are not invariant with respect to such
reductions.

In Section~\ref{sec:mehlercook}, we prove the existence of wave
operators under rather weak assumptions on the perturbative potential
$V$. In the case $n_-=n_+=1$ for instance, the repulsive effects due
to $-x_1^2$ overwhelm the other effects: confinement due to $+x_2^2$
and drift due to the Stark potential.

In Section~\ref{gena2}, we give the asymptotic completeness if
$H_{0}$ has no Stark effect and no Schr\"{o}dinger part (this means
that $n_{+} + n_{-} =n$). The hypothesis on $V$ are similar to
\eqref{va1}--\eqref{va2}.

In Section~\ref{subsP.4}, we construct the asymptotic velocity under
the previous hypothesis. As the free Hamiltonian $H_{0}$ is a sum of commuting
self-adjoint ``one-dimensional'' operators, the existence of
asymptotic velocities in each space direction is a corollary of
Theorem~\ref{Th.P.1} applied to the one-dimensional case. For the
Hamiltonian $H$, the asymptotic velocity of
Theorem~\ref{Th.P.1} exists also and is equal to $P_{\ell}^+$, where 
$\omega_\ell=\max_{1\leq j\leq n_-}\omega_j$, and $P_{\ell}^+$ is
the asymptotic velocity in the direction $x_\ell$.

\subsection{Existence of wave operators}
\label{sec:mehlercook}

In this section, we prove the existence of wave operators for
perturbations of $H_0$ by Cook's method. We consider the perturbation
$H=H_0+V$, where $V(x)$ is a real-valued potential which can be
decomposed as 
\begin{equation}  \label{tx25}
V (x) = V_{1} (x) + V_{2} (x) + W (x) ,
\end{equation}
where
\begin{equation}  \label{tx24}
V_{j} \in L^{p_{j}} (\R^{n} ; \R) \text{, for } j=1,2, \text{ with }
\left\{ 
\begin{aligned}
&2 \leq p_{j} < \infty  &\text{ if } n \leq 3 ,  \\
&2 < p_{j} < \infty  &\text{ if } n = 4 ,  \\ 
&n/2 \leq p_{j} < \infty  &\text{ if } n \geq 5 ,
\end{aligned}
\right.
\end{equation}
and $W$ is a sum of terms in $L^{\infty} (\R^{n})$ satisfying a.e.
\begin{equation}\label{hyp:W}
|W(x)| \lesssim \bigg( \prod_{j=1}^{n_-} \< \ln \< x \> \>^{-\beta_j}
 \bigg) \bigg( \prod_{j=n_-+n_++1}^{n_{-} + n_{+} + n_{E}} \<
 x_j\>^{-\beta_j /2} \bigg) \bigg( \prod_{j=n_-+n_+ n_{E} +1}^{n} \<
 x_j\>^{-\beta_j} \bigg) , 
\end{equation}
with $\beta_j\geq 0$ and $\sum \beta_j >1$. Notice that the $V_j$'s do
not contain pointwise information.

\begin{theo}
\label{theo:cook}
$(i)$ Suppose that the quadratic part of $U$ has \emph{at most} one
(simple) positive eigenvalue ($n_+\leq 1$), and \emph{at least} one
negative eigenvalue ($n_-\geq 1$). Let $V$ satisfying the previous
assumptions. Then $H=H_0+V$ admits a unique self-adjoint extension,
and the following strong limits exist in $L^2(\R^n)$, 
\begin{equation*}
s-\lim_{t\to \pm \infty}e^{itH} e^{-it H_0}\, .
\end{equation*}
$(ii)$ If the quadratic part of $U$ has \emph{at least} one negative
eigenvalue ($n_-\geq 1$) and $V$ satisfies the previous assumptions
with $V_1$ and $V_2$ compactly supported, then the same conclusions
hold. \\ 
$(ii)$ If $W$ satisfies \eqref{hyp:W} and $V_1=V_2\equiv 0$, then the
same conclusions hold.  
\end{theo}

The self-adjointness property follows from Faris--Lavine
Theorem (\cite{ReedSimon2}, see also Section~\ref{ele}). Before the proof
of this result, a few remarks are in order.

\begin{rema}
This result shows that wave operators exist even for very slowly
decaying potentials. Potentials decaying even more slowly could be
included (involving $\ln (\ln |x|)$ for instance); like for
Theorem~\ref{comple}, we do not seek so
general results, and rather focus on the method. Notice that in the
first case, singular potentials, like 
\begin{equation*}
V(x)=\frac{1}{|x|^\delta} = \frac{1}{|x|^\delta}{\bf 1}_{|x|\leq 1} + 
\frac{1}{|x|^\delta}{\bf 1}_{|x|> 1}+ 0\, ,
\end{equation*}
are allowed, provided that $\delta <\min(2,n/2)$. This includes the case of
Coulomb potentials in space dimension $n\geq 3$. 
\end{rema}

\begin{rema}
The dynamics associated to $H_0$ is known
explicitly (see \eqref{eq:mehlergen}), and cannot be compared to
that of $-\Delta$.   
\end{rema}

\begin{rema}
If the quadratic part of $U$ has more than one positive
eigenvalue, results similar to the first point of the theorem 
can be proved, provided that $U$ has at
least one negative eigenvalue. This will be clear from the proof
below, as well as the reasons why we did not wish to state too general
a result. 
\end{rema}

Resuming the notations of Section~\ref{sec:mehler}, the
dilation operator ${\mathcal D}_t$ is crucial. As mentioned in
Section~\ref{sec:mehler}, a 
formula similar to \eqref{eq:factor} is available for
$e^{it\Delta}$. The factor ${\mathcal M}_t$ in that case is different,
but still of modulus one, while ${\mathcal D}_t$ corresponds to
dilations of size $2t$ instead of $g( 2t)$. This is closely related to
the properties of the classical trajectories. The operator ${\mathcal
D}_t$ enables us to prove the existence of wave operators with Cook's 
method.

It is of course sufficient to study the case $t\to +\infty$, the case
$t\to -\infty$ being similar. A density argument shows that
Theorem~\ref{theo:cook} follows from:  
\begin{lem}\label{lem:cook}
Under the assumptions of Theorem~\ref{theo:cook}, for any $\varphi =
\varphi_1\otimes\ldots\otimes\varphi_n$, with $\varphi_j \in \S(\R)$,
there exists a unique $\phi\in L^2(\R^n)$ such that  
\begin{equation*}
\left\| e^{it\H}e^{-it\H_0}\varphi -\phi\right\|_{L^2} \Tend t
{+\infty} 0 \, .
\end{equation*} 
\end{lem}

\begin{proof}[Proof of Lemma~\ref{lem:cook}]
Following Cook's method, we compute
\begin{equation*}
\frac{d}{dt}e^{it\H}e^{-it\H_0}\varphi= i e^{it\H}V e^{-it\H_0}\varphi\, .
\end{equation*}
Taking the $L^2$--norm, 
\begin{equation*}
\left\|\frac{d}{dt}e^{it\H}e^{-it\H_0}\varphi\right\|_{L^2}= 
\left\|Ve^{-it\H_0}\varphi\right\|_{L^2}\, .
\end{equation*}
Under the assumptions of Theorem~\ref{theo:cook}, it is sufficient to
prove that the maps 
\begin{equation*}
t \mapsto \left\|W e^{-it\H_0}\varphi\right\|_{L^2}\quad \text{and}\quad  
t \mapsto \left\|V_j e^{-it\H_0}\varphi\right\|_{L^2},\ j=1,2, 
\end{equation*}
are integrable on $[1,+\infty[$. Let $t\geq 1$. 
Since the operator ${\mathcal D}_t$ is unitary on $L^2$, we have, from
\eqref{eq:factor} and H\"older's inequality,  
\begin{equation}\label{eq:estfond}
\begin{aligned}
\left\|V_je^{-it\H_0}\varphi\right\|_{L^2}
&=\left\|V_j {\mathcal D}_t {\mathcal F}{\mathcal
M}_t\varphi\right\|_{L^2}  \\ 
&=\left\|{\mathcal D}_t V_j(\cdot g_1(2t),\ldots,\cdot g_n(2t) ) 
{\mathcal F}{\mathcal M}_t\varphi\right\|_{L^2}\\
&=\left\|V_j(\cdot g_1(2t),\ldots,\cdot g_n(2t) ) {\mathcal
F}{\mathcal M}_t\varphi\right\|_{L^2}\\ 
&\leq \left\|V_j(\cdot g_1(2t),\ldots,\cdot g_n(2t)
)\right\|_{L^{p_j}}\left\| {\mathcal 
F}{\mathcal M}_t\varphi\right\|_{L^{q_j}}\\  
&\lesssim \(\frac{1}{g_{n_-+1}(2 t)}
\prod_{k=1}^{n_-} \frac{\omega_k}{\sinh (2\omega_k t)}\)^{1/p_j}
\left\|V_j\right\|_{L^{p_j}}\left\|
{\mathcal M}_t\varphi\right\|_{L^{q_j'}}\, ,\\   
\end{aligned}
\end{equation}
where $\frac{1}{p_j}+\frac{1}{q_j}=\frac{1}{2}$ and the last estimate
stems from Hausdorff--Young inequality ($q_j'$ denotes the H\"older
conjugate exponent of $q_j$).

In the first case of Theorem~\ref{theo:cook},
the function $g_{n_-+1}$ is a sinus
if the quadratic part of $U$ has one positive eigenvalue ($n_+=1$),
and is linear otherwise.

The exponential decay of
$1/g_1(2t)$ is enough to ensure the integrability of the right hand
side of \eqref{eq:estfond}.
The worst possible situation with our assumptions is $n_+=1$, where
\eqref{eq:estfond} yields, since $n_-\geq 1$, 
\begin{equation*}
\left\|V_je^{-itH_0}\varphi\right\|_{L^2}
\lesssim \(\frac{1}{\sin(2\omega_{n_-+1} t)\sinh (2\omega_1 t)}
\)^{1/p_j}
\left\|V_j\right\|_{L^{p_j}}\left\|
\varphi\right\|_{L^{q_j'}}\, . 
\end{equation*}
From assumption (\ref{tx24}), $p_j\geq 2$, and the map 
$t\mapsto 1/ (\sin(2\omega_{n_-+1} t)\sinh(2\omega_1 t))^{1/p_j}$ is
integrable on $[1,+\infty[$. Since $L^{q_j'} \subset L^1\cap L^2$, the
lemma is proved for the $V_j$'s parts, in the first case of
Theorem~\ref{theo:cook}.
\esp

In the second case of Theorem~\ref{theo:cook}, we assume in addition
that the $V_j$'s are 
compactly supported: $\operatorname{supp} V_j \subset \{|x|\leq R\}=:B$
. In that case, they are $\H_0$--bounded (with 
bound zero); the assumptions on $p_j$ are such that $V_j$ is
$\Delta$--bounded, and the polynomial $U$ is bounded on the support of
$V_j$. To take advantage of this, write
$$\H_0 =  -\d_1^2 -\omega_1^2 x_1^2 + \widetilde\H_0\, ,$$
where $ \widetilde\H_0$ takes the last $(n-1)$ variables into
account, and use a factorization like \eqref{eq:factor} in the first
variable only. Mimicking the above computations with $x$ replaced by
$x_1$ yields
\begin{equation*}
\left\|V_j e^{-it\H_0}\varphi\right\|_{L^2} =\left\|V_j
e^{-it\H_0}\varphi\right\|_{L^2(B)} = \left\|V_j\(\cdot 
g_1(2t),\cdot,\ldots, \cdot\)
\F_1 {\mathcal M}^1_t e^{-it\widetilde\H_0}\varphi\right\|_{L^2(B_t)} , 
\end{equation*}
where $\F_1$ stands for the Fourier transform with respect to the
first variable, and 
$${\mathcal M}^1_t = \exp\( i 
x_1^2\frac{h_1(2t)}{2 g_1(2t)}\).$$
Notice that these two operators
commute with $\widetilde\H_0$. We also denoted $B_t= \{
x_1^2g_1(2t)^2+x_2^2+\ldots +x_n^2\leq R^2\}$. 
From H\"older's inequality, this last
term is estimated by
\begin{equation*}
 \left\|V_j\(\cdot
g_1(2t),\cdot,\ldots, \cdot\) \right\|_{L^{p_j}(B_t)} \left\|\F_1 {\mathcal
M}^1_t e^{-it\widetilde\H_0}\varphi\right\|_{L^{q_j}(B_t)},
\end{equation*}
where $\frac{1}{2}= \frac{1}{p_j}+\frac{1}{q_j}$. The first term is
equal to $(g_1(2t))^{-1/p_j} \|V_j\|_{L^{p_j}}$; since we assume
$n_-\geq 1$, $g_1$ has an exponential growth, and this term  is
integrable. It suffices to show that the second term is bounded. From
our assumption on $p_j$, $H^2(\R^n)\subset L^{q_j}(\R^n)$, and
\begin{equation*}
\begin{aligned}
\left\|\F_1 {\mathcal
M}^1_t e^{-it\widetilde\H_0}\varphi\right\|_{L^{q_j}(B_t)}& \lesssim
\left\|\F_1 {\mathcal 
M}^1_t e^{-it\widetilde\H_0}\varphi\right\|_{L^2} +\left\|\Delta\(\F_1
{\mathcal  
M}^1_t e^{-it\widetilde\H_0}\varphi\)\right\|_{L^2(B_t)}\\
&\lesssim \left\|\varphi\right\|_{L^2} + \left\|\( \d_1^2 -\widetilde
\H_0\)\(\F_1 {\mathcal  
M}^1_t e^{-it\widetilde\H_0}\varphi\)\right\|_{L^2(B_t)}\, , 
\end{aligned}
\end{equation*}
since $U$ is bounded on $B_t$, uniformly with respect to $t\geq 1$
($B_t\subset B$ for $t\geq 1$). Finally, we have
\begin{equation*}
\left\|\( \d_1^2 -\widetilde
\H_0\)\(\F_1 {\mathcal  
M}^1_t e^{-it\widetilde\H_0}\varphi\)\right\|_{L^2(\R^n)}=  \left\|
\(x_1^2 + \widetilde \H_0\)\varphi\right\|_{L^2(\R^n)}.
\end{equation*}
Therefore, the lemma is proved for the $V_j$'s parts, in the second case of
Theorem~\ref{theo:cook}.
\esp

For the last component of $V$ (the function $W$), we make no
assumption on $n_-$ or $n_+$. This is where we use the tensor product
structure for $\varphi$. The idea is to proceed as above, except for
the components $(x_{n_{-} + n_{+} +1}, \ldots , x_{n})$, for which we
proceed ``as usual''. We denote 
\begin{equation}   \label{tx33}
H_{0}^{j} = \left\{
\begin{aligned}
-&\Delta_{x_{j}} - \omega_{j}^{2} x_{j}^{2} \quad  &&\text{ for } 1
\leq j \leq n_{-}, \\ 
-&\Delta_{x_{j}} + \omega_{j}^{2} x_{j}^{2} \quad  &&\text{ for }
n_{-} < j \leq n_{-} + n_{+}, \\ 
-&\Delta_{x_{j}} + E_{j} x_{j} \quad  &&\text{ for } n_{-} + n_{+} < j
\leq n_{-} + n_{+} +n_{E}, \\ 
-&\Delta_{x_{j}} \quad  &&\text{ for } n_{-} + n_{+} +n_{E} < j \leq n,
\end{aligned}
\right. 
\end{equation}
and we have $\displaystyle e^{-it\H_0} \varphi = e^{-it\H_0^{1}}
\varphi_1 \otimes \ldots \otimes e^{-it\H_0^{n}} \varphi_n$.

Since $e^{-it H_{0}^{j}}$ is unitary, we have, for $j=n_{-}+1, \ldots , n_{-} +
n_{+}$:
\begin{equation}  \label{tx29}
\begin{split}
\left\| W e^{-it\H_0} \varphi \right\|_{L^{2}} \lesssim &
\prod_{j=1}^{n_{-}} \left\| \< \ln \< x_{j} \> \>^{-\beta_{j}}
e^{-it\H_{0}^{j}} \varphi_{j} \right\|_{L^{2}} \prod_{j=n_{-} + n_{+}
+1}^{n_{-} + n_{+} + n_{E}} \left\| \< x_{j} \>^{-\beta_{j} /2}
e^{-it\H_{0}^{j}} \varphi_{j} \right\|_{L^{2}} \times    \\ 
&\times \prod_{j=n_{-} + n_{+} + n_{E} + 1}^{n} \left\| \< x_{j}
\>^{-\beta_{j}} e^{-it\H_{0}^{j}} \varphi_{j} \right\|_{L^{2}}, 
\end{split}
\end{equation}
and we estimate each term in the product.

For $1\leq j\leq n_-$, we have, from \eqref{eq:factor} in the one
dimensional case, 
\begin{align*}
\left\| \< \ln \< x_{j} \> \>^{-\beta_{j}} e^{-it\H_{0}^{j}}
\varphi_{j} \right\|_{L^{2}}  =& \left\| \< \ln \< x_{j} \>
\>^{-\beta_{j}} {\mathcal M}^{j}_t {\mathcal D}^{j}_t {\mathcal F}_{j}
{\mathcal M}^{j}_t \varphi_{j} \right\|_{L^{2}}   \\ 
=& \left\| \< \ln \< x_{j} g_{j} (2t) \> \>^{-\beta_{j}}  {\mathcal
F}_{j} {\mathcal M}^{j}_t \varphi_{j} \right\|_{L^{2}}, 
\end{align*}
where ${\mathcal M}^{j}_t$, ${\mathcal D}^{j}_t$ and ${\mathcal
F}^{j}$ are given by 
\eqref{eq:factor} in dimension $1$. Denoting $\widetilde{\varphi}_{j}
= {\mathcal F}_{j} {\mathcal M}^{j}_t \varphi_{j}$ and replacing
$g_{j} (2t)$ with $e^{2\omega_{j} t}$, we get 
\begin{align*}
\left\| \< \ln \< x_{j} g_{j} (2t) \> \>^{-\beta_{j}}
\widetilde{\varphi}_{j} \right\|_{L^{2}}^{2} =\, &\int_{|x_j| >e^{-\omega_j
t}}\left| \< \ln \< x_{j} g_{j} (2t) \> \>^{-\beta_{j}} \widetilde
\varphi^j_t \right|^2dx_j\\
&  + \int_{|x_j| \leq e^{-\omega_j t}} \left|
\< \ln \< x_{j} g_{j} (2t) \> \>^{-\beta_{j}} \widetilde \varphi^j_t
\right|^2dx_j \\ 
\lesssim \, & \frac{1}{t^{2\beta_j} }\|\widetilde
\varphi^j_t\|_{L^2}^2 + \int {\bf 1}_{|x_j| \leq e^{-\omega_j
t}}\left|\widetilde \varphi^j_t \right|^2dx_j \\ 
\lesssim\,  & \frac{1}{t^{2\beta_j} }\|\varphi_j\|_{L^2}^2 +
e^{-\omega_j t/2}\|\widetilde \varphi^j_t\|_{L^4}^2 , 
\end{align*}
where we used Cauchy--Schwarz inequality. From Hausdorff--Young
inequality, we have:
\begin{equation}\label{eq:10:25}
\left\| \< \ln \< x_{j} \> \>^{-\beta_{j}} e^{-it\H_{0}^{j}}
\varphi_{j} \right\|_{L^2}^2 \lesssim
\frac{1}{t^{2\beta_j}}\|\varphi_j\|_{L^2}^2 + e^{-\omega_jt/2} \|
\varphi_j \|_{L^{4/3}}^2. 
\end{equation}

For the case $n_-+n_+ < j \leq n_-+n_+ + n_{E}$, we simply recall the
approach of \cite{AvronHerbst}. From Avron--Herbst formula (see
e.g. \cite{AvronHerbst,Cycon}),  
\begin{equation*}
\( e^{-it H_{0}^{j} } \varphi_{j} \) ( x_{j} ) = e^{-i \big( t E_{j}
x_{j} + \frac{t^3}{3} E_{j}^2 \big)} \( e^{it \Delta_{x_{j}}}
\varphi_{j} \) \( x_{j} + t^2 E_{j} \) . 
\end{equation*}
Using Avron--Herbst formula, the term we have to estimate reads
\begin{equation*}
\left\| \< x_{j}  \>^{-\beta_j/2} e^{-it H_{0}^{j} } \varphi_{j}
\right\|_{L^2} = \left\| \< x_{j} -t^2 E_{j} \>^{-\beta_j/2} e^{it
\Delta_{x_{j}}} \varphi_{j} \right\|_{L^2} . 
\end{equation*}
By a density argument, we may assume $\varphi_{j} \in {\mathcal F_{j}}
(C^\infty_0 (\R))$ ($\operatorname{supp} {\mathcal F}_{j}
(\varphi_{j}) \subset \{|\xi | \leq R \}$ for some positive $R$). For
$| x_{j} |<3Rt$, the drift caused by Stark effect accelerates the
particle:  
\begin{equation}
\left\| \< x_{j} -t^2 E_{j} \>^{-\beta_j/2} e^{it \Delta_{x_{j}}}
\varphi_{j} \right\|_{L^{2}(| x_{j} |<3Rt)} \lesssim | t^2
|^{-\beta_{j} /2} \| \varphi_{j} \| \lesssim t^{-\beta_{j}} . 
\end{equation}
For $| x_{j} |>3Rt$, a nonstationary phase argument and
$\operatorname{supp} {\mathcal F}_{j} (\varphi_{j}) \subset
\{|\xi|\leq R\}$ show that 
\begin{align}
\big\| e^{it \Delta_{x_{j}}} \varphi_{j} \big\|_{L^{2} (| x_{j}
|>3Rt)} =& \left\| \frac{1}{\sqrt{2 \pi}}  \int e^{-it \xi^2 +i x_{j}
\xi} {\mathcal F}_{j} (\varphi_{j}) (\xi) d \xi \right\|_{L^{2} (|
x_{j} |>3Rt)}  \nonumber \\ 
=& {\mathcal O} (t^{-\infty}).
\end{align}
These two estimates yield:
\begin{equation} \label{eq:10:53}
\left\| \< x_{j}  \>^{-\beta_j/2} e^{-it H_{0}^{j} } \varphi_{j}
\right\|_{L^2} \lesssim t^{-\beta_{j}} . 
\end{equation}

We now study the term with $n_-+n_+ + n_{E} < j \leq n$. By a density
argument, we can assume that $\varphi_{j}$ satisfies 
\begin{align*}
\varphi_{j} (x_{j}) = \iint e^{i(x_{j} -y) \xi} \chi(\xi) \psi(y) dy
d\xi ,  
\end{align*}
where $\psi\in C_0^\infty(\R)$ ($\operatorname{supp}\psi \subset
\{|y|<R\}$), $\chi\in C^\infty_0(\R, [0,1])$ and $\chi = 0$ for
$|\xi|\leq c$. We get 
\begin{align*}
\< x_{j} \>^{-\beta_{j}} e^{-it\H_{0}^{j}} \varphi_{j} = \< x_{j}
\>^{-\beta_{j}} \iint e^{-it \xi^2 +i(x_{j} -y) \xi} \chi(\xi) \psi(y)
dy d\xi . 
\end{align*}
Obviously, we have 
\begin{equation}  \label{tx30}
\left\| \< x_{j} \>^{-\beta_{j}} e^{-it\H_{0}^{j}} \varphi_{j}
\right\|_{L^2(|x| > ct )} \lesssim t^{-\beta_j}\| \varphi_{j} \|_{L^2}
\, . 
\end{equation}
For $|x|<ct$, recall that on the support of $\psi(y)\chi(\xi)$,
$|\xi|>c$. Differentiating the phase, and noting that 
\begin{align*}
| 2t\xi +x-y|\geq |2t\xi +x|-R \geq 2t|\xi|-|x| -R \geq ct-R\, ,
\end{align*}
a nonstationary phase argument shows that, for large $t$,  
\begin{equation}   \label{tx31}
\left\| \< x_{j} \>^{-\beta_{j}} \iint e^{-it \xi^2 +i(x_{j} -y) \xi}
\chi(\xi) \psi(y) dy d\xi \right\|_{L^2(|x|<ct)} = \O \( t^{-\infty}
\) \, . 
\end{equation}
Combining \eqref{tx30} and \eqref{tx31}, we get, for $n_-+n_+ + n_{E}
< j \leq n$, 
\begin{equation}   \label{eq:10:37}
\left\| \< x_{j} \>^{-\beta_{j}} e^{-it\H_{0}^{j}} \varphi_{j}
\right\|_{L^{2}} \lesssim t^{-\beta_{j}} . 
\end{equation}

Gathering \eqref{tx29}, \eqref{eq:10:25} and
\eqref{eq:10:37} together yields 
\begin{equation*}
\left\| W e^{-it\H_0} \varphi_{j} \right\|_{L^2} \lesssim t^{-\sum
\beta_j}.  
\end{equation*}
Since we assumed $\sum \beta_j>1$, the lemma follows.
\end{proof}

If the quadratic part of $U$ has more than one positive
eigenvalue, the problem may become intricate for the estimates related
to the $V_j$'s. 
For instance, if the quadratic part of $U$  has
one positive eigenvalue, whose order is $n_+\geq 2$, then
\eqref{eq:estfond} becomes
\begin{equation*}
\left\|V_je^{-itH_0}\varphi_+\right\|_{L^2}
\lesssim \(\frac{1}{|\sin(2\omega_{n_-+1} t)|^{n_+}\sinh (2\omega_1 t)}
\)^{1/p_j}
\left\|V_j\right\|_{L^{p_j}}\left\|
\varphi_+\right\|_{L^{q_j'}}\, . 
\end{equation*}
and if $n_+=2$ and $n= 3$ (see the Assumption~\ref{tx24}), one has to
adapt the  
assumption on the power $p_j$ for this map to be integrable near
$+\infty$: the value $p_j=2$ is not allowed, because for
that value, the above map is not even locally integrable.  

Reasoning the same way, we notice that if the quadratic part of $U$
has several distinct positive eigenvalues, then arithmetic properties
of these eigenvalues will have to be taken into account. We leave out
the discussion at this stage.

\subsection{Asymptotic completeness}
\label{gena2}

In this part, we assume that $n_{-} + n_{+} =n$. Then
\begin{equation*}
H_0 = -\Delta - \sum_{k=1}^{n_{-}} \omega_{k}^{2} x_{k}^{2} +
\sum_{k=n_{-} + 1}^{n} \omega_{k}^{2} x_{k}^{2}. 
\end{equation*}
Like for Theorem~\ref{comple}, we assume that $V$ is a real-valued
function with 
\begin{equation}  \label{va4}
V (x) = V_{1} (x) + V_{2} (x),
\end{equation}
where
\begin{equation}  \label{va3}
V_{1}  \text{ is a compactly supported measurable function, and }
\Delta \text{--compact}, 
\end{equation}
and $V_{2} \in L^{\infty} (\R^{n};\R)$ satisfies the short range
condition 
\begin{equation}    \label{SRo}
|V_{2} (x) | \lesssim \< \ln \< x_{-} \> \>^{-1-\varepsilon}
 ,\quad\text{a.e. }x\in\R^n\, , 
\end{equation}
for some $\varepsilon >0$, and $x=(x_{-},x_{+})\in \R^{n_-} \times
\R^{n_+}$. Notice that there is  propagation only in the $x_{-}$
direction. It is therefore reasonable to impose decay only in that
direction. Denote 
\begin{equation*} 
H^-_{0} =- \Delta_{x_{-}} - (\omega_{-} x_{-})^{2}\quad ;
 \quad H^+_{0}=-\Delta_{x_{+}} + (\omega_{+} x_{+})^{2},
\end{equation*}
with $\omega_{-} = {\rm diag} (\omega_{1}, \ldots , \omega_{n_{-}})$,
$\omega_{+} = {\rm diag} (\omega_{n_{-} +1 }, \ldots ,
\omega_{n})$. Let $N_{\pm}=-\Delta_{x_{\pm}} + x_{\pm}^2$. As in
Section~\ref{doma}, we can show that the operators $(H_{0}^{\pm} ,
D(N_{\pm}))$ and $(H,D(N))$ are essentially self-adjoint ($N=N_2$ is
the harmonic oscillator on $\R^n$). The result
of this section is:

\begin{theo}   \label{Th.gena2}
Assume that $V$ satisfies \eqref{va4}--\eqref{SRo}. Then there exist
\begin{gather}   \label{Wo1}
s-\lim_{t\rightarrow\infty} e^{itH} e^{-itH_{0}},   \\
s-\lim_{t\rightarrow\infty} e^{itH_{0}} e^{-itH} {\bf 1}^c (H).  \label{Wo2}
\end{gather}
If we denote \eqref{Wo1} by $\Omega^+$, then \eqref{Wo2} equals
$(\Omega^+)^*$ and we have 
\begin{align*}
(\Omega^+)^*\Omega^+={\bf 1},\quad \Omega^+(\Omega^+)^*= {\bf 1}^c(H).
\end{align*}
\end{theo}

\begin{rema}
From the following discussion, it is clear that the conditions
$\omega_j>0$ for all $j$ (resp. $n_-+n_+=n$) are crucial for the
proof. If $\omega_j=0$ for one $j$, then the equations \eqref{sss},
\eqref{es2} below fail to be true. For the same reason, we cannot include
linear terms like $E\cdot x$. 
\end{rema}

\begin{proof}
Since the proof is very similar to the case $n_-=n$ and $\omega_j=1$ for
all $j$, we will be very concise. The following points have to be
addressed: 
\begin{enumerate}
\item Definition of the conjugate operator $A$,

\item The regularity $H_{0}\in C^2(A)$,

\item The Mourre estimate for $H_{0}$,

\item The regularity $H \in C^{1+\delta} (A)$ and the Mourre estimate
for $H$, 

\item Replacement of the conjugate operator $A$ by $\< \ln\<x_{n_-}\>
\>$, 

\item Proof of the asymptotic completeness.
\end{enumerate}

(1) As in \eqref{nnn}, we choose for the conjugate operator
$A = {\rm Op} (a (x,\xi))$, with
\begin{equation*}
a(x,\xi) = \ln \< \xi_{-} + \omega_{-} x_{-} \> - \ln \< \xi_{-}
- \omega_{-} x_{-} \> .
\end{equation*}
One can show as in Lemma~\ref{nel01} that $(A,D(N))$ is essentially self-adjoint. As in Proposition~\ref{pse}, we can prove that, for $\psi \in C_0^{\infty} (\R)$ with $\psi =1$ near $0$, 
\begin{equation}   \label{pseudo}
(H +i)^{-1} = (H +i)^{-1} {\rm Op} \bigg( \psi \bigg( \frac{\xi^{2} -
( \omega_{-} x_{-} )^{2} + ( \omega_{+} x_{+} )^{2}}{(\xi^{2} + x^{2}
+1)^{\beta}} \bigg) \bigg) +{\mathcal O}(1){\rm Op}(r), 
\end{equation}
with $r \in S \left( \< x , \xi \>^{- \beta} , g_{0} \right)$ for
$1/2< \beta \leq 1$. If the support of $\psi$ is small enough, we have
\begin{equation}  \label{sss}
\xi_{-}^{2} + \xi_{+}^{2} + x_{+}^{2} \lesssim x_{-}^{2} +1\, ,
\end{equation}
on the support of $\psi \Big( \frac{\xi^{2} - ( \omega_{-} x_{-} )^{2} + ( \omega_{+} x_{+} )^{2}}{\xi^{2} + x^{2} +1} \Big)$.

(2) First, note that $[H_{0},A]$ and $[ [H_{0},A] ,A]$ are bounded on $L^{2} (\R^{n})$, by the same arguments as in Lemma~\ref{comutH0}  and \ref{doublecom}. Notice that $N=N_-+N_+$ and $D(N)=D(N_-)\cap D(N_+)$. Clearly, $D(N_+)=D(H_{0}^+)$. Since $[ H_{0}^+ , H_{0} ]=0$, we have:
\begin{equation*}
(z-H_0)^{-1}: D(N_+)\rightarrow D(N_+) .
\end{equation*}
To prove  $\displaystyle (z-H_0)^{-1}: D(N_-)\rightarrow D(N_-)$,
we can proceed as in the proof of Lemma~\ref{rere}; we use
\begin{equation*}
(z-H_0)^{-1} = \int_0^{\infty} e^{-itH_{0}^+} e^{-it H_{0}^-} e^{itz}dt,
\end{equation*}
and $[ N_-,H_{0}^+ ] =0$.

(3) As in Lemma~\ref{m20}, one can show that, for $\chi \in
    C_{0}^{\infty} (\R)$, 
\begin{equation*}
\chi(H_{0}) [iH_{0} , A] \chi(H_{0}) = \chi(H_{0}) \sum_{j=1}^{n_{-}}
    2 \omega_{j} \( \frac{(D_{j} + \omega_{j}  x_{j})^{2}}{\<
    D_{j} + \omega_{j}  x_{j} \>^{2}} + \frac{(D_{j} - \omega_{j}
    x_{j})^{2}}{\< D_{j} - \omega_{j}  x_{j} \>^{2}} \)
    \chi(H_{0}) .
\end{equation*}
Using G\aa rding inequality, we get, for any $\mu >0$, 
\begin{align*}
\chi(H_{0}) [iH_{0} , A] \chi(H_{0}) \geq (2 \hat{\omega} -\mu) \chi^{2} (H_{0}) + \chi(H_{0}) R \chi(H_{0}) , 
\end{align*}
where $\hat{\omega} = \min_{j \in \{1, \ldots ,n_- \} } \omega_j$, and
$R$ is a pseudo-differential operator whose symbol is decreasing in
$(x_{-}, \xi_{-})$. Using \eqref{pseudo} and \eqref{sss}, this decay
becomes a decay in $(x, \xi)$ on the energy level, and then 
\begin{align*}
\chi(H_{0}) [iH_{0} , A] \chi(H_{0}) \geq (2 \hat{\omega} -\mu)
\chi^{2} (H_{0}) + \chi(H_{0}) K \chi(H_{0}), 
\end{align*}
with $K$ compact. If the support of $\chi$ is sufficiently small, we
therefore obtain: 
\begin{equation*}
\chi(H_{0}) [iH_{0} , A] \chi(H_{0}) \geq (2 \hat{\omega} -\mu)
\chi^{2} (H_{0}). 
\end{equation*}

(4) Using \eqref{sss}, one can show that $V A$ is compact from
    $D(H_{0})$ to $L^{2} (\R^{n})$, since the decay of $V_{2}$
    in $x_{-}$ \eqref{SRo} yields decay in all the variables. Thus,
    the Mourre estimate and the regularity $C^{1+\delta} (A)$ for $H$
    can be obtained as in Section~\ref{sec:mourre}.

(5) We apply the same arguments as in Section~\ref{mimi}. We have to use that
\begin{equation}  \label{es2}
\< \xi_{-} - x_{-} \> \lesssim \< x_{-} \>\quad ;\quad
\< \xi_{-} + x_{-} \> \lesssim \< x_{-} \>,
\end{equation}
on the energy levels, and \eqref{sss}. Then we show that the assumptions
of Lemma~\ref{GeNi3} are fulfilled.

(6) The proof of the asymptotic completeness is exactly as in
    Section~\ref{subsec:demo}, using the minimal velocity estimate and
    the fact that $\chi(H)-\chi(H_{0})$ is compact for
    $\chi\in C_0^{\infty} (\R)$. 
\end{proof}

\subsection{Asymptotic velocity}
\label{subsP.4}

We assume that the hypothesis of Theorem~\ref{Th.gena2} are satisfied. Let
${\mathcal H}^k=L^2(\R)$ for  $k\in \{ 1, \ldots ,n \}$. Clearly
${\mathcal H} =  \otimes_{k=1}^n{\mathcal H}^k$. We
write, as in \eqref{tx33}, 
\begin{equation}  \label{P.4.1}
H_{0} = \sum_{j=1}^n H_{0, \omega_{j}}^j \, ,\quad \text{ with } \quad
H_{0, \omega_{j}}^j = -\Delta_{x_{j}} \pm \omega_j^2 x_j^2. 
\end{equation}
Obviously
\begin{equation}    \label{P.4.2}
[ H_{0, \omega_{j}}^j,H_{0, \omega_{k}}^k ]=0.
\end{equation}
If we use this separation of variables and apply Theorem~\ref{Th.P.1}
to the one-dimensional case, we obtain asymptotic velocities in each
space direction. Let $\V_j = [ H , \ln \< x_j \> ]$. Like for $\V$ we can
show that $(\V_j,D(N))$ is well defined as an operator, and essentially
self-adjoint. We denote again $\V_j$ its self-adjoint extension.

\begin{theo}[Asymptotic Velocities]
\label{Th.P.13}
There exists a vector $\vec{P}^+=(P_{1}^+, \ldots ,P_{n}^+)$ of
bounded self-adjoint commuting operators $P_{j}^+$ which commute with
$H$, such that:
\begin{itemize}
\item[(i)] $\displaystyle
\vec{P}^+=s-C_{\infty}-\lim_{t\rightarrow\infty}e^{itH} 
\( \frac{\ln\<x_1\>}{t}, \ldots ,\frac{\ln\<x_{n}\>}{t} \)
e^{-itH}$.
\item[(ii)] For  $ j \in \{1, \ldots ,n_- \}$,  the spectrum of
$P_{j}^+$ is given by $\displaystyle\sigma (
P_{j}^+)= \left\{  
\begin{aligned}
\{ & 0,2\omega_j \} && \text{ if } \sigma_{pp} (H) \neq \emptyset ,  \\
\{ & 2 \omega_j \} && \text{ if } \sigma_{pp} (H) = \emptyset .
\end{aligned}
\right. $\\
For $ j\in\{n_-+1, \ldots ,n\}$ ,  the spectrum of
$P_{j}^+$ is reduced to $0$: $\sigma(P_{j}^+) = \{ 0 \}$.
\esp

\item[(iii)] For  $  j \in \{1, \ldots ,n_- \}$, $  {\bf 1}_{\{ 0
\}} (P_{j}^+) = {\bf 1}^{pp} (H)$.
\esp

\item[(iv)] For any $J \in C_{\infty} (\R)$, $\displaystyle J(P_{j}^+)
{\bf 1}_{\R\setminus\{0\}} (P_{j}^+) = s-\lim_{t \rightarrow \infty}
e^{itH} J(v_j) e^{-itH} {\bf 1}_{\R \setminus \{ 0 \}} (P_{j}^+)$.
\end{itemize}
\end{theo}

\begin{proof}
We denote $H_{\omega} = H$, and $H_{0,\omega} = H_{0}$. We prove
Theorem~\ref{Th.P.13} in two steps.
\esp

\noindent {\bf First step.} We assume $V\equiv 0$, that is $H_{\omega} =
H_{0,\omega}$.
\esp

$(i)$ We first treat the case $n_-=n=1$. If $\omega=(-1)$, then the
claim follows from Theorem~\ref{Th.P.1}. For
$\omega=(-\omega_1^2)$, set: 
\begin{equation*}
(\D v)(x)=\frac{1}{\omega_1^{1/4}} v \(\frac{x}{\sqrt{\omega_1}} \).
\end{equation*}
Then $\D:L^2(\R)\rightarrow L^2(\R)$ is unitary, and we have
\begin{equation}  \label{P.4.3}
H_{0,\omega}=\omega_1 \D^*H_{0,(-1)}\D.
\end{equation}
Therefore:
\begin{align}
s-\lim_{t\rightarrow\infty} e^{itH_{0,\omega}} J \(
\frac{\ln\<x\>}{t} \) e^{-itH_{0,\omega}} =&
\D^*s-\lim_{t\rightarrow\infty} e^{itH_{0,(-1)}} J \( \frac{\ln
\< \frac{x}{\sqrt{\omega_1}} \> \omega_1}{t} \)
e^{-itH_{0,(-1)}}\D  \nonumber   \\ 
=& \D^*s-\lim_{t\rightarrow\infty}e^{itH_{0,(-1)}}J \(
\frac{\ln\<x\>\omega_1}{t} \) e^{-itH_{0,(-1)}} \D , 
\label{P.4.4}
\end{align}
because $J \( \frac{ \ln \< \frac{x}{\sqrt{\omega_1}} \>
\omega_1}{t} \) - J \( \frac{\ln\<x\>\omega_1}{t} \) =
{\mathcal O}(t^{-1})$. Thus the result for general $\omega$ follows
from the result for $\omega=(-1)$. 
Let now $n_-=0$ and $n_+=1$. Then:
\begin{equation}  \label{P.4.5}
s-\lim_{t\rightarrow\infty} e^{itH_{0,\omega}} J \bigg(
\frac{\ln\<x\>}{t} \bigg) e^{-itH_{0,\omega}}=J(0) , 
\end{equation}
because $L^2(\R)$ possesses a basis of eigenfunctions of
$H_{0,\omega}$. The general case follows from the one-dimensional
cases using \eqref{P.4.1}, \eqref{P.4.2}. We denote $\vec{P}_{0}^+$ the
vector associated to $H_{0,\omega}$. 
\esp

$(ii)$ First note that
\begin{align*}
J(P_{0,j}^+)=s-\lim_{t\rightarrow\infty}e^{itH_{0,\omega_j}^j} J
\bigg( \frac{\ln\<x_j\>}{t} \bigg) e^{-itH_{0,\omega_j}^j}. 
\end{align*}
The result on the spectrum follows from \eqref{P.4.4}, \eqref{P.4.5},
and from Theorem~\ref{Th.P.1}, $(ii)$. 
\esp

$(iii)$ Both operators are zero.
\esp

$(iv)$ By $(ii)$, $P^+_{0,j}$ depends only on $\omega_j$, and we 
note it $P^+_{0,\omega_j}$. We have
$P^+_{0,\omega_j}=\omega_jP^+_{0,1}$. If $j\in\{n_-+1, \ldots ,n \}$,
both operators are zero. For $j\in\{1, \ldots ,n_-\}$, we have 
\begin{align*}
J(P_{0,\omega_j}^+)=&J(\omega_jP_{0,1}^+)=\D^*J(\omega_jP^+_{0,1})\D   
=\D^*s-\lim_{t\rightarrow\infty}e^{itH_{0,1}^j} 
J(\omega_j\V_j)e^{-itH_{0,1}^j}\D
\\ 
=&s-\lim_{t\rightarrow\infty} e^{itH_{0,\omega_j}^j}
\D^*J(\omega_j \V_j)\D e^{-itH_{0,\omega_j}^j} 
\end{align*}
Here we have used that $P_{0,1}^+=2$, and thus
$\D P_{0,1}^+\D^*=P_{0,1}^+$, as well as \eqref{P.4.3}. To prove $(iv)$, it
is thus sufficient to show that: 
\begin{equation*}
\begin{aligned}
& (\D^*J(\omega_j\V_j)\D-J(\V_j)) \chi(H_{0,\omega_j}^j) \text{ is compact
on $L^2(\R)$},  \\ 
&\text{ resp. } (J(\omega_j\V_j)-\D J(\V_j)\D^*) \chi(H_{0,1}^j) \text{ is
compact on $L^2(\R)$}, 
\end{aligned}
\end{equation*}
for any $\chi\in C_0^{\infty}(\R)$. All operators have to be
understood as operators acting on $L^2(\R)$. Let
$\widetilde{\V}_j=\D \V_j\D^*$. The operator $\widetilde{\V}_j$ is a
pseudo-differential operator, with symbol
${\widetilde{\v}}_j(x_j,\xi_j)=\frac{x_j\xi_j}{\big\< \frac{x_j}{\sqrt{\omega_j}} \big\>^2}$.
Recall from Proposition~\ref{pse} that 
\begin{align*}
\chi(H_{0,1}^j) = {\rm Op} \bigg( \psi \bigg( \frac{ \xi_j^2 - x_j^2}{
\xi_j^2 + \< x_j \>^2} \bigg) \bigg) \chi(H_{0,1}^j)+R_j , 
\end{align*}
with $N_jR_j$ bounded ($N_j= -\Delta_{x_{j}} + \<x_j\>^2$), and $\psi
\in C_0^{\infty}(\R)$ with $\psi=1$ in a neighborhood of zero. Clearly,
$(J(\omega_j \V_j)-J(\widetilde{\V}_j)))R_j$ is compact, and it remains to
show that: 
\begin{equation*}
\(J\(\omega_j\V_j\)-J\(\widetilde{\V}_j\)\){\rm Op} \( \psi \(
\frac{\xi_j^2 - x_j^2}{\xi_j^2 +\< x_j \>^2} \) \) \quad \text{is
compact.} 
\end{equation*}
By the Helffer--Sj\"ostrand formula, it is sufficient to
show that for any $z\in
\C\setminus(\sigma(\V_j)\cup\sigma(\widetilde{\V}_j))$ ,
\begin{align}
\(z-\widetilde{\V}_j\)^{-1}& \(\widetilde{\V}_j-\omega_j\V_j\)
\(z-\omega_j\V_j\)^{-1} 
{\rm Op} \( \psi \( \frac{\xi_j^2- x_j^2}{\xi_j^2 + \< x_j
\>^2} \) \) =   \nonumber  \\ 
& \begin{aligned}      \label{P.4.6}
&=\, \(z-\widetilde{\V}_j\)^{-1}\(\widetilde{\V}_j- \omega_j\V_j\) {\rm
Op} \( \psi 
\( \frac{\xi_j^2 - x_j^2}{\xi_j^2 + \<x_j\>^2} \) \)
(z-\omega_j\V_j)^{-1}   \\ 
\ &+\(z-\widetilde{\V}_j\)^{-1}\(\widetilde{\V}_j-\omega_j\V_j\)
\(z-\omega_j\V_j\)^{-1} 
\omega_j \left[ {\rm Op} \( \psi \( \frac{\xi_j^2 -
x_j^2}{\xi_j^2 + \< x_j \>^2} \) \) , \V_j \right]
\(z-\omega_j\V_j\)^{-1} , 
\end{aligned}
\end{align}
is compact. We have: 
\begin{align*}
\left| \( {\widetilde{\v}}_j(x_j,\xi_j) - \omega_j
{\v}_j(x_j,\xi_j) \) \psi \( \frac{\xi_j^2 - x_j^2}{\xi_j^2
+ \<x_j\>^2} \) \right| =& \left| x_j \xi_j \(
\frac{1-\omega_j}{\< x_j \>^2 \< \frac{x_j}{\sqrt{\omega_j}}
\>^2} \) \psi \( \frac{\xi_j^2 - x_j^2}{\xi_j^2 + x_j^2}
\) \right|    \\ 
\lesssim& \left| \frac{1-\omega_j}{\< \frac{x_j}{\sqrt{\omega_j}}
\>^2} \psi \( \frac{\xi_j^2 - x_j^2}{\xi_j^2+\<x_j\>^2} \)
\right| , 
\end{align*}
and each derivative of this symbol satisfies the same estimate. Thus
the first term in \eqref{P.4.6} is compact, by the pseudo-differential
calculus. Next we compute: 
\begin{align*}
\left[ i\V_j , {\Op} \( \psi \( \frac{\xi_j^2 -
x_j^2}{\xi_j^2 + \<x_j\>^2} \) \) \right] = {\rm Op} (c_1) +
{\rm Op} (c_2) , 
\end{align*}
with
\begin{align*}
c_1 = \psi ' \( \frac{\xi_j^2 - x_j^2}{\xi_j^2 + \< x_j \>^2}
\) \left\{ {\v}_j , \frac{\xi_j^2 - x_j^2}{\xi_j^2 + \<x_j\>^2}
\right\} = \psi' \( \frac{\xi_j^2 - x_j^2}{\xi_j^2 + \< x_j \>^2}
\) \tilde{c}_1. 
\end{align*}
We have $\tilde{c}_1 \in S \left( \< x \>^{-2} , g_2 \right)$, $c_1 \in S \left( \< x \>^{-1} \< \xi \>^{-1}, g_2 \right)$ and $c_2 \in S \left( \<x\>^{-3} \<\xi\>^{-1} , g_2
\right)$. Thus ${\rm Op}(c_1)$ and ${\rm Op}(c_2)$ are compact by the
pseudo-differential calculus. 
\esp

\noindent {\bf Second step.} General
case.
\esp

Let $J \in C_{\infty} ( \R^{n} )$. We have 
\begin{align*}
s-\lim_{t\rightarrow\infty} e^{itH_{\omega}} J \bigg(
\frac{\ln\<x_1\>}{t} , \ldots ,\frac{\ln\<x_{n}\>}{t} \bigg)
e^{-itH_{\omega}} = \Omega^+ J (\vec{P}_{0}^+) (\Omega^+)^*+ J(0) {\bf
1}^{pp} (H_{\omega}). 
\end{align*}
The existence of $\vec{P}^+$ follows from the existence of
$\vec{P}_{0}^+$. Specializing $\tilde{J} (x_1, \ldots , x_n ) =
J(x_j)$ we obtain furthermore 
\begin{align*}
J(P_{j}^+) = \Omega^+ J ( P_{0,j}^+) (\Omega^+)^*+ J(0) {\bf 1}^{pp}
(H_{\omega}). 
\end{align*}
Then $(ii)$--$(iv)$ follow from this formula and the results on
$P_{0,j}^+$, as in the proof of Theorem~\ref{Th.P.1}. 
\end{proof}

One can ask whether the construction of Theorem~\ref{Th.P.1} works
also in the more general case, and what is the possible link between
the vector $\vec{P}^+$ and $P^+$. The answer is given by the following
theorem.

\begin{theo}   \label{Th.P.14}
Let $\omega_\ell = \displaystyle \max_{1\leq j \leq n_-  }
\omega_j$. There exists 
\begin{align*}
P^+ = s-C_{\infty}-\lim_{t\rightarrow\infty} e^{itH} \frac{ \ln \< x
\>}{t} e^{-itH}, \quad \text{and we have }P^+=P_{\ell}^+\, .
\end{align*}
\end{theo}

\begin{proof}
{\bf First step.} We assume $V\equiv 0$, that is $H = H_{0}$. 

We already know
that $P_{l}^+ = 2 \omega_l$. Thus we only have to show that for $J \in
C_0^{\infty}(\R)$,
\begin{equation}  \label{P.4.9}
s-\lim_{t \rightarrow \infty} e^{itH_{0}} J \( \frac{\ln\<x\>}{t}
\) e^{-itH_{0}} = J(2 \omega_\ell ). 
\end{equation}
We can suppose $\ell=n_-$ and $\omega_1 \le \omega_2 \le \ldots \le
\omega_k < \omega_{k+1} = \ldots = \omega_{n_-}$. Let
$\varepsilon>0$. For $j \in \{ 1, \ldots ,k+1 \}$, we choose
$\widetilde{J}_j \in C_0^{\infty} ([ 2 \omega_j - \varepsilon , 2
\omega_j + \varepsilon ])$, with $\widetilde{J}_j =1$ near
$2\omega_j$. For $j \in \{ n_-+1 , \ldots ,n \}$, we choose
$\widetilde{J}_j \in C_0^{\infty} ([- \varepsilon , \varepsilon ])$
with $\widetilde{J}_j = 1$ near $0$. Then by Propositions~\ref{P.P.6},
\ref{P.P.7} and the separability of the variables, we have 
\begin{align*}
e^{itH_{0}} J \bigg( & \frac{ \ln \< x \>}{t} \bigg) e^{-itH_{0}}   \\
=& e^{itH_{0}} J \bigg( \frac{\ln \< x \>}{t} \bigg)
\widetilde{J}_{k+1} \bigg( \frac{\ln \left\< x_{k+1}, \ldots ,x_{n_-}
\right\>}{t} \bigg) \prod_{ \begin{array}{c}
\scriptstyle{j\in\{1, \ldots ,n\}}  \\
\scriptstyle{j\not \in \{ k+1, \ldots ,n_- \}}
\end{array} }
 \widetilde{J}_j \bigg( \frac{\ln \< x_j \>}{t} \bigg)
e^{-itH_{0}} + R(t), 
\end{align*}
with $s-\lim R(t)=0$. It is clearly sufficient to show
\begin{equation} \label{gl1}
\begin{split}
s-\lim_{t \rightarrow \infty} e^{itH_{0}} \bigg( & J \bigg(
\frac{\ln\<x\>}{t} \bigg) - J \bigg( \frac{\ln \left\< x_{k+1}, \ldots
,x_n \right\>}{t} \bigg) \bigg) \times   \\ 
\times& \widetilde{J}_{k+1} \bigg( \frac{\ln \left\< x_{k+1}, \ldots
,x_{n_-} \right\>}{t} \bigg) \prod_{ \begin{array}{c}
\scriptstyle{j\in\{1, \ldots ,n\}}  \\
\scriptstyle{j\not \in \{ k+1, \ldots ,n_- \}}
\end{array} } \widetilde{J}_j \bigg( \frac{\ln \< x_j \>}{t}
\bigg) e^{-itH_{0}}=0. 
\end{split}
\end{equation}
We have
\begin{align}
\bigg| J \bigg( \frac{\ln \< x \>}{t} \bigg) - J \bigg( \frac{\ln
\left\< x_{k+1}, \ldots ,x_{n_-} \right\>}{t} \bigg) \bigg| \lesssim&
\frac{\ln \< x \>}{t} - \frac{\ln \left\< x_{k+1}, \ldots ,x_{n_-}
\right\>}{t}  \nonumber    \\ 
\lesssim& \frac{1}{t}\ln \Bigg( 1+ \sum_{ \begin{array}{c}
\scriptstyle{j\in\{1, \ldots ,n\}}  \\
\scriptstyle{j\not \in \{ k+1, \ldots ,n_- \}}
\end{array} } \frac{x_j^2}{\left\< x_{k+1}, \ldots ,x_{n_-}
\right\>^2} \Bigg)  \label{gl2} 
\end{align}
For $j\le k$, we have on $\supp\widetilde{J}_j$
\begin{equation}    \label{gl3}
\frac{\ln \< x_j \>^2}{t} \le 4\omega_j +2 \varepsilon \Longrightarrow
x_j^2\le e^{(4\omega_j+2\varepsilon)t}-1 . 
\end{equation}
We have, on $\supp \widetilde{J}_{k+1}$,
\begin{equation}    \label{gl4}
\frac{\ln \left\< x_{k+1}, \ldots ,x_{n_-} \right\>^2}{t} \ge 4
\omega_j-2 \varepsilon \Rightarrow \left\< x_{k+1}, \ldots ,x_{n_-}
\right\>^2 \ge e^{(4 \omega_n -2 \varepsilon )t} . 
\end{equation}
For $j\ge n_-+1$ we have, on $\supp\widetilde{J}_j$,
\begin{equation}     \label{gl3a}
x_j^2 \le e^{2 \varepsilon t}-1.
\end{equation}
Gathering \eqref{gl2}--\eqref{gl3a} together, we obtain
\begin{align*}
\Bigg| \( J \( \frac{\ln \< x \>}{t} \) - J \(
\frac{\ln \left\< x_{k+1}, \ldots ,x_{n_-} \right\>}{t} \) \)
 \widetilde{J}_{k+1} &\( \frac{\ln \left\< x_{k+1}, \ldots
,x_{n_-} \right\>}{t} \)\times\\
&\times \prod_{ \begin{array}{c}
\scriptstyle{j\in\{1, \ldots ,n\}}  \\
\scriptstyle{j\not \in \{ k+1, \ldots ,n_- \}}
\end{array} } \widetilde{J}_j \( \frac{\ln \< x_j \>}{t} 
\) \Bigg|    
\lesssim\, f(t)\, ,
\end{align*}
where $f(t)$ is defined by 
\begin{equation*}
f(t)=  \frac{1}{t}\ln \( 1+ \sum_{j=1}^k
 \frac{e^{(4\omega_j+2\varepsilon)t}-1}{e^{(4\omega_n-2\varepsilon)t}}
 + \sum_{j=n_-+1}^n \frac{e^{2\varepsilon
 t}-1}{e^{(4\omega_n-2\varepsilon)t}} \).
\end{equation*}
If $\varepsilon$ is small enough, $\displaystyle \lim_{t \rightarrow \infty}
f(t)=0$. This yields \eqref{gl1}.

\noindent{\bf Second step.} General case. 
We have 
\begin{align*}
s-\lim_{t \rightarrow \infty} e^{itH} J \bigg( \frac{\ln \< x \>}{t}
\bigg) e^{-itH} =&\Omega^+J(P_0^+)(\Omega^+)^* + J(0) {\bf 1}^{pp}
(H) 
=\Omega^+J(P_{0,\ell}^+)(\Omega^+)^* + J(0) {\bf 1}^{pp} (H)\\
=&\Omega^+J(2\omega_\ell)(\Omega^+)^* + J(0) {\bf 1}^{pp} (H)
=J(2\omega_\ell) {\bf 1}^c (H) + J(0) {\bf 1}^{pp} (H).
\end{align*}
Thus $P^+$ exists, and 
\begin{equation}   \label{P.4.13}
J(P^+) = J(2\omega_\ell) {\bf 1}^c (H) + J(0) {\bf 1}^{pp} (H).
\end{equation}
Let ${\mathcal H}={\mathcal H}_{pp} \oplus {\mathcal H}_c$. Then by
\eqref{P.4.13}, 
$P^+=0 \oplus 2 \omega_\ell = P^+_{\ell}$,
which proves the theorem.
\end{proof}

\bibliographystyle{amsplain}
\bibliography{bchm}

\end{document}